\documentclass[12pt]{amsart}
\usepackage{amsmath,amsthm,amsfonts,amssymb,graphicx,psfrag,dsfont}

\usepackage{hyperref}
\usepackage[alphabetic,initials]{amsrefs}

\theoremstyle{plain}
\newtheorem*{thmu}{Theorem}

\newtheorem{thm}{Theorem}
\newtheorem{prop}{Proposition}[section]

\newtheorem{cor}[prop]{Corollary}
\newtheorem{lemma}[prop]{Lemma}

\newtheorem*{question}{Question}

\theoremstyle{definition}
\newtheorem{example}[prop]{Example}

\newtheorem{defn}[prop]{Definition}
\newtheorem*{notation}{Notation}

\theoremstyle{remark}
\newtheorem{remark}[prop]{Remark}

\newcommand{\Aut}{\operatorname{Aut}}
\newcommand{\aut}{\mathfrak{aut}}

\newcommand{\coker}{\operatorname{coker}}

\newcommand{\cov}{\operatorname{cov}}
\newcommand{\Crit}{\operatorname{Crit}}
\newcommand{\Diff}{\operatorname{Diff}}

\newcommand{\End}{\operatorname{End}}

\newcommand{\Hom}{\operatorname{Hom}}

\renewcommand{\Im}{\operatorname{Im}}
\newcommand{\im}{\operatorname{im}}
\newcommand{\ind}{\operatorname{ind}}

\newcommand{\interior}{\operatorname{int}}

\newcommand{\loc}{\operatorname{loc}}

\newcommand{\MB}{\textsl{MB}}
\newcommand{\muCZ}{\mu_{\operatorname{CZ}}}
\newcommand{\ord}{\operatorname{ord}}

\newcommand{\sing}{\operatorname{sing}}

\newcommand{\virdim}{\operatorname{vir-dim}}

\newcommand{\wind}{\operatorname{wind}}

\newcommand{\SL}{\operatorname{SL}}

\newcommand{\Spp}{\operatorname{Sp}}

\renewcommand{\AA}{{\mathbb A}}
\newcommand{\CC}{{\mathbb C}}
\newcommand{\DD}{{\mathbb D}}

\newcommand{\NN}{{\mathbb N}}

\newcommand{\RR}{{\mathbb R}}
\newcommand{\TT}{{\mathbb T}}
\newcommand{\ZZ}{{\mathbb Z}}
\newcommand{\bB}{{\mathcal B}}
\newcommand{\cC}{{\mathcal C}}

\newcommand{\eE}{{\mathcal E}}

\newcommand{\hH}{{\mathcal H}}
\newcommand{\jJ}{{\mathcal J}}
\newcommand{\lL}{{\mathcal L}}
\newcommand{\mM}{{\mathcal M}}

\newcommand{\oO}{{\mathcal O}}

\newcommand{\tT}{{\mathcal T}}
\newcommand{\uU}{{\mathcal U}}

\newcommand{\1}{\mathds{1}}
\newcommand{\p}{\partial}
\renewcommand{\dbar}{\bar{\partial}}
\newcommand{\Cinftyloc}{C^\infty_{\loc}}

\numberwithin{equation}{section}

\title[Transversality and Orbifolds of Holomorphic Curves]{Automatic
Transversality and Orbifolds of Punctured Holomorphic Curves in Dimension Four}
\author{Chris Wendl}
\address{ETH Z\"urich \\ Departement Mathematik, HG G38.1 \\ 
R\"amistrasse 101 \\
8092 Z\"urich \\ 
Switzerland}
\email{wendl@math.ethz.ch}
\urladdr{http://www.math.ethz.ch/~wendl/}
\thanks{Research partially supported by an NSF Postdoctoral Fellowship
(DMS-0603500) and a DFG grant (CI 45/2-1).}

\begin{document}

\begin{abstract}
We derive a numerical criterion for $J$--holomorphic curves 
in $4$--dimensional symplectic cobordisms to achieve 
transversality without any genericity assumption.  This generalizes results of
Hofer-Lizan-Sikorav \cite{HoferLizanSikorav} and Ivashkovich-Shevchishin
\cite{IvashkovichShevchishin} to allow
punctured curves with boundary that generally need not be
somewhere injective or immersed.  As an application, we combine this with
the intersection theory of punctured holomorphic curves to
prove that certain geometrically natural moduli spaces are
globally smooth orbifolds, consisting generically
of embedded curves, plus unbranched multiple covers that form isolated
orbifold singularities.
\end{abstract}

\maketitle

\tableofcontents

\section{Introduction}

Applications of pseudoholomorphic curves in symplectic
$4$--manifolds and contact $3$--manifolds often depend on the rather special
transversality properties that exist in this low-dimensional setting.
Unlike the general situation, where
the moduli space is smooth only at somewhere 
injective curves and only for generic data, certain moduli
spaces in dimension~$4$ are smooth for all data as long as the right numerical
criteria are satisfied.  For example, suppose $(W,J)$ is any almost
complex $4$--manifold, $(\Sigma,j)$ is a closed Riemann surface of genus~$g$
and $u : (\Sigma,j) \to (W,J)$ is a pseudoholomorphic curve.  The
following result was first mentioned by Gromov \cite{Gromov}, and later
given a complete proof by Hofer-Lizan-Sikorav:

\begin{thmu}[\cite{HoferLizanSikorav}]
If $u$ is embedded and $c_1(u^*TW,J) > 0$, then the moduli space of 
unparametrized pseudoholomorphic
curves near $u$ is a smooth manifold of dimension $2 c_1(u^*TW,J) +
2g - 2$.
\end{thmu}

Observe that the assumptions in the theorem do not require any data to
be generic: rather, the criterion $c_1(u^*TW) > 0$ implies regularity
for uniquely $4$--dimensional reasons that are loosely related to positivity of
intersections.  The dimension of the moduli space is then equal to its
so-called \emph{virtual dimension}, also called the \emph{index} of $u$,
defined as $\ind(u) = 2 c_1(u^*TW) + 2g - 2$.  Thus $c_1(u^*TW) > 0$ is
equivalent to the condition $\ind(u) > 2g - 2$, which leads one to
summarize results of this type with the motto,
\emph{``the moduli space is smooth if the index is sufficiently large.''}
Exactly how large the index needs to be depends on the genus: this is
the reason why almost all applications of such results (including the one
in this paper) principally involve curves of genus zero.

Versions of the theorem above for compact immersed holomorphic
curves with boundary were proved in \cite{HoferLizanSikorav}, 
and similar results for immersed
punctured curves in symplectizations of contact $3$--manifolds also
appeared in \cites{HWZ:props3,Wendl:thesis}.  The reason for dealing with
\emph{immersed} curves in particular was that one could then describe
a neighborhood of $u$ in the moduli space using sections of its normal
bundle and thus reduce the linearization to the so-called
\emph{normal Cauchy-Riemann operator}.  The key fact about
this operator is that its domain is a space of sections on a complex
\emph{line bundle}, thus the zeroes of these sections can be counted and
related to the same topological invariants that appear in the index
formula, giving rise to constraints on the kernel and cokernel.
A generalization for closed holomorphic curves with critical points
was carried out in \cite{IvashkovichShevchishin}, where the normal bundle
was replaced by a \emph{normal sheaf}.

In this paper, we establish a transversality criterion that 
generalizes all of the above results, applying to arbitrary
$J$--holomorphic curves with totally real boundary and cylindrical ends
in $4$--dimensional symplectic cobordisms.

One of the advantages of this approach to transversality is that it applies
to more than just somewhere injective curves: in \S\ref{sec:application}, 
we will
describe a setting in which our criterion, combined with some nontrivial
intersection theory, implies that certain moduli spaces are smooth
orbifolds, which consist mostly of embedded holomorphic curves but also have 
isolated singularities consisting of unbranched
multiple covers over embedded curves.  These moduli spaces
arise quite naturally in a geometric setting: they are the building blocks
of $J$--holomorphic foliations, cf.~\cites{HWZ:foliations,Wendl:OTfol}.

\subsection{The setting}

Let $n \ge 2$.  In all of what follows, $(W,J)$ will denote 
a $2n$--dimensional 
almost complex manifold with noncompact cylindrical ends, which approach
$(2n-1)$--manifolds $M_\pm$ equipped with stable Hamiltonian structures.  
We now recall the precise definitions.

We use the term \emph{stable Hamiltonian structure} to mean the collection of
data that were introduced in \cite{SFTcompactness} as the appropriate setting
for pseudoholomorphic curves in cylindrical manifolds.  Namely, such a
structure $\hH = (\xi,X,\omega,J)$ on a
$(2n-1)$--manifold $M$ consists of the following data:\footnote{The inclusion
of $J$ in the data is somewhat nonstandard but convenient for our purposes.
The data $(\xi,X,\omega)$ are equivalent to the definition of a 
\emph{framed Hamiltonian structure} stated in \cite{EliashbergKimPolterovich},
with the exception that the latter requires $\omega$ to be exact; here it
need only be closed.}
\begin{itemize}
\item $\xi$ is a smooth cooriented hyperplane distribution on $M$
\item $\omega$ is a smooth closed $2$--form on $M$ which restricts to a 
 symplectic structure on the vector bundle $\xi \to M$
\item $X$ is a smooth vector field which is transverse to $\xi$, satisfies
 $\omega(X,\cdot) \equiv 0$, and whose flow preserves $\xi$
\item $J$ is a smooth complex structure on the bundle $\xi \to M$, compatible
 with $\omega$ in the sense that $\omega(\cdot, J\cdot)$ defines a
 bundle metric
\end{itemize}
Note that, as a consequence of these definitions, the flow $\varphi_X^t :
M \to M$ of $X$ 
also preserves the symplectic structure $\omega|_{\xi}$, and the special
$1$--form $\lambda$ associated to $\xi$ and $X$ by the conditions
$$
\lambda(X) \equiv 1,
\qquad
\ker\lambda \equiv \xi,
$$
satisfies $d\lambda(X,\cdot) \equiv 0$.  The \emph{symplectization}
$\RR\times M$ now admits a natural $\RR$--invariant almost complex
structure $\tilde{J}$, defined by the conditions
$$
\tilde{J} \p_a = X,
\qquad
\tilde{J}|_{\xi} = J
$$
where $a$ denotes the coordinate on the $\RR$--factor and
$\p_a \in T(\RR \times M)$ is the corresponding unit vector field.

Recall that a $T$--periodic orbit $x : \RR \to M$ is \emph{nondegenerate}
if the linearized return map $d\varphi_X^T(x(0))|_{\xi_{x(0)}}$ does not
have~$1$ as an eigenvalue.  More generally, a 
\emph{Morse-Bott manifold of $T$--periodic orbits} is a submanifold
$N \subset M$ tangent to $X$ such that $\varphi_X^T|_N$ is the identity,
and for all $p \in N$,
$$
T_p N = \ker\left( d\varphi_X^T(p) - \1 \right).
$$
We will say that an orbit with period $T$ is \emph{Morse-Bott} if it 
is contained in a Morse-Bott manifold of $T$--periodic orbits; note that
this manifold could be a circle, meaning the orbit is nondegenerate.
Moreover, $X$ itself is said to be Morse-Bott (or nondegenerate) if every 
periodic orbit of $X$ is Morse-Bott (or nondegenerate).

We now fix two closed $(2n-1)$--manifolds $M_\pm$ with stable Hamiltonian
structures $\hH_\pm = (\xi_\pm,X_\pm,\omega_\pm,J_\pm)$ and associated
data $\lambda_\pm$ and $\tilde{J}_\pm$, as well as an almost complex
$2n$--manifold $(W,J)$ which decomposes
$$
W = E_- \cup_{M_-} W_0 \cup_{M_+} E_+
$$
so that
\begin{itemize}
\item
$W_0$ is a compact $2n$--manifold with boundary $\p W_0 = M_- \sqcup M_+$
\item
$(E_-,J) \cong ((-\infty,0] \times M_-,\tilde{J}_-)$ and
$(E_+,J) \cong ([0,\infty) \times M_+,\tilde{J}_+)$
\end{itemize}
Fix also a totally real submanifold $L \subset W$.

Near $\p E_\pm \subset E_\pm$, the data $\hH_\pm$ define natural symplectic
forms $\omega_\pm + d(a\lambda_\pm)$ which can be extended (non-uniquely)
over $E_\pm$.  Then given any symplectic form $\omega$ on $W_0$ that attaches
smoothly to $\omega_\pm + d(a\lambda_\pm)$ at $\p W_0$, we denote by
$\jJ_\omega(W,\hH_+,\hH_-)$ the space of 
almost complex structures $J$ on $W$ that are compatible with $\omega$ on
$W_0$
and satisfy the conditions above.\footnote{The symplectic form $\omega$ will 
play almost no role in anything that follows, but becomes important in
applications, e.g.~it yields compactness results as in \cite{SFTcompactness}.}

We will consider pseudoholomorphic (or \emph{$J$--holomorphic}) curves
$$
u : (\dot{\Sigma},j) \to (W, J),
$$
where $\dot{\Sigma} = \Sigma \setminus\Gamma$, $(\Sigma,j)$ is a compact
connected Riemann surface with boundary, 
$\Gamma\subset \interior\Sigma$ is a finite set of interior 
punctures,\footnote{For brevity we're leaving out the case 
of punctures on the boundary, though this can presumably be handled
by similar methods.}
and by definition $u$ satisfies the nonlinear Cauchy-Riemann equation
$T u \circ j = J \circ T u$ and boundary condition 
$u(\p\Sigma) \subset L$.  We also will assume $u$ is 
\emph{asymptotically cylindrical}, which means the following. Partition
the punctures into \emph{positive} and \emph{negative} subsets
$$
\Gamma = \Gamma^+ \cup \Gamma^-,
$$
and at each $z \in \Gamma^\pm$, choose a biholomorphic identification
of a punctured neighborhood of $z$ with the half-cylinder $Z_\pm$, where
$$
Z_+ = [0,\infty) \times S^1
\qquad
\text{ and }
\qquad
Z_- = (-\infty,0] \times S^1.
$$
Then writing $u$ near the puncture in cylindrical coordinates
$(s,t)$, for $|s|$ sufficiently large, it satisfies an asymptotic formula 
of the form
$$
u \circ \varphi(s,t) = \exp_{(Ts,x(Tt))} h(s,t) \in E_\pm.
$$
Here $T > 0$ is a constant, $x : \RR \to M_\pm$ is a $T$--periodic orbit of 
$X_\pm$, the exponential map is defined with respect to any $\RR$--invariant
metric on $\RR\times M_\pm$, $h(s,t) \in \xi_{x(Tt)}$ goes to~$0$ uniformly 
in~$t$ as $s \to \pm\infty$ and $\varphi : Z_\pm \to Z_\pm$ is a smooth
embedding such that
$$
\varphi(s,t) - (s + s_0, t + t_0) \to 0
$$
as $s \to \pm\infty$ for some constants $s_0 \in \RR$, $t_0 \in S^1$.
We will denote by $\gamma_z$ the $T$--periodic orbit 
parametrized by $x$, and call it the \emph{asymptotic orbit} of~$u$ at the 
puncture~$z$.  With this asymptotic behavior in mind, it is convenient
to think of $(\dot{\Sigma},j)$ as a Riemann surface with cylindrical
ends, and we will sometimes refer to neighborhoods of the punctures as 
\emph{ends} of~$\dot{\Sigma}$.  As is well known (cf.~\cites{Hofer:weinstein,
HWZ:props1}), the asymptotically cylindrical holomorphic curves in $(W,J)$ 
are precisely those which satisfy a certain \emph{finite energy} condition, 
though we will not need this fact here.

Denote by $\mM := \mM(J,L)$ the moduli space of equivalence classes of
asymptotically cylindrical $J$--holomorphic curves in $W$ with boundary 
on~$L$; here an equivalence class is defined by the data
$(\Sigma,j,\Gamma,u)$ where $\Gamma$ is considered to be an \emph{ordered}
set, and we define $(\Sigma,j,\Gamma,u) \sim
(\Sigma',j',\Gamma',u')$ if there exists a biholomorphic map
$\varphi : (\Sigma,j) \to (\Sigma',j')$ taking $\Gamma$ to $\Gamma'$ with
the ordering preserved, such that $u = u' \circ \varphi$.  We shall often
abuse notation and write $u \in \mM$ or $(\Sigma,j,\Gamma,u) \in \mM$ 
when we mean $[(\Sigma,j,\Gamma,u)]
\in \mM$.  The moduli space has a natural topology defined by 
$\Cinftyloc$--convergence
on $\dot{\Sigma}$ and uniform convergence up to the ends.  For any
$u \in \mM$, denote by $\mM_u$ the connected component of $\mM$ containing~$u$.

It is often interesting to consider subspaces of $\mM$ defined by imposing
constraints on the asymptotic behavior at some of the punctures.
\begin{defn}
For a given punctured surface $\dot{\Sigma} = \Sigma \setminus
(\Gamma^+ \cup \Gamma^-)$, let $\mathbf{c}$ denote a choice of periodic
orbit $\gamma^{\mathbf{c}}_z$ in $M_\pm$ for some subset of punctures $z \in \Gamma^\pm$.
We call $\mathbf{c}$ a choice of \emph{asymptotic constraints},\footnote{One 
can impose more stringent constraints as well, e.g.~on
the rate at which $u$ converges to its asymptotic orbits; such constraints
are treated in \cites{Wendl:BP1,Wendl:compactnessRinvt}.  Another
possibility is to allow marked points that map to specified points in the
image, perhaps with cusps of prescribed order, as in 
\cites{Barraud:courbes,Francisco:thesis}.
We omit all these possibilities here for the sake of brevity.} 
and refer to each puncture $z$ for which $\mathbf{c}$ specifies an orbit 
$\gamma_z^{\mathbf{c}}$ as a \emph{constrained} puncture.
\end{defn}
For any choice of domain $\dot{\Sigma}$ and asymptotic constraints
$\mathbf{c}$, we can consider the constrained moduli space
$$
\mM^{\mathbf{c}} \subset \mM
$$
consisting of curves $u : \dot{\Sigma} \to W$ that approach the specified 
orbit $\gamma_z^{\mathbf{c}}$ at each of the constrained punctures $z \in \Gamma$,
and arbitrary orbits at the unconstrained punctures.  The constraints define
another partition of $\Gamma$,
$$
\Gamma = \Gamma_C \cup \Gamma_U
$$
into the sets of constrained and unconstrained punctures respectively.
The positive and negative subsets within each of these will be denoted by
$\Gamma_C^\pm$ and $\Gamma_U^\pm$.

If the asymptotic orbits of $u$ are all Morse-Bott, then the
so-called \emph{virtual dimension} of $\mM_u^\mathbf{c}$ is given by the 
\emph{Fredholm index}
\begin{equation}
\label{eqn:index}
\ind(u ; \mathbf{c}) = (n - 3) \chi(\dot{\Sigma}) + 2 c_1^\Phi(u^*TW) +
\mu^\Phi(u ; \mathbf{c})
\end{equation}
where $c_1^\Phi(u^*TW)$ is the relative first Chern number of 
$(u^*TW,J) \to \dot{\Sigma}$ with respect to a suitable choice of 
trivialization $\Phi$ along the ends and boundary, and 
$\mu^\Phi(u ; \mathbf{c})$
is a sum of Conley-Zehnder indices of the asymptotic orbits and
a Maslov index at the boundary with respect to $\Phi$; a precise definition
will be given in \S\ref{subsec:nonlinear}.

As we shall review in more detail in \S\ref{sec:normal}, the nonlinear 
Cauchy-Riemann equation can be expressed as a smooth 
section of a Banach space bundle
$$
\dbar_J : \bB \to \eE : (j,u) \mapsto Tu + J \circ Tu \circ j,
$$
such that a neighborhood of any non-constant 
$(\Sigma,j,\Gamma,u)$ in $\mM^{\mathbf{c}}$ is 
in one-to-one correspondence with $\dbar_J^{-1}(0) / \Aut(\dot{\Sigma},j)$,
where the group $\Aut(\dot{\Sigma},j)$ of biholomorphic maps $(\Sigma,j) \to
(\Sigma,j)$ fixing $\Gamma$ acts on pairs $(j',u') \in
\dbar_J^{-1}(0)$ by
$$
\varphi \cdot (j',u') = (\varphi^*j', u' \circ \varphi).
$$
It is then standard to say that 
$(\Sigma,j,\Gamma,u) \in \mM$ is \emph{regular} if it represents a 
transverse intersection with the zero-section, i.e. the linearization
$$
D\dbar_J(j,u) : T_{(j,u)}\bB \to \eE_{(j,u)}
$$
is surjective.  We will give a precise definition in 
\S\ref{subsec:nonlinear} once the functional analytic setup is in place.  
Observe that if $u : \dot{\Sigma} \to W$ is not constant, then the action of
$\Aut(\dot{\Sigma},j)$ induces a natural inclusion of its Lie algebra
$\aut(\dot{\Sigma},j)$ into $\ker D\dbar_J(j,u)$.
For the sake of completeness, we will present in 
\S\ref{subsec:nonlinear} a proof of the following standard folk theorem:

\setcounter{thm}{-1}

\begin{thm}
\label{thm:orbifoldFolk}
Assume $u : (\dot{\Sigma},j) \to (W,J)$ is a non-constant curve in 
$\mM^\mathbf{c}$ with only Morse-Bott asymptotic orbits.  If $u$ is regular, 
then a neighborhood of $u$ in $\mM^\mathbf{c}$ naturally admits the structure 
of a smooth orbifold of dimension $\ind(u ; \mathbf{c})$, whose isotropy 
group at $u$ is
$$
\Aut(u) := \{ \varphi \in \Aut(\dot{\Sigma},j) |\ 
\text{$u = u\circ \varphi$} \},
$$
and there is a natural isomorphism
$$
T_u\mM^\mathbf{c} = \ker D\dbar_J(j,u) / \aut(\dot{\Sigma},j).
$$
\end{thm}

In particular, regularity implies that $\mM^\mathbf{c}$ is a manifold near 
$u$ if $u$ is somewhere injective, and in general the isotropy group for 
an orbifold singularity
has order bounded by the covering number of~$u$.  Note that in contrast to
the standard theory of $J$--holomorphic curves (cf.~\cite{McDuffSalamon:Jhol}),
we shall in this paper be especially interested in cases where $u$ achieves
regularity despite being multiply covered, so the moduli space is smooth
but may be an orbifold rather than a manifold.

In the case $\dim W = 4$,
another number that turns out to play an important role is the so-called
\emph{normal first Chern number} $c_N(u ; \mathbf{c}) \in \frac{1}{2}\ZZ$, 
which can be defined most simply via the formula
\begin{equation}
\label{eqn:normalChernIndex}
2 c_N(u ; \mathbf{c}) = \ind(u ; \mathbf{c}) - 2 + 2g + \#\Gamma_0(\mathbf{c}) 
+ \#\pi_0(\p\Sigma).
\end{equation}
Here $g$ is the genus of $\Sigma$ and 
$\Gamma_0(\mathbf{c}) \subset \Gamma$ is the subset of punctures for which the
asymptotic orbit has even Conley-Zehnder index (this is the correct definition
if all orbits are nondegenerate; in the Morse-Bott case the definition is
more complicated and may depend on the asymptotic constraints, 
see \S\ref{subsec:nonlinear}).
We will be able to give a better motivated definition in \S\ref{subsec:four} 
using the linear theory in \S\ref{sec:linear}, 
but for now, the significance of $c_N(u ; \mathbf{c})$ can be illustrated 
by considering
the case where $\Sigma$ is closed and $\Gamma = \emptyset$.  Then a
combination of \eqref{eqn:index} and \eqref{eqn:normalChernIndex} yields
the relation
$c_N(u ; \mathbf{c}) = c_1(u^*TW) - \chi(\Sigma)$,
so $c_N(u ; \mathbf{c})$ is the first Chern number of the normal bundle if 
$u$ is immersed.  This is the appropriate philosophical interpretation of
$c_N(u ; \mathbf{c})$ in general, 
as will become obvious from further considerations.

As a final piece of preparation, note that since a non-constant holomorphic
curve $u : \dot{\Sigma} \to W$ is necessarily immersed near the ends,
it can have at most finitely many critical points.  Indeed, as we will
review in \S\ref{subsec:normal},
the bundle $u^*TW \to \dot{\Sigma}$ admits a natural holomorphic structure
such that the section
$$
du \in \Gamma(\Hom_\CC(T\dot{\Sigma},u^*TW))
$$
is holomorphic; its critical points are thus isolated and have positive
order, which we denote by $\ord(du ; z)$ for any $z \in \Crit(u)$.
The quantity
\begin{equation}
\label{eqn:Zdu}
Z(du) := \sum_{z \in du^{-1}(0) \cap \interior{\dot{\Sigma}}} \ord(du ; z)
+ \frac{1}{2} \sum_{z \in du^{-1}(0) \cap \p\Sigma} \ord(du ; z)
\end{equation}
is therefore a finite nonnegative half-integer (or integer if $\p\Sigma =
\emptyset$), and it equals zero if and only if $u$ is immersed.

\subsection{Local and global transversality results}

We now state the main result of this paper.  The following will be a
convenient piece of shorthand notation: if $\p\Sigma \ne \emptyset$,
then for given constants $c \in \RR$ 
and $G \ge 0$, define the nonnegative integer
\begin{equation}
\label{eqn:K}
\begin{split}
K(c,G) = \min \{ k + \ell \ |\ & \text{$k, \ell$ nonnegative integers,}\\
& \text{$k \le G$ and $2k + \ell > 2c$} \}.
\end{split}
\end{equation}
If $\p\Sigma = \emptyset$ we modify this definition slightly by
requiring the integer $\ell$ to be \emph{even}.  Note that in most 
applications known to the author, it will turn out that $c < 0$, so
$K(c,G) = 0$.

\begin{thm}
\label{thm:criterion}
Suppose $\dim W = 4$ and $(\Sigma,j,\Gamma,u) \in \mM^\mathbf{c}$ is a
non-constant curve with only Morse-Bott asymptotic orbits.  If
\begin{equation}
\label{eqn:criterion}
\ind(u ; \mathbf{c}) > c_N(u ; \mathbf{c}) + Z(du),
\end{equation}
then $u$ is regular.  Moreover when this condition is not satisfied, 
we have the following bounds on the dimension of
$\ker D\dbar_J(j,u)$: if $\ind(u ; \mathbf{c}) \le 2Z(du)$ then
\begin{multline*}
2Z(du) \le \dim\left( \ker D\dbar_J(j,u) / \aut(\dot{\Sigma},j) \right) \\
\le 2Z(du) + K(c_N(u ; \mathbf{c}) 
- Z(du),\#\Gamma_0(\mathbf{c}))
\end{multline*}
and if $2Z(du) \le \ind(u ; \mathbf{c})$, then
\begin{multline*}
\ind(u ; \mathbf{c}) \le \dim\left( \ker D\dbar_J(j,u) /
\aut(\dot{\Sigma},j) \right) \\
\le \ind(u ; \mathbf{c}) +
K(c_N(u ; \mathbf{c}) + Z(du) - \ind(u ; \mathbf{c}), \#\Gamma_0(\mathbf{c})).
\end{multline*}
\end{thm}

\begin{remark}
Plugging in the definition of $c_N(u ; \mathbf{c})$ and the index formula, 
the condition \eqref{eqn:criterion} is equivalent to
$$
\ind(u ; \mathbf{c}) > 2g + \#\Gamma_0(\mathbf{c}) + 
\#\pi_0(\p\Sigma) - 2 + 2Z(du),
$$
or
$$
2 c_1^\Phi(u^*TW) + \mu^\Phi(u ; \mathbf{c}) + \#\Gamma_1(\mathbf{c}) > 2Z(du),
$$
where $\Gamma_1(\mathbf{c}) := \Gamma \setminus \Gamma_0(\mathbf{c})$.  
These are direct generalizations of the criteria in
\cites{HoferLizanSikorav,Wendl:thesis,IvashkovichShevchishin}.
\end{remark}

\begin{remark}
An important special case of the dimension bound, which we will use in the
application, appears when $c_N(u ; \mathbf{c}) < Z(du)$: then
$K(c_N(u ; \mathbf{c}) - Z(du),\#\Gamma_0(\mathbf{c})) = 0$, so 
$\dim\ker\left( D\dbar_J(j,u) \right) = 2Z(du)$, its smallest possible value.
\end{remark}

Results of this type have been used previously for a variety of applications,
including disk filling and deformation arguments in contact $3$--manifolds 
\cites{Hofer:weinstein,HWZ:tight3sphere,Wendl:OTfol}, and the symplectic
isotopy problem \cites{Shevchishin:isotopy,Sikorav:gluing}.
In the last section of this paper, we will use our generalization to
prove a somewhat surprising global structure theorem
for certain geometrically natural moduli spaces of holomorphic curves in
$4$--dimensional symplectic cobordisms.  

To motivate this, consider for
a moment the case of a closed holomorphic curve 
$u : \Sigma \to W$ that satisfies the criterion $\ind(u) > c_N(u) +
Z(du)$.  We know then that $\mM$ is smooth in some neighborhood of $u$,
but ideally one would like to know that the entire connected component
$\mM_u$ is smooth.  In general this will not be true, as other curves
in $\mM_u$ may have more critical points and thus fail to satisfy the
criterion.  One favorite way to evade this issue is by assuming that $u$ is
\emph{embedded}: then the adjunction formula (cf.~\cite{McDuffSalamon:Jhol}) 
guarantees that all somewhere injective curves $u' \in \mM_u$ are also
embedded, hence $Z(du') = 0$ and the criterion is satisfied.  
The catch is that unless one imposes additional 
restrictive conditions on the homology class
$[u] \in H_2(W)$, not every curve in $\mM_u$ need be somewhere injective:
a sequence of embedded curves may converge to a branched cover,
which will not always be regular since its branch points are critical.
Thus $\mM_u$ may fail to be globally smooth if it contains branched covers, 
and there is no way to avoid this in general.\footnote{Note that multiply
covered curves also pose a problem in the standard transversality theory,
but for completely different reasons.}

The surprising fact is that if we impose an additional, rather natural 
intersection-theoretic
condition on $u$, then the multiple covers that arise turn out to be
``harmless'': even the multiple covers are regular, and $\mM_u$ is thus 
globally smooth.  The condition in question
arises from the study of $J$--holomorphic foliations: in particular, we
focus on punctured embedded curves $u : \dot{\Sigma} \to W$ that exist in
$1$-- or $2$--dimensional families (with respect to some constraints
$\mathbf{c}$) and have the property of never intersecting
their neighbors, i.e.~these families foliate either an open set or a
hypersurface containing $u(\dot{\Sigma}) \subset W$.  A complete
characterization of such curves is given in \cite{Wendl:BP1} and will be
reviewed in \S\ref{sec:application}; we refer to them as
\emph{stable, nicely embedded} curves.  If $u$ is such a curve, then it
automatically satisfies the criterion of Theorem~\ref{thm:criterion}, thus
the local structure of $\mM_u^\mathbf{c}$ near $u$ is well understood, 
but one still has the global question:

\begin{question}
Can a sequence of stable, nicely embedded curves converge to
a multiple cover?
\end{question}
If the answer is no, then $\mM_u^\mathbf{c}$ is a smooth manifold, and
we'll show that this is indeed the case whenever $W$ is an $\RR$--invariant 
symplectization (with generic $J$) 
or a closed symplectic manifold.  In general, it turns out that
multiple covers \emph{can} appear, but only if they are \emph{immersed},
in which case the regularity criterion is still satisfied.  The proof
of this fact will make use of our transversality arguments for non-immersed
curves, establishing in effect that any component of $\mM^\mathbf{c}$
containing such a non-immersed multiple cover can consist only of multiple
covers.  The result is:

\begin{thm}
\label{thm:orbifold}
For generic $J$, if $u \in \mM^\mathbf{c}$ is a stable, nicely embedded curve,
then every curve in $\mM_u^\mathbf{c}$ is regular: in particular
$\mM_u^\mathbf{c}$ naturally admits the structure of a smooth 
orbifold of dimension $ind(u ; \mathbf{c}) \in \{1,2\}$, with
only isolated singularities.
Moreover, all curves in $\mM_u^\mathbf{c}$ are embedded
except for a discrete subset, consisting of unbranched multiple
covers over stable, nicely embedded index~$0$ curves, and the images of 
any two curves in $\mM_u^\mathbf{c}$ are either identical or disjoint.
\end{thm}

This will follow from a more general result (Theorem~\ref{thm:orbifoldStrong})
proved in \S\ref{sec:application}, which applies also to parametrized
moduli spaces under a generic homotopy of almost complex structures.
As a simple corollary, we observe the two aforementioned 
cases where the answer to the question posed above is no:

\begin{cor}
\label{cor:noMultiples}
For the curve $u : \dot{\Sigma} \to W$ in Theorem~\ref{thm:orbifold},
suppose that either
\begin{itemize}
\item $\dot{\Sigma}$ is a closed Riemann surface (without punctures), or
\item $(W,J) = (\RR\times M,\tilde{J})$ is the symplectization of a
$3$--manifold with stable Hamiltonian structure $\hH = (X,\xi,\omega,J)$,
where $J$ is generic.
\end{itemize}
Then every curve in $\mM_u^\mathbf{c}$ is embedded, thus
$\mM_u^\mathbf{c}$ is a manifold.
\end{cor}
\begin{proof}
For the $\RR$--invariant symplectization $(\RR\times M,\tilde{J})$, a 
multiple cover $u = v \circ \varphi$ would require a somewhere injective
curve $v$ of index~$0$, which doesn't exist if $J$ is generic.  The
reasoning in the closed case is different: it depends on the fact that,
as we'll show in \S\ref{subsec:immersed}, stable nicely embedded curves 
always have genus zero.  Then $\varphi$ must be a holomorphic map
$S^2 \to S^2$ with no branch points, contradicting the Riemann-Hurwitz
formula.
\end{proof}

Unbranched multiple covers can and do appear in general
if $\dot{\Sigma}$ has punctures and $(W,J)$ is a non-cylindrical
manifold, e.g.~a nontrivial symplectic cobordism.  We will show an
example at the end of \S\ref{sec:application}, where the resulting collection
of curves actually foliates~$W$.

The phenomenon illustrated by Theorem~\ref{thm:orbifold} contrasts with
the more general study of holomorphic curves, e.g.~in Symplectic Field
Theory, where transversality can only be achieved in general by
abstract perturbations.  Such perturbations usually destroy many of the nice 
geometric properties of holomorphic curves---such as positivity of 
intersections---but the philosophy here is that for curves that are especially
nice in some geometric sense, precisely these nice properties make abstract 
perturbations
unnecessary.  In particular, the theorem is part of a larger program outlined
in \cite{Wendl:compactnessRinvt}, to prove that the compactified moduli
spaces of curves that can occur in foliations always have a nice global
structure: in principle, after proving a suitable compactness theorem for
this ``nice'' class of curves, transversality should always follow 
``for free''.  Such results are necessary tools
in the general theory of $J$--holomorphic foliations, as one would like
to prove that these foliations can always be 
carried through under various types of homotopies and stretching arguments.
The situation is already well understood in the $\RR$--invariant case due to
\cite{Wendl:compactnessRinvt}, and Theorem~\ref{thm:orbifold} may
be seen as a partial result in the direction of generalizing that 
compactness theorem to symplectic cobordisms.  
(See Example~\ref{ex:foliation} and
Remark~\ref{remark:bigCompactness} for an idea
of what such a generalization might look like.)

\subsection{Outline of the proofs}
\label{subsec:outline}

The technical backbone of Theorem~\ref{thm:criterion} is the analysis of
the \emph{normal Cauchy-Riemann operator} $\mathbf{D}^N_u$ associated to
any holomorphic curve $u : \dot{\Sigma} \to W$.  As we will recall,
this is well defined 
even if $u$ has critical points, because there always exists a splitting
$$
u^*TW = T_u \oplus N_u
$$
such that $(T_u)_z$ is the image of $T_z\dot{\Sigma}$ under $du$ at all
regular points~$z$.  The domain of $\mathbf{D}^N_u$ is then a space of
sections of $N_u$, a complex line bundle.  We describe the required 
linear theory of such operators in \S\ref{sec:linear}, giving criteria
that guarantee surjectivity of $\mathbf{D}^N_u$ as well as bounds on the
dimension of its kernel.  

The next step is then to relate the operator $\mathbf{D}^N_u$ to the 
nonlinear problem.  In the immersed case, the traditional approach
(cf.~\cites{HoferLizanSikorav,HWZ:props3})
is to set up the nonlinear problem to detect $J$--invariant maps
that can be expressed as sections of the normal bundle of a given
solution~$u$, in which case the linearization is equivalent to 
$\mathbf{D}^N_u$.  This is
no longer possible when $u$ has critical points; Ivashkovich and Shevchishin 
in \cite{IvashkovichShevchishin} dealt with this difficulty by replacing the 
normal bundle with a \emph{normal sheaf} and proving that $u$ is regular
if and only if $\mathbf{D}_u^N$ is surjective.  Our approach takes some
inspiration from theirs but is
less algebraic and more analytical in flavor, as we avoid any
reference to sheaves and exact sequences in favor of Banach space splittings 
and Fredholm operators.  Unlike \cites{HoferLizanSikorav,HWZ:props3},
we treat the nonlinear problem in the way that is standard for
arbitrary dimensions, as a section
$$
\dbar_J : \tT \times \bB \to \eE
$$
of a suitable Banach space bundle, where $\bB$ is a (globally defined)
Banach manifold of maps $\dot{\Sigma} \to W$ (including reparametrizations) 
and $\tT$ is a (locally defined)
finite dimensional space of complex structures parametrizing on open
subset in the
Teichm\"uller space of $\dot{\Sigma}$.  We will use the splitting
$u^*TW = T_u \oplus N_u$ and some properties of the standard Cauchy-Riemann
operator on $\Gamma(T\dot{\Sigma})$ to give a precise relation 
(Theorem~\ref{thm:normalCR}) between the kernels and images
of $D\dbar_J(j,u)$ and
$\mathbf{D}_u^N$ in arbitrary dimensions.  A consequence is the fact that
each of these operators is surjective if and only if the other is.

As for the proof of Theorem~\ref{thm:orbifold}: assume $u_n$ is a 
sequence of stable, nicely embedded curves converging to a multiple cover
$u = v \circ \varphi$, where $v$ is somewhere injective.
We observe first that the embedded curves $u_n$ necessarily satisfy the 
criterion of Theorem~\ref{thm:criterion}, so this will remain true for the 
limit $u$ unless it acquires critical points.  The main task then is to
show that $u$ is immersed, and the kernel bounds in Theorem~\ref{thm:criterion}
for non-immersed curves turn out to be a useful tool in proving this.
The first step is to show that the underlying simple curve $v$ is embedded
and has index~$0$: this follows by a careful application of the intersection 
theory of punctured holomorphic curves, which we review at the beginning of
\S\ref{sec:application}. 
Note that this is the only point in the argument at which we assume $J$
to be generic: it's necessary to obtain a lower bound on the index of $v$ and
thus on its related intersection invariants,
but it will not be required in proving transversality for~$u$.
With this established,
critical points of $u$ arise only from branch points of the cover $\varphi$,
hence $Z(du) = Z(d\varphi)$, i.e.~the ramification number of~$\varphi$.
Now the dimension bound in Theorem~\ref{thm:criterion} turns out to
imply that a neighborhood of $u$ in $\mM_u^\mathbf{c}$ 
``lives inside a space of dimension at most $2 Z(du)$''; we will make this
statement precise later using the implicit function theorem.
But if $Z(du) > 0$, then the space of holomorphic branched covers homotopic
to $\varphi$ is nontrivial and has precisely this dimension, which yields
a $2Z(du)$--dimensional 
smooth submanifold of $\mM_u^\mathbf{c}$ containing~$u$.
It follows that this describes a neighborhood of $u$ in $\mM_u^\mathbf{c}$, 
so any sequence of curves converging to $u$ must then have the form 
$u_n = v \circ \varphi_n$, i.e.~they are all multiple covers with the same
image, and this is a contradiction.

\subsubsection*{Acknowledgments}

Many thanks to Denis Auroux, Kai Cieliebak, Oliver Fabert, Helmut Hofer, 
Sam Lisi, Klaus Mohnke, Sewa Shevchishin and Richard Siefring for useful 
conversations.

\section{Cauchy-Riemann type operators on bundles}
\label{sec:linear}

\subsection{Generalities}
\label{subsec:general}

Let $(\Sigma,j)$ 
be a compact Riemann surface with genus~$g$, $m \ge 0$~boundary
components, and a finite set of positive/negative interior punctures
$\Gamma = \Gamma^+ \cup \Gamma^- \subset \interior\Sigma$, with the
corresponding punctured surface denoted by $\dot{\Sigma} = \Sigma \setminus
\Gamma$.  Regarding $\dot{\Sigma}$ as a surface with cylindrical ends
$\{\uU_z\}_{z \in \Gamma^\pm}$
biholomorphic to the half-cylinders $Z_\pm$, it admits a natural
compactification $\overline{\Sigma}$ obtained by replacing
$[0,\infty) \times S^1$ by $[0,\infty]\times S^1$ and
$(-\infty,0] \times S^1$ by $[-\infty,0] \times S^1$.  The compactified
space is naturally a topological $2$--manifold with boundary
$$
\p\overline{\Sigma} = \p\Sigma \sqcup \bigcup_{z \in \Gamma} \delta_z,
$$
where for each $z \in \Gamma^\pm$, $\delta_z \cong \{\pm\infty\} \times S^1$
denotes the corresponding ``circle at infinity''.  Note that in making this
definition we've chosen cylindrical coordinates $(s,t) \in Z_\pm$ over
each end $\{\uU_z\}_{z \in \Gamma^\pm}$, and we will continue to use these
coordinates whenever convenient.  The definitions of $\overline{\Sigma}$ and
$\delta_z$ do not depend on this choice, and in fact the resulting
identification of each $\delta_z$ with $S^1 = \RR / \ZZ$ is unique up to
a constant shift.

Let $E \to \overline{\Sigma}$ be a complex vector bundle of rank~$n$
whose restriction
to $\dot{\Sigma}$ and each of the circles $\delta_z$ has a smooth structure.
Assume moreover that $E$ is given a Hermitian structure over each end $\uU_z$.
By an \emph{admissible trivialization} of $E$ near $z \in \Gamma^\pm$, 
we mean a smooth unitary bundle isomorphism $\Phi : E|_{\uU_z} \to 
Z_\pm \times \RR^{2n}$ (where $\RR^{2n}$ is identified with $\CC^n$), which
covers the coordinate map $\uU_z \to Z_\pm$ and
extends continuously to a smooth unitary trivialization 
$E|_{\delta_z} \to S^1 \times \RR^{2n}$.  An
\emph{asymptotic operator} at $z \in \Gamma$ is then
a bounded real linear operator
$$
\mathbf{A}_z : H^1(E|_{\delta_z}) \to L^2(E_{\delta_z})
$$
whose expression with respect to any admissible trivialization takes the
form
$$
H^1(S^1,\RR^{2n}) \to L^2(S^1,\RR^{2n}) : 
\eta \mapsto -J_0 \dot{\eta} - S \eta;
$$
here $J_0 = i$ is the standard complex structure on $\RR^{2n} = \CC^n$ and
$S = S(t)$ is any smooth loop of symmetric $2n$--by--$2n$ 
matrices.  This defines
an unbounded self-adjoint operator on the complexification of
$L^2(E|_{\delta_z})$.  We say that $\mathbf{A}_z$ is \emph{nondegenerate}
if its spectrum $\sigma(\mathbf{A}_z)$ does not contain~$0$.

Define the standard $\dbar$--operator for smooth functions on $\dot{\Sigma}$ by
$$
\dbar : C^\infty(\dot{\Sigma},\CC) \to \Omega^{0,1}(\dot{\Sigma}) : f
\mapsto df + i\ df \circ j.
$$
For any two complex vector bundles $E$ and $E'$ over the same base, we denote 
by $\Hom_\CC(E,E')$ and $\overline{\Hom}_\CC(E,E')$ the corresponding bundles 
of complex linear and antilinear maps $E \to E'$ respectively.  There are also
the corresponding endomorphism bundles $\End_\CC(E) := \Hom_\CC(E,E)$ and
$\overline{\End}_\CC(E) := \overline{\Hom}_\CC(E,E)$.

\begin{defn}
\label{defn:CR}
A (smooth, real linear) \emph{Cauchy-Riemann type operator} on $E$ is a
first-order linear differential operator
$$
\mathbf{D} : \Gamma\left(E|_{\dot{\Sigma}}\right) \to 
\Gamma\left(\overline{\Hom}_\CC(T\dot{\Sigma},E|_{\dot{\Sigma}})\right)
$$
such that for every smooth section $v : \dot{\Sigma} \to E$ and smooth
function $f : \dot{\Sigma} \to \RR$,
$$
\mathbf{D} (fv) = (\dbar f) v + f (\mathbf{D} v).
$$
Given an asymptotic operator $\mathbf{A}_z$ at $z \in \Gamma^\pm$, we
will say that $\mathbf{D}$ is \emph{asymptotic to $\mathbf{A}_z$} if
its expression in an admissible trivialization $\Phi$ near $z$ takes the form
$$
(\mathbf{D} v)(s,t) = \p_s v(s,t) + J_0 \p_t v(s,t) + S(s,t) v(s,t),
$$
where $S(s,t)$ is a smooth family of real-linear transformations on
$\RR^{2n}$ which converges uniformly as $s \to \pm\infty$ to a smooth
loop of symmetric matrices $S(t)$, such that
$$
-J_0 \frac{d}{dt} - S(t)
$$
is the coordinate expression for $\mathbf{A}_z$ with respect to~$\Phi$.
\end{defn}

Define the Banach space $W^{k,p}(E)$ to consist of sections
$v : \dot{\Sigma} \to E$ of class $W^{k,p}_{\text{loc}}$ such that in
any choice of admissible trivialization near each puncture $z \in \Gamma^\pm$,
the corresponding map
$Z_\pm \to \RR^{2n}$ is of class $W^{k,p}$.  If $\ell \subset E|_{\p\Sigma}$
is a smooth totally real submanifold, define the subspace
$$
W^{k,p}_\ell(E) = \{ v \in W^{k,p}(E)\ |\ v(\p\Sigma) \subset \ell \}.
$$
Observe that 
$\overline{\Hom}_\CC(T\dot{\Sigma},E)$ also admits a natural extension
over $\overline{\Sigma}$, and the combination of the coordinates $(s,t)$ 
with the trivialization $\Phi$ near $z \in \Gamma$ also gives rise to a 
trivialization of $\overline{\Hom}_\CC(T\dot{\Sigma},E)$.  Using this
we can define the Banach spaces 
$W^{k,p}(\overline{\Hom}_\CC(T\dot{\Sigma},E))$.  We will
generally write $W^{0,p}$ as $L^p$.

Now fix a smooth totally real subbundle $\ell \subset E|_{\p\Sigma}$, and
asymptotic operators $\mathbf{A}_z$ for each $z \in \Gamma$,
denoting the collection of all these operators by $\mathbf{A}_\Gamma$.
Let $\mathbf{D}$ be a
Cauchy-Riemann type operator that is asymptotic to
$\mathbf{A}_z$ for each $z \in \Gamma$.  We will then be interested in the
bounded linear operator
$$
\mathbf{D} : W^{1,p}_{\ell}(E) \to L^p(\overline{\Hom}_\CC(T\dot{\Sigma},E)).
$$
This is a Fredholm operator if all the $\mathbf{A}_z$ are nondegenerate,
and its index is determined by a variety of 
topological quantities which we shall recall next.

Fix a set of admissible trivializations near each puncture $z \in \Gamma$
as well as smooth complex trivializations of $E|_{\p\Sigma}$, denoting
the collection of all these choices by~$\Phi$.  One can then define
the \emph{relative first Chern number} $c_1^\Phi(E) \in \ZZ$.  If
$E$ is a line bundle, then $c_1^\Phi(E)$ is defined simply by counting
zeroes of a generic smooth section $\dot{\Sigma} \to E$ that extends
continuously over $\overline{\Sigma}$ and is a nonzero constant with 
respect to $\Phi$ on $\p\overline{\Sigma}$.  For
higher rank bundles, $c_1^\Phi(E)$ can be defined axiomatically via the
direct sum property and the assumption that it matches the ordinary
first Chern number if $\dot{\Sigma}$ is closed.

For each connected component $C \subset \p\Sigma$, the totally
real subbundle $\ell|_{C} \subset E|_C$ has a \emph{Maslov index}
$\mu^\Phi\left(E|_{C},\ell|_{C}\right)$, 
and we shall denote the sum of these by $\mu^\Phi(E,\ell)$.

Finally for each puncture $z \in \Gamma^\pm$, the asymptotic operator
$\mathbf{A}_z$, expressed as $-J_0 \p_t - S(t)$ with respect to~$\Phi$,
gives rise to a linear Hamiltonian flow in $\RR^{2n}$ via the equation
$$
\dot{\eta}(t) = J_0 S(t) \eta(t).
$$
If $\mathbf{A}_z$ is nondegenerate, then the 
resulting path of symplectic matrices $\Psi(t) \in \Spp(n)$ ends at
a matrix $\Psi(1)$ which does not have~$1$ as an eigenvalue, 
so it has a well defined
\emph{Conley-Zehnder index} which we denote by $\muCZ^\Phi(\mathbf{A}_z)$.
All of this together allows us to define the \emph{total Maslov index}
$$
\mu^\Phi(E,\ell,\mathbf{A}_\Gamma) :=
\mu^\Phi(E,\ell) + \sum_{z \in \Gamma^+} \muCZ^\Phi(\mathbf{A}_z) -
\sum_{z \in \Gamma^-} \muCZ^\Phi(\mathbf{A}_z).
$$
The Fredholm index of $\mathbf{D}$ is then given by the following
generalization of the Riemann-Roch formula:
\begin{equation}
\label{eqn:RiemannRoch}
\ind(\mathbf{D}) = n \chi(\dot{\Sigma}) + 2 c_1^\Phi(E) +
\mu^\Phi(E,\ell,\mathbf{A}_\Gamma).
\end{equation}
This follows from the formula for the case $\p\Sigma = \emptyset$
proved in \cite{Schwarz}, together with a gluing/doubling argument;
cf.~\cite{Wendl:thesis}.  Note that all dependence on $\Phi$ in the
right hand side of \eqref{eqn:RiemannRoch} cancels out.

Let us briefly review the useful generalization of the above that 
arises by considering Banach spaces
with exponential weights.  Pick numbers $\boldsymbol{\delta}_z \in \RR$ for
each $z \in \Gamma$ and denote the collection of these by
$\boldsymbol{\delta}_\Gamma = \{ \boldsymbol{\delta}_z \}_{z\in\Gamma}$.  Then we define
$$
W^{k,p,\boldsymbol{\delta}_\Gamma}(E)
$$
to be the space of $W^{k,p}_{\text{loc}}$ sections $v : \dot{\Sigma} \to E$
such that in an admissible trivialization near each $z \in \Gamma^\pm$,
the function $Z_\pm \to \RR^{2n} : (s,t) \mapsto e^{\pm\boldsymbol{\delta}_z s} v(s,t)$
is of class~$W^{k,p}$.  This imposes an exponential decay condition at each
puncture where $\boldsymbol{\delta}_z > 0$, or a bound on exponential growth if $\boldsymbol{\delta}_z
< 0$.  There are now obvious definitions for the spaces
$W^{1,p,\boldsymbol{\delta}_\Gamma}_\ell(E)$ and 
$L^{p,\boldsymbol{\delta}_\Gamma}(\overline{\Hom}_\CC(T\dot{\Sigma},E))$, so that
the Cauchy-Riemann type operator $\mathbf{D}$ defines a bounded linear map
$$
\mathbf{D} : W^{1,p,\boldsymbol{\delta}_\Gamma}_\ell(E) \to 
L^{p,\boldsymbol{\delta}_\Gamma}(\overline{\Hom}_\CC(T\dot{\Sigma},E)).
$$
It is simple to show (cf.~\cites{HWZ:props3,Wendl:BP1}) that 
this map is conjugate to another Cauchy-Riemann operator
$\mathbf{D}_{\boldsymbol{\delta}_\Gamma} : W^{1,p}_\ell(E) \to
L^p(\overline{\Hom}_\CC(T\dot{\Sigma},E))$,
which is asymptotic at $z \in \Gamma^\pm$ to 
$\mathbf{A}_z \pm \boldsymbol{\delta}_z$; denote the latter collection of operators
by $\mathbf{A}_\Gamma \pm \boldsymbol{\delta}_\Gamma$.
The operator on the weighted space is thus
Fredholm if and only if $\mp\boldsymbol{\delta}_z \not\in \sigma(\mathbf{A}_z)$
for all $z \in \Gamma^\pm$, and
its index can then be read off again from \eqref{eqn:RiemannRoch},
but with $\mathbf{A}_\Gamma \pm \boldsymbol{\delta}_\Gamma$ replacing
$\mathbf{A}_\Gamma$.  Note in particular that if all $\mathbf{A}_z$ are 
nondegenerate and all $\boldsymbol{\delta}_z$ are sufficiently close to~$0$, then the 
weighting does not change the index of~$\mathbf{D}$.

\subsection{The line bundle case}
\label{subsec:linebundle}

For the rest of this section we assume $n=1$, so each asymptotic operator
is equivalent to an unbounded self-adjoint operator on $L^2(S^1,\RR^2)$
of the form $\mathbf{A} = - J_0 \frac{d}{dt} - S(t)$, whose eigenfunctions
can be assigned winding numbers.
For $\lambda \in \sigma(\mathbf{A})$ define $w(\lambda) \in \ZZ$ to be 
the winding number of any nontrivial section in the $\lambda$--eigenspace of 
$\mathbf{A}$; this number depends only on $\lambda$, by a result in
\cite{HWZ:props2}.  Moreover, it is shown in the same paper that
$w(\lambda)$ is an increasing function of $\lambda$ which takes every
integer value exactly twice (counting multiplicity).  We define
\begin{equation}
\label{eqn:defAlpha}
\begin{split}
\alpha_-(\mathbf{A}) &= \max \{ w(\lambda) \ |\ 
\text{$\lambda \in \sigma(\mathbf{A})$, $\lambda < 0$} \}, \\
\alpha_+(\mathbf{A})  &= \min \{ w(\lambda) \ |\ 
\text{$\lambda \in \sigma(\mathbf{A})$, $\lambda > 0$} \}, \\
p(\mathbf{A})      &= \alpha_+(\mathbf{A}) - \alpha_-(\mathbf{A}),
\end{split}
\end{equation}
so if $\mathbf{A}$ is nondegenerate, $p(\mathbf{A}) \in \{0,1\}$.
By another result in \cite{HWZ:props2}, these winding numbers
are related to the Conley-Zehnder index by
\begin{equation}
\label{eqn:CZwinding}
\muCZ(\mathbf{A}) = 2\alpha_-(\mathbf{A}) + p(\mathbf{A}) =
2\alpha_+(\mathbf{A}) - p(\mathbf{A}).
\end{equation}
This entire discussion applies also to the operators $\mathbf{A}_z$ once
trivializations $\Phi$ are specified; we thus denote these winding numbers
by $\alpha_\pm^\Phi(\mathbf{A}_z)$, and observe
that $p(\mathbf{A}_z) \in \{0,1\}$ does not depend on $\Phi$.  The latter is
the \emph{parity} of the puncture $z \in \Gamma$, defining a partition of
$\Gamma$ into sets of \emph{even} and \emph{odd} punctures, 
denoted $\Gamma_0$ and $\Gamma_1$ respectively.

Define the $\frac{1}{2}\ZZ$--valued 
\emph{adjusted first Chern number} of $(E,\ell,\mathbf{A}_\Gamma)$ by
\begin{equation}
\label{eqn:c1adjusted}
c_1(E,\ell,\mathbf{A}_\Gamma) = c_1^\Phi(E) + 
\frac{1}{2} \mu^\Phi(E,\ell) +
\sum_{z \in \Gamma^+} \alpha_-^\Phi(\mathbf{A}_z) -
\sum_{z \in \Gamma^-} \alpha_+^\Phi(\mathbf{A}_z),
\end{equation}
and observe that this does not depend on $\Phi$.  Using \eqref{eqn:CZwinding}
and the index formula, it is easy to show that
\begin{equation}
\label{eqn:2c1}
2 c_1(E,\ell,\mathbf{A}_\Gamma) = \ind(\mathbf{D}) - 2 + 2g + \#\Gamma_0 + m.
\end{equation}
Note that $c_1(E,\ell,\mathbf{A}_\Gamma)$ is necessarily an integer if
$\p\Sigma = \emptyset$.

The adjusted first Chern number has the following interpretation
which justifies its name.  If $v \in \ker\mathbf{D}$ is a nontrivial
section, then the equation $\mathbf{D}v = 0$ together with the similarity 
principle implies that $v$ has only isolated zeroes, all of positive
order.  Moreover, by arguments in \cites{HWZ:props1,Siefring:asymptotics},
$v$ satisfies an asymptotic formula of the form
\begin{equation}
\label{eqn:linearAsymptotic}
v(s,t) = e^{\lambda s}(e_\lambda(t) + r(s,t))
\end{equation}
in admissible trivializations near each puncture $z \in \Gamma^\pm$, where
$\lambda \in \sigma(\mathbf{A}_z)$ satisfies $\pm\lambda < 0$, $e_\lambda \in
\Gamma(E|_{\delta_z})$ is a section in the corresponding 
eigenspace and the remainder
$r(s,t)$ goes to zero as $s \to \pm\infty$.  It follows that $v(s,t)$ has
only finitely many zeroes, and near $z \in \Gamma^\pm$ it has
a well defined \emph{asymptotic winding number} $\wind_z^\Phi(v) \in \ZZ$,
which is bounded from above by $\alpha_-^\Phi(\mathbf{A}_z)$ if $z \in \Gamma^+$,
or from below by $\alpha_+^\Phi(\mathbf{A}_z)$ if $z \in \Gamma^-$.  
We use this to define the \emph{asymptotic vanishing} of $v$:
$$
Z_\infty(v) = \sum_{z \in \Gamma^+} \left[ \alpha_-^\Phi(\mathbf{A}_z) - 
\wind^\Phi_z(v)\right] + \sum_{z \in \Gamma^-} 
\left[ \wind^\Phi_z(v) - \alpha_+^\Phi(\mathbf{A}_z)\right].
$$
Define also the $\frac{1}{2}\ZZ$--valued count of zeroes,
$$
Z(v) = \sum_{z \in v^{-1}(0) \cap \interior{\dot{\Sigma}}} \ord(v ; z)
+ \frac{1}{2} \sum_{z \in v^{-1}(0) \cap \p\Sigma} \ord(v ; z),
$$
where the order of a zero on the boundary is defined by a doubling argument
described in the appendix.  
Now a simple computation using these definitions and
Prop.~\ref{prop:halfc1} yields the relation
\begin{equation}
\label{eqn:c1counts}
Z(v) + Z_\infty(v) = c_1(E,\ell,\mathbf{A}_\Gamma).
\end{equation}
Observe that both terms on the left hand side are manifestly
nonnegative.

The next result is the main objective of this section.  Recall from
\eqref{eqn:K} the nonnegative integer $K(c,G)$.

\begin{prop}
\label{prop:transversality}
\ \\
\begin{enumerate}
\item
In the case $\ind(\mathbf{D}) \le 0$, $\mathbf{D}$ is injective
if $c_1(E,\ell,\mathbf{A}_\Gamma) < 0$, and otherwise
$$
\dim\ker\mathbf{D} \le 
K\left(c_1(E,\ell,\mathbf{A}_\Gamma),\#\Gamma_0\right).
$$
\item
In the case $\ind(\mathbf{D}) \ge 0$, $\mathbf{D}$ is surjective if
$\ind(\mathbf{D}) > c_1(E,\ell,\mathbf{A}_\Gamma)$, and otherwise
$$
\ind(\mathbf{D}) \le \dim\ker\mathbf{D} \le \ind(\mathbf{D}) +
K\left(c_1(E,\ell,\mathbf{A}_\Gamma) - \ind(\mathbf{D}), \#\Gamma_0\right).
$$
\end{enumerate}
\end{prop}
\begin{proof}
The argument rests crucially on \eqref{eqn:c1counts}, together with
the observation that if $z \in \Gamma_0$, then the space of 
eigenfunctions with negative eigenvalue
and winding $\alpha_-^\Phi(\mathbf{A}_z)$ is $1$--dimensional, as is the
space with positive eigenvalue and winding $\alpha_+^\Phi(\mathbf{A}_z)$.

We prove the result first for the case $\ind(\mathbf{D}) \le 0$.  If 
$c_1(E,\ell,\mathbf{A}_\Gamma) < 0$, then $\mathbf{D}$ is clearly injective, 
else \eqref{eqn:c1counts} would force any nontrivial section 
$v \in \ker\mathbf{D}$ to have either $Z(v) < 0$ or $Z_\infty(v) < 0$.
To establish the dimension bound for $\ker\mathbf{D}$ when
$c_1(E,\ell,\mathbf{A}_\Gamma) \ge 0$, choose 
any nonnegative integers $k$ and $n$ such that $k \le \#\Gamma_0$
and $2k + n > 2 c_1(E,\ell,\mathbf{A}_\Gamma)$.  In this situation we can
construct an injective homomorphism from $\ker\mathbf{D}$ into a real
vector space of dimension $n + k$.  Indeed, if $\p\Sigma \ne \emptyset$,
pick $n$ distinct points $\zeta_1,\ldots,\zeta_n \in \p\Sigma$, and choose also
$k$ distinct even punctures $z_1,\ldots,z_k \in \Gamma_0$.  For each of the
$z_j$, the correspondence $v \mapsto e_\lambda$ coming from the asymptotic 
formula \eqref{eqn:linearAsymptotic} defines
a linear map from $\ker\mathbf{D}$ into the $1$--dimensional vector space 
$V_{z_j}$ consisting of eigenfunctions of $\mathbf{A}_{z_j}$ with winding 
equal to $\alpha_\mp^\Phi(\mathbf{A}_{z_j})$.  We can define this map so that
it takes the value $0 \in V_{z_j}$ if and only if the eigenfunction in
\eqref{eqn:linearAsymptotic} has a different winding number.  Using these
maps and the evaluation of $v$ at the points $\zeta_j \in \p\Sigma$, we
obtain a homomorphism
$$
\Psi : \ker\mathbf{D} \to \ell_{\zeta_1} \oplus \ldots \oplus \ell_{\zeta_n}
\oplus V_{z_1} \oplus \ldots \oplus V_{z_k}.
$$
The claim is that $\Psi$ is injective, and thus $\dim\ker\mathbf{D}
\le n + k$.  Indeed, suppose $v \in \ker\mathbf{D}$ is a nontrivial
section with $\Psi(v) = 0$.  Then the asymptotic winding of $v$ differs
from $\alpha_\mp^\Phi(\mathbf{A}_{z_j})$ at each of the punctures $z_j$,
implying $Z_\infty(v) \ge k$.  Similarly, $v$ has boundary zeroes at
$\zeta_1,\ldots,\zeta_n$, contributing at least $n/2$ to $Z(v)$, hence
$$
c_1(E,\ell,\mathbf{A}_\Gamma) = Z(v) + Z_\infty(v) \ge \frac{n}{2} + k,
$$
which contradicts our assumptions on $n$ and~$k$.

A minor modification to this argument is needed if $\p\Sigma = \emptyset$.
We must now assume $n$ is even, and choose distinct interior points
$\zeta_1,\ldots,\zeta_{n/2} \in \dot{\Sigma}$, using evaluation at these
points to define the homomorphism
$$
\Psi : \ker\mathbf{D} \to E_{\zeta_1} \oplus \ldots \oplus E_{\zeta_{n/2}}
\oplus V_{z_1} \oplus \ldots \oplus V_{z_k}.
$$
The right hand side is again a vector space of real dimension $n + k$,
and the same argument as above shows that $\Psi$ is injective.

To deal with the case $\ind(\mathbf{D}) \ge 0$, we consider the
formal adjoint $\mathbf{D}^*$ (cf.~\cite{Schwarz}).  This can be regarded
as a Cauchy-Riemann type operator
on the bundle $F := \overline{\Hom}_\CC(T\dot{\Sigma},E) \to \dot{\Sigma}$
with an appropriate totally real boundary condition $\ell^*$ and
asymptotic operators $\mathbf{A}^*_z$, which have the same parity as
$\mathbf{A}_z$.  It satisfies 
\begin{equation}
\label{eqn:adjointIndex}
\ind(\mathbf{D}^*) = -\ind(\mathbf{D}),
\qquad
\dim\ker\mathbf{D}^* = \dim\coker\mathbf{D},
\end{equation}
and applying
\eqref{eqn:2c1} to $\mathbf{D}$ and $\mathbf{D}^*$ together, we find
$$
c_1(E,\ell,\mathbf{A}_\Gamma) - c_1(F,\ell^*,\mathbf{A}^*_\Gamma) =
\frac{1}{2}\left[ \ind(\mathbf{D}) - \ind(\mathbf{D}^*) \right]
= \ind(\mathbf{D}).
$$
Then the condition $c_1(F,\ell^*,\mathbf{A}^*_\Gamma) < 0$ is 
satisfied if and only if
$\ind(\mathbf{D}) > c_1(E,\ell,\mathbf{A}_\Gamma)$, and this implies
$\mathbf{D}$ is surjective.  If $\ind(\mathbf{D}) \le 
c_1(E,\ell,\mathbf{A}_\Gamma)$,
then $c_1(F,\ell^*,\mathbf{A}^*_\Gamma) \ge 0$ and we can apply 
the above estimate to $\dim\ker\mathbf{D}^*$, giving
\begin{equation*}
\begin{split}
\dim\ker\mathbf{D} &= \ind(\mathbf{D}) + \dim\coker\mathbf{D} =
\ind(\mathbf{D}) + \dim\ker\mathbf{D}^* \\
&\le \ind(\mathbf{D}) + 
K\left(c_1(F,\ell^*,\mathbf{A}^*_\Gamma),\#\Gamma_0\right) \\   
&= \ind(\mathbf{D}) + 
K\left(c_1(E,\ell,\mathbf{A}_\Gamma) - \ind(\mathbf{D}),\#\Gamma_0\right).
\end{split}
\end{equation*}
\end{proof}

\begin{remark}
The proof of Theorem~\ref{thm:orbifold} requires only the very simplest
case of this dimension bound, namely that $\ker\mathbf{D}$ is trivial when
$c_1(E,\ell,\mathbf{A}_\Gamma) < 0$.  As that proof will demonstrate, 
however, such bounds can sometimes be useful in cases where
$\mathbf{D}$ is not surjective, so perhaps the more general dimension bound
will eventually find similar application.
\end{remark}

\section{The normal operator for a holomorphic curve}
\label{sec:normal}

In this section we will give the precise definition of regularity and
show that it is equivalent to the surjectivity of a certain Cauchy-Riemann
operator on a generalized normal bundle.  The precise relation between this
operator and the concept of regularity is stated in \S\ref{subsec:splitting}
as Theorem~\ref{thm:normalCR}, and
in \S\ref{subsec:four} we apply the linear transversality theory from 
\S\ref{subsec:linebundle} to show that Theorem~\ref{thm:criterion} follows
as an easy corollary.

Throughout the following, we fix a compact, connected and oriented surface 
$\Sigma$ of genus $g \ge 0$ with $m \ge 0$ boundary components, and a finite
set $\Gamma \subset \interior{\Sigma}$, writing $\dot{\Sigma} =
\Sigma \setminus\Gamma$.

\subsection{Teichm\"uller slices and Cauchy-Riemann operators}
\label{subsec:Teichmueller}

We begin by collecting some classical facts about moduli spaces of Riemann
surfaces which can be related to the analysis of
Cauchy-Riemann type operators.

Let $\jJ(\Sigma)$ denote the space of smooth complex structures on
$\Sigma$ that induce the given orientation, and denote by
$\Diff_+(\Sigma,\Gamma)$ the group of orientation preserving 
diffeomorphisms on $\Sigma$ 
that fix $\Gamma$, and $\Diff_0(\Sigma,\Gamma) \subset \Diff_+(\Sigma,\Gamma)$
those which are homotopic to the identity.
Both of these groups act on $\jJ(\Sigma)$ by $(\varphi,j) \mapsto
\varphi^* j$, and
the \emph{Teichm\"uller space} of $\dot{\Sigma}$ is a smooth 
finite-dimensional manifold defined as
$$
\tT(\dot{\Sigma}) = \jJ(\Sigma) / \Diff_0(\Sigma,\Gamma).
$$
Its quotient by the mapping
class group $M(\dot{\Sigma}) = \Diff_+(\Sigma,\Gamma) / \Diff_0(\Sigma,\Gamma)$
gives the \emph{moduli space of Riemann surfaces} (with genus $g$, $m$
boundary components and $\#\Gamma$ interior marked points)
$$
\mM(\dot{\Sigma}) = \tT(\dot{\Sigma}) / M(\dot{\Sigma}) =
\jJ(\Sigma) / \Diff_+(\Sigma,\Gamma),
$$
which is in general an orbifold of the same dimension.  We say
$\dot{\Sigma}$ is \emph{stable} if $\chi(\dot{\Sigma}) < 0$, in which case
$$
\dim \mM(\dot{\Sigma}) = 6g - 6 + 3m + 2\#\Gamma = -3\chi(\dot{\Sigma}) -
\#\Gamma,
$$ 
and the \emph{automorphism group}
$$
\Aut(\dot{\Sigma},j) = \{ \varphi \in \Diff_+(\Sigma,\Gamma)\ |\ 
\varphi^*j = j \}
$$
is finite for any choice of $j\in \jJ(\Sigma)$ 
(though its order may depend on~$j$).  Let $\DD \subset \CC$ denote the
closed unit disk, $\AA = [0,1] \times S^1$ the compact annulus
and $\TT^2 = \RR^2 / \ZZ^2$ the $2$--dimensional torus.
For our purposes, the non-stable cases to be considered are the following:
\begin{enumerate}
\item $\mM(S^2) = \{[i]\}$ and $\dim\Aut(S^2,i) = 6$.
\item $\mM(\CC) = \{[i]\}$ and $\dim\Aut(\CC,i) = 4$.
\item $\mM(\DD) = \{[i]\}$ and $\dim\Aut(\DD,i) = 3$.
\item $\mM(\RR\times S^1) = \{[i]\}$ and $\dim\Aut(\RR\times S^1,i) = 2$.
\item $\mM(\DD\setminus\{0\}) = \{[i]\}$ 
and $\dim\Aut(\DD\setminus\{0\},i) = 1$.
\item $\dim \mM(\AA) = 1$ and $\dim\Aut(\AA,j) = 1$.
\item $\dim \mM(\TT^2) = 2$ and $\dim\Aut(\TT^2,j) = 2$.
\end{enumerate}
For all but the last case, the mapping class group $M(\dot{\Sigma})$ is
trivial and thus $\mM(\dot{\Sigma}) = \tT(\dot{\Sigma})$ is a manifold.
Observe also that if $\dot{\Sigma}$ is not stable,
\begin{equation}
\label{eqn:TeichmuellerFredholm}
\dim\Aut(\dot{\Sigma},j) - \dim\mM(\dot{\Sigma}) = 3\chi(\dot{\Sigma})
+ \#\Gamma.
\end{equation}
Fixing $p > 2$, the latter is the Fredholm index of the standard linear
Cauchy-Riemann operator
$$
\mathbf{D}^{\Sigma}_\Gamma : W^{1,p}_{T(\p\Sigma)}(T\Sigma ; \Gamma) \to
L^p(\overline{\End}_\CC(T\Sigma)),
$$
where $W^{1,p}_{T(\p\Sigma)}(T\Sigma ; \Gamma)$ 
is the space of $W^{1,p}$--smooth vector fields
$Y : \Sigma \to T\Sigma$ satisfying $Y(\p\Sigma) \subset T(\p\Sigma)$ and
$Y|_{\Gamma} = 0$.

\begin{lemma}
\label{lemma:cokernel}
For all choices of $(\Sigma,j,\Gamma)$,
$\dim \coker(\mathbf{D}^{\Sigma}_\Gamma) = \dim\mM(\dot{\Sigma})$.
\end{lemma}
\begin{proof}
This may be regarded as a standard piece of Teichm\"uller theory in the
stable case (cf.~\cite{Tromba}), and also follows by using 
a simplified version of the argument in the proof of
Prop.~\ref{prop:transversality} to show that $\mathbf{D}^\Sigma_\Gamma$
is injective.  Here one must account also for the condition
$Y|_{\Gamma} = 0$, which ensures $Z(Y) \ge \#\Gamma$, thus it suffices to
observe that the adjusted first Chern number is strictly less than
$\#\Gamma$.  In the non-stable case,
a similar argument shows that $\dim\ker(\mathbf{D}^\Sigma_\Gamma) \le
\dim\Aut(\dot{\Sigma},j)$, and by interpreting $\mathbf{D}^\Sigma_\Gamma$
as the linearization of a nonlinear operator
$\dbar_j \varphi = T\varphi + j \circ T\varphi \circ j$, one sees that
$\ker(\mathbf{D}^\Sigma_\Gamma)$ contains $\aut(\dot{\Sigma},j)$,
giving an inequality in the other direction, hence
$$
\dim\ker(\mathbf{D}^\Sigma_\Gamma) = \dim\Aut(\dot{\Sigma},j).
$$
The result then follows from \eqref{eqn:TeichmuellerFredholm}.
\end{proof}

Given $j \in \jJ(\Sigma)$ and the corresponding Cauchy-Riemann operator
$\mathbf{D}^\Sigma_\Gamma$, pick a complement of 
$\im(\mathbf{D}^\Sigma_\Gamma)$, i.e.~a subspace $C \subset
L^p(\overline{\End}_\CC(T\Sigma))$ such that
$$
\im(\mathbf{D}^\Sigma_\Gamma) \oplus C = L^p(\overline{\End}_\CC(T\Sigma)).
$$
By approximation, we may assume every section in $C$ is smooth and
vanishes on a neighborhood of $\Gamma$.  We can then choose a small
neighborhood $\oO \subset C$ of~$0$ and define the map
\begin{equation}
\label{eqn:Phi}
\Phi : \oO \to \jJ(\Sigma) : y \mapsto \left( \1 + \frac{1}{2} j y\right)
j \left( \1 + \frac{1}{2} j y \right)^{-1},
\end{equation}
which has the properties $\Phi(0) = j$ and
$$
\left.\frac{\p}{\p t} \Phi(ty)\right|_{t=0} = y,
$$
thus it is injective if $\oO$ is sufficiently small.  The image 
$$
\tT := \Phi(\oO) \subset \jJ(\Sigma)
$$
is thus a smooth manifold of dimension
$\dim C = \dim \tT(\dot{\Sigma})$, with $T_j\tT = C$ and consisting of
smooth complex structures close to $j$ that are identical to~$j$ 
on some fixed neighborhood of $\Gamma$.
It parametrizes a neighborhood of $[j]$ in $\tT(\dot{\Sigma})$,
i.e.~the projection $\jJ(\Sigma) \to \tT(\dot{\Sigma})$ restricts to a
diffeomorphism from $\tT$ onto a neighborhood of~$[j]$.  This provides an
explicit construction of the following general object:

\begin{defn}
\label{defn:TeichmuellerSlice}
Given $j \in \jJ(\Sigma)$, we define a \emph{Teichm\"uller slice through $j$}
to be any smooth family $\tT \subset \jJ(\Sigma)$ parametrized by an
injective map $\uU \to \jJ(\Sigma)$, where $\uU$ is an open subset of
$\RR^{\dim\tT(\dot{\Sigma})}$, such that all $j' \in \tT$ are identical
on some fixed neighborhood of $\Gamma$, and
$\im(\mathbf{D}^\Sigma_\Gamma) \oplus T_j\tT = 
L^{p}(\overline{\End}_\CC(T\Sigma))$.
\end{defn}

\begin{lemma}
\label{lemma:Ginvariant}
If $(\dot{\Sigma},j)$ is stable, then there exists a Teichm\"uller slice $\tT$ 
through $j$ that is invariant under the group action
$\Aut(\dot{\Sigma},j) \times \jJ(\Sigma) \to \jJ(\Sigma) :
(\varphi,j') \mapsto \varphi^*j'$.
\end{lemma}
\begin{proof}
The automorphism group $G := \Aut(\dot{\Sigma})$ is finite and consists of
biholomorphic maps on $(\Sigma,j)$ that fix~$\Gamma$.  Each point in
$\Gamma$ then has a $G$--invariant neighborhood biholomorphically equivalent 
to the standard unit disk, on which $G$ acts by rational rotations.
Let $g$ denote a metric on $\Sigma$ that is invariant under the action of~$G$;
such a metric can be constructed by starting from the Poincar\'{e} metric
on $\dot{\Sigma}$ and interpolating this with flat rotation-invariant
metrics on disk-like neighborhoods of each point in~$\Gamma$.  Then $g$
induces a bundle metric on $\overline{\End}_\CC(T\Sigma) \to \Sigma$ and
a corresponding $G$--invariant $L^2$--inner product 
$\langle\ ,\ \rangle_{L^2}$ on the space of sections of this bundle.

To prove the lemma, it suffices to construct a $G$--invariant complement 
$C \subset L^p(\overline{\End}_\CC(T\Sigma))$ of 
$\im(\mathbf{D}_\Gamma^\Sigma)$ that consists of smooth sections vanishing
near~$\Gamma$: then an appropriate Teichm\"uller slice can be defined via
\eqref{eqn:Phi} since $\varphi^*j = j$ implies
$\Phi(\varphi^*y) = \varphi^*\Phi(y)$ for any $\varphi \in G$.
Observe that $\im(\mathbf{D}^\Sigma_\Gamma)$ itself is
$G$--invariant, since $\varphi^*j = j$ also implies
$\mathbf{D}^\Sigma_\Gamma(\varphi^*Y) = \varphi^*(\mathbf{D}^\Sigma_\Gamma Y)$
for all $Y \in W^{1,p}_{T(\p\Sigma)}(T\Sigma ; \Gamma)$.  Now using the
$G$--invariant $L^2$--product chosen above, define a complement $C_0$ as the
$L^2$--orthogonal complement of $\im(\mathbf{D}_\Gamma^\Sigma)$, i.e.
$$
C_0 = \left\{ y \in L^p(\overline{\End}_\CC(T\Sigma))\ \big|\ \left\langle
\mathbf{D}_\Gamma^\Sigma Y , y \right\rangle_{L^2} = 0 \text{ for all
$Y \in W^{1,p}_{T(\p\Sigma)}(T\Sigma ; \Gamma)$} \right\}.
$$
This space is $G$--invariant due to the $G$--invariance of 
$\im(\mathbf{D}_\Gamma^\Sigma)$ and $\langle\ ,\ \rangle_{L^2}$, and by elliptic
regularity for weak solutions of the formal adjoint equation, it consists 
only of smooth sections.  Now choosing $G$--invariant disk-like neighborhoods
of the points in $\Gamma$, we can obtain the desired complement $C$ by
multiplying the sections in $C_0$ by $G$--invariant cutoff functions that
vanish near~$G$ and equal~$1$ outside a sufficiently small neighborhood
of~$\Gamma$.
\end{proof}

For the two non-stable cases in which $\tT(\dot{\Sigma})$ is nontrivial,
it will be convenient to have explicit global Teichm\"uller slices.
If $\dot{\Sigma} = \AA = [0,1] \times S^1$, for each $\tau > 0$ define the
diffeomorphism $\varphi_\tau : \AA \to [0,\tau]
\times S^1 : (s,t) \mapsto (\tau s,t)$
and let $\tT_\AA$ denote the collection of complex structures
$\{ \varphi_\tau^*i \}_{\tau > 0}$.  This parametrizes the entirety of
$\tT(\AA)$ (which equals $\mM(\AA)$ since the mapping class
group is trivial), and also gives a natural identification of every
$\Aut(\AA,\varphi_\tau^*i)$ with $S^1$, acting on $\AA$ by translation of the
second factor.
If $\dot{\Sigma} = \TT^2 = \RR^2 / \ZZ^2$, we define $\tT_{\TT^2}$ to be the 
space
of all \emph{constant} complex structures on $\RR^2 = \CC$ that are compatible
with the standard orientation; clearly these descend to $\TT^2$, and they
also parametrize the entirety of $\tT(\TT^2)$.  Then for each
$j \in \tT_{\TT^2}$, the subgroup
$$
\Aut_0(\TT^2,j) := \Aut(\TT^2,j) \cap \Diff_0(\TT^2)
$$
can be identified naturally with $\TT^2$, acting by translations.
Choosing a base point $p = [(0,0)] \in \TT^2$, the stabilizer of 
$[j] \in \tT(\TT^2)$ under the action of
$M(\TT^2) = \SL(2,\ZZ)$ is meanwhile the finite subgroup
$$
\Aut(\TT^2,j ; p) := \{ \varphi \in \Aut(\TT^2,j)\ |\ \varphi(p) = p \}
= \{ A \in \SL(2,\ZZ)\ |\ A j = j A \},
$$
and $\Aut(\TT^2,j)$ is the semidirect product of
$\Aut(\TT^2,j ; p)$ with $\Aut_0(\TT^2,j) = \TT^2$.
Note in particular that for any $j \in \tT_{\TT^2}$, this group acts
by affine transformations on $\RR^2$ (descending to $\TT^2$), and
the action $(\varphi,j') \mapsto \varphi^*j'$ therefore preserves
$\tT_{\TT^2}$.

The following will be useful for technical reasons in our analysis of
the relationship between $\mathbf{D}_u$ and its normal component.
\begin{lemma}
\label{lemma:vanishAtCritical}
For any $j \in \jJ(\Sigma)$ and finite set $K \subset \dot{\Sigma}$, 
there exists a Teichm\"uller slice $\tT$ through $j$ such that 
every $j' \in \tT$ is identical to $j$ on some fixed neighborhood
of $K \cup \Gamma$.
\end{lemma}
\begin{proof}
It suffices to construct $C = T_j\tT$ so that every $y \in C$ vanishes
near $K \cup \Gamma$.  This can be done using cutoff functions to replace
a basis of any given complement with one that vanishes in such a neighborhood;
the new basis can be made $L^p$--close to the old one if the neighborhood
is sufficiently small.
\end{proof}

For any Teichm\"uller slice $\tT$ through $j$, the operator
$$
\mathbf{L}_\Gamma :
T_j \tT \oplus W^{1,p}(T\Sigma ; \Gamma) \to L^p(\overline{\End}_\CC(T\Sigma))
: (y,Y) \mapsto jy + \mathbf{D}^\Sigma_\Gamma Y
$$
is clearly surjective; indeed, since $\mathbf{D}^\Sigma_\Gamma$ is
complex linear, $\mathbf{L}_\Gamma(y,jY) = 
j\left( y + \mathbf{D}^\Sigma_\Gamma Y \right)$, and the target space
is spanned by $T_j\tT$ and $\im(\mathbf{D}^\Sigma_\Gamma)$.
For the analysis in the following
sections it will be useful to derive a corresponding statement for the
standard Cauchy-Riemann operator on a Riemann surface with
ends.  We will recall in the next section the construction of
certain Banach manifolds containing
asymptotically cylindrical maps $\dot{\Sigma} \to W$.  In the simple case
$W = \dot{\Sigma}$, the tangent space to such a Banach manifold $\bB^\Sigma$ 
at the identity map $\1 : \dot{\Sigma} \to \dot{\Sigma}$ can be written as
$$
T_{\1} \bB^{\Sigma} = W^{1,p,\boldsymbol{\delta}}_{T(\p\Sigma)}(T\dot{\Sigma})
\oplus V^{\Sigma}_\Gamma,
$$
where $\boldsymbol{\delta} > 0$ is a small weight applying at every end and
$V^{\Sigma}_\Gamma \subset \Gamma(T\dot{\Sigma})$ is a 
$2\#\Gamma$--dimensional space of
smooth sections that are supported near infinity and constant in some
fixed choice of cylindrical
coordinates near each end.  
The natural nonlinear Cauchy-Riemann operator defines a
section of a Banach space bundle over $\bB^\Sigma$, 
whose linearization at $\1$ is
the usual linear Cauchy-Riemann operator given by the holomorphic
structure of $T\dot{\Sigma} \to \dot{\Sigma}$, denoted here by
$$
\mathbf{D}^\Sigma : W^{1,p,\boldsymbol{\delta}}_{T(\p\Sigma)}(T\dot{\Sigma})
\oplus V^\Sigma_\Gamma \to
L^{p,\boldsymbol{\delta}}(\overline{\End}_\CC(T\dot{\Sigma})).
$$
Now since every $y \in T_j\tT$ is smooth and vanishes near $\Gamma$, 
there is a natural inclusion of $T_j\tT \subset
L^{p,\boldsymbol{\delta}}(\overline{\End}_\CC(T\dot{\Sigma}))$,
as well as a bounded linear map $T_j\tT \to 
L^{p,\boldsymbol{\delta}}(\overline{\End}_\CC(T\dot{\Sigma})) : y \mapsto jy$.

\begin{lemma}
\label{lemma:Teichmueller}
The operator
\begin{equation*}
\begin{split}
\mathbf{L} :
T_j \tT \oplus \left( W^{1,p,\boldsymbol{\delta}}_{T(\p\Sigma)}(T\dot{\Sigma}) \oplus
V^\Sigma_\Gamma \right) &\to L^{p,\boldsymbol{\delta}}(\overline{\End}_\CC(T\dot{\Sigma})) 
\\
(y,\eta) &\mapsto j y + \mathbf{D}^\Sigma \eta
\end{split}
\end{equation*}
is surjective.
\end{lemma}
\begin{proof}
Applying the linear theory in \S\ref{sec:linear}, we find that
$\ind(\mathbf{D}^\Sigma) = \ind(\mathbf{D}^\Sigma_\Gamma)$ and hence
$\ind(\mathbf{L}) = \ind(\mathbf{L}_\Gamma)$.  Now in light of
the natural inclusion 
$$
W^{1,p,\boldsymbol{\delta}}_{T(\p\Sigma)}(T\dot{\Sigma}) \oplus
V^\Sigma_\Gamma \hookrightarrow W^{1,p}_{T(\p\Sigma)}(T\Sigma ; \Gamma),
$$
we have $\ker(\mathbf{L}) \subset \ker(\mathbf{L}_\Gamma)$,
but since $\dim\ker(\mathbf{L}) \ge \ind(\mathbf{L})$ in
general and $\dim\ker(\mathbf{L}_\Gamma) = \ind(\mathbf{L}_\Gamma)$
by the remarks above, it follows that
$$
\dim\ker(\mathbf{L}) \le \dim\ker(\mathbf{L}_\Gamma) =
\ind(\mathbf{L}_\Gamma) = \ind(\mathbf{L}),
$$
implying $\mathbf{L}$ is surjective.
\end{proof}

\begin{cor}
\label{cor:directSum}
For any Teichm\"uller slice $\tT$ through $j$,
$L^{p,\boldsymbol{\delta}}(\overline{\End}_\CC(T\dot{\Sigma})) =
\im(\mathbf{D}^\Sigma) \oplus T_j\tT$.
\end{cor}

\subsection{Functional analytic setup}
\label{subsec:nonlinear}

To state the definition of regularity, we begin by reviewing
the nonlinear functional analytic setup used in
\cites{Dragnev,Bourgeois:thesis} for asymptotically cylindrical
maps $u : \dot{\Sigma} \to W$ with nondegenerate or Morse-Bott 
asymptotic orbits.  Fix the surface $\Sigma$, punctures
$\Gamma = \Gamma^+ \cup \Gamma^- \subset \interior{\Sigma}$,
and asymptotic constraints $\mathbf{c}$.  
Recall that the latter choice partitions
$\Gamma$ into a set of \emph{constrained} and \emph{unconstrained} punctures
$\Gamma = \Gamma_C \cup \Gamma_U$,
and assigns to each $z \in \Gamma^\pm_C$ an orbit $P_z$ of $X_\pm$,
which we will assume is Morse-Bott, and we will denote its period by $T_z$.
For each $z \in \Gamma_U^\pm$, we instead choose an arbitrary Morse-Bott 
manifold of periodic orbits in $M_\pm$, denoted again by $P_z$,
with period $T_z > 0$.  By a slight abuse of notation, each $P_z$ may be
regarded both as a submanifold $P_z \subset M_\pm$ and as a set
of $T_z$--periodic orbits $\gamma \in P_z$ (sometimes with only one
element).  Denote this collection of choices for all punctures $z \in \Gamma$
by $P_\Gamma$.  We shall then consider a Banach manifold consisting
of asymptotically cylindrical maps $u : \dot{\Sigma} \to W$ whose
asymptotic orbits $\gamma_z$ for $z \in \Gamma$ satisfy $\gamma_z \in P_z$.

Before explaining the Banach manifold, we digress for a moment to
define some important invariants that enter into the index formula.
Recall that any $T$--periodic orbit $\gamma$ of $X_\pm$ has an associated
\emph{asymptotic operator} $\mathbf{A}_\gamma$, defined on sections of
the bundle $\xi_\pm$ along~$\gamma$.  One can write it down by choosing
a parametrization $x : S^1 \to M_\pm$ of $\gamma$ with 
$\lambda(\dot{x}) \equiv T$, and defining $\mathbf{A}_\gamma :
\Gamma(x^*\xi_\pm) \to \Gamma(x^*\xi_\pm)$ by
$$
\mathbf{A}_\gamma v = -J_\pm (\nabla_t v - T \nabla X_\pm)
$$
for any symmetric connection $\nabla$ on $M_\pm$.  This gives an unbounded
self-adjoint operator on $L^2(x^*\xi_\pm)$ of the form considered in
\S\ref{sec:linear}, and it is nondegenerate if and only if the orbit is
nondegenerate, in which case we define the Conley-Zehnder index
$\muCZ^\Phi(\gamma) = \muCZ^\Phi(\mathbf{A}_\gamma)$ for any choice of
trivialization $\Phi$ on $x^*\xi_\pm$.  If $\gamma$ is degenerate,
then $\mathbf{A}_\gamma$ can be perturbed to a nondegenerate
asymptotic operator by adding any number $\boldsymbol{\epsilon} 
\in \RR \setminus
-\sigma(\mathbf{A})$, and we thus define the \emph{perturbed} 
Conley-Zehnder index
$$
\muCZ^\Phi(\gamma + \boldsymbol{\epsilon}) := 
\muCZ^\Phi(\mathbf{A}_\gamma + \boldsymbol{\epsilon}),
$$
and its \emph{parity}
$$
p(\gamma + \boldsymbol{\epsilon}) = 
\begin{cases}
0 & \text{ if $\muCZ^\Phi(\gamma + \boldsymbol{\epsilon})$ is even,}\\
1 & \text{ if $\muCZ^\Phi(\gamma + \boldsymbol{\epsilon})$ is odd,}
\end{cases}
$$
which does not depend on~$\Phi$.
Observe that if $\gamma$ is nondegenerate and $\boldsymbol{\epsilon}$ 
is sufficiently close to 
zero, then $\muCZ^\Phi(\gamma + \boldsymbol{\epsilon}) = 
\muCZ^\Phi(\gamma)$ since 
$\sigma(\mathbf{A}_\gamma)$ is discrete.  More generally, one can
see from the relationship between the Conley-Zehnder index and spectral
flow (cf.~\cite{RobbinSalamon}) that for sufficiently small 
$\boldsymbol{\epsilon} > 0$,
\begin{equation}
\label{eqn:spectralFlow}
\muCZ^\Phi(\gamma  -\boldsymbol{\epsilon}) - 
\muCZ^\Phi(\gamma + \boldsymbol{\epsilon}) =
\dim\ker(\mathbf{A}_\gamma).
\end{equation}
In particular if $\gamma$ belongs to a Morse-Bott family $P$, then
the right hand side of \eqref{eqn:spectralFlow} is $\dim P - 1$,
and $\muCZ^\Phi(\gamma \pm \boldsymbol{\epsilon})$ remains unchanged if we
move $\gamma$ to a different orbit in the same family.

If $M_\pm$ are $3$--dimensional, then $\xi_\pm$ have complex rank 
one, so recalling the definitions in \S\ref{subsec:linebundle}, 
we can associate to any $T$--periodic orbit $\gamma$ of $X_\pm$
and real number $\boldsymbol{\epsilon}$ the so-called
\emph{extremal winding numbers}
$$
\alpha_\pm^\Phi(\gamma + \boldsymbol{\epsilon}) := 
\alpha_\pm^\Phi(\mathbf{A}_\gamma + \boldsymbol{\epsilon}),
$$
or for the case $\boldsymbol{\epsilon} = 0$, simply
$\alpha_\pm^\Phi(\gamma) = \alpha_\pm^\Phi(\mathbf{A}_\gamma)$.
We will refer to the eigenfunctions of $\mathbf{A}_\gamma$ involved in this
definition as \emph{extremal eigenfunctions} at $\gamma$ if 
$\boldsymbol{\epsilon} = 0$, or more generally \emph{extremal eigenfunctions
with respect to $\boldsymbol{\epsilon}$}.
Now if $\boldsymbol{\epsilon} \not\in -\sigma(\mathbf{A}_\gamma)$,
\eqref{eqn:CZwinding} gives
\begin{equation*}
\begin{split}
\muCZ^\Phi(\gamma + \boldsymbol{\epsilon}) &= 2\alpha_\mp^\Phi(\gamma +
\boldsymbol{\epsilon}) \pm p(\gamma + \boldsymbol{\epsilon}) \\
p(\gamma + \boldsymbol{\epsilon}) &= 
\alpha_+^\Phi(\gamma + \boldsymbol{\epsilon}) 
- \alpha_-^\Phi(\gamma + \boldsymbol{\epsilon}) \in \{0,1\}.
\end{split}
\end{equation*}
Choosing $\boldsymbol{\delta} > 0$ arbitrarily small, it will also 
be convenient to define
\begin{equation}
\label{eqn:nu}
\nu_\pm(\gamma) = \alpha_\pm^\Phi(\gamma - \boldsymbol{\delta}) -
\alpha_\pm^\Phi(\gamma + \boldsymbol{\delta}),
\end{equation}
which equals~$0$ whenever $\gamma$ is nondegenerate, and is otherwise
either~$0$ or~$1$.\footnote{For orbits in two-dimensional families,
the numbers $\nu_\pm(\gamma)$ are closely related to the \emph{sign}
of a Morse-Bott surface, as defined in \cite{Wendl:BP1}.}

\begin{notation}
Fix a number $\boldsymbol{\delta} > 0$, which we will generally assume to
be as small as may be needed.
Suppose $\Gamma = \Gamma^+ \cup \Gamma^-$ is a set of punctures and
$\mathbf{c}$ is a set of asymptotic constraints, defining
constrained and unconstrained subsets $\Gamma_C,\Gamma_U \subset \Gamma$ 
respectively.  We then associate to each puncture $z \in \Gamma$ 
a real number
\begin{equation}
\label{eqn:cz}
\mathbf{c}_z := 
\begin{cases}
\boldsymbol{\delta} & \text{ if $z \in \Gamma_C$,} \\
-\boldsymbol{\delta} & \text{ if $z \in \Gamma_U$.}
\end{cases}
\end{equation}
\end{notation}

For asymptotically cylindrical maps $u : \dot{\Sigma} \to W$ subject to
constraints $\mathbf{c}$, we will use
the following notational conventions throughout.  The asymptotic orbit of $u$ 
at a puncture
$z \in \Gamma^\pm$ will be called $\gamma_z$, with asymptotic operator
$\mathbf{A}_z := \mathbf{A}_{\gamma_z}$, and the collection of these
for all punctures will be denoted by $\gamma_\Gamma$ and
$\mathbf{A}_\Gamma$ respectively.  Denote the
corresponding collection of \emph{perturbed} asymptotic operators
$\{\mathbf{A}_z \pm \mathbf{c}_z\}_{z\in\Gamma^\pm}$ 
by $\mathbf{A}_\Gamma \pm \mathbf{c}_\Gamma$, noting that the sign
choice must always match the sign of the puncture.  For $i \in \{0,1\}$,
let
$$
\Gamma_i^\pm(\mathbf{c}) = \{ z \in \Gamma^\pm\ |\ 
p(\gamma_z \pm \mathbf{c}_z) = i \}
$$
and $\Gamma_i(\mathbf{c}) = \Gamma_i^+(\mathbf{c}) \cup 
\Gamma_i^-(\mathbf{c})$.  This defines a partition of $\Gamma$ into so-called
\emph{even} and \emph{odd} punctures with respect to the constraints.
Note that when $\gamma_z$ is nondegenerate, the parity of $z$ is simply
the even/odd parity of $\muCZ^\Phi(\gamma_z)$; in general however this
distinction depends on not just the orbit
and constraints, but also the \emph{sign} of the puncture.

Choose $p > 2$ and define the Banach manifold
$$
\bB := \bB^{1,p,\boldsymbol{\delta}}(\dot{\Sigma},W ; L,P_\Gamma)
$$
to consist of maps $\dot{\Sigma} \to W$ of class $W^{1,p}_{\text{loc}}$
which satisfy $u(\p\Sigma) \subset L$ and have asymptotically cylindrical 
behavior approaching the orbits $\{P_z\}_{z\in\Gamma}$ 
at the corresponding punctures: the latter means in particular
that using cylindrical coordinates $(s,t) \in Z_\pm$ near $z \in \Gamma^\pm$,
there exists an orbit $\gamma_z \in P_z$ 
with parametrization $x : \RR \to M_\pm$
and a constant $s_0$ such that for sufficiently large~$|s|$,
$$
u(s + s_0,t) = \exp_{\tilde{x}(s,t)} h(s,t),
$$
where $\tilde{x}(s,t) := (T_z s,x(T_z t))
\in E_\pm \subset \RR\times M_\pm$ and $h \in \Gamma(\tilde{x}^*T E_\pm)$
is of weighted Sobolev class $W^{1,p,\boldsymbol{\delta}}$ on $Z_\pm$.
The tangent space $T_u \bB$ can then be written as
$$
T_u \bB = W_\Lambda^{1,p,\boldsymbol{\delta}}(u^*TW) \oplus V_\Gamma \oplus X_\Gamma,
$$
where the summands are defined as follows.
The subscript $\Lambda$ refers to the totally real subbundle
\begin{equation}
\label{eqn:Lambda}
\Lambda := \left(u|_{\p\Sigma}\right)^*TL \to \p\Sigma,
\end{equation}
so that 
sections $v \in W_\Lambda^{1,p,\boldsymbol{\delta}}(u^*TW)$ are required to satisfy the
boundary condition $v(\p\Sigma) \subset \Lambda$, as well as decaying in
accordance with the small exponential weight $\boldsymbol{\delta} > 0$ at each end.  
The other two summands are both finite dimensional vector spaces consisting
of sections $\dot{\Sigma} \to u^*TW$ that are supported near infinity
and asymptotically equal to constant vectors in some choice of $\RR$--invariant
coordinates near the asymptotic orbit.  In particular, $V_\Gamma$ has
dimension $2 \#\Gamma$ and contains vector fields that are parallel
to the orbit cylinders $\tilde{x}(s,t) = (Ts,x(Tt))$ near infinity, while
the vector fields in $X_\Gamma$ are trivial whenever $P_z$ is a fixed
orbit and otherwise parallel to the Morse-Bott manifolds $P_z$, thus
\begin{equation}
\label{eqn:dimXGamma}
\dim X_\Gamma = \sum_{z \in \Gamma} (\dim P_z - 1) =
\sum_{z \in \Gamma_U} \dim\ker(\mathbf{A}_z).
\end{equation}

Fixing a complex structure $j$ on $\Sigma$, there is a Banach space bundle 
$\eE \to \bB$ whose fibers are spaces of complex antilinear bundle maps
$$
\eE_u = L^{p,\boldsymbol{\delta}}(\overline{\Hom}_\CC(T\dot{\Sigma},u^*TW)),
$$
and the nonlinear Cauchy-Riemann operator defines a smooth section
$$
\dbar_J : \bB \to \eE : u \mapsto Tu + J \circ Tu \circ j,
$$
whose zeroes are parametrizations of asymptotically cylindrical 
pseudoholomorphic curves $u : (\dot{\Sigma},j) \to (W,J)$.  The linearization
of $\dbar_J$ at a zero $u$ defines a linear Cauchy-Riemann type operator,
\begin{equation}
\label{eqn:linearCR}
\begin{split}
\mathbf{D}_u : \Gamma(u^*TW) &\to \Gamma(\overline{\Hom}_\CC(T\dot{\Sigma},
u^*TW)) \\
v &\mapsto \nabla v + J \circ \nabla v \circ j + (\nabla_v J) \circ Tu \circ j,
\end{split}
\end{equation}
where $\nabla$ is any symmetric connection on~$W$.  As a 
bounded linear operator $T_u \bB \to \eE_u$, $\mathbf{D}_u$ is Fredholm.
To write down its index,
let $\Phi$ be an arbitrary choice of trivialization for 
$u^*TW$ along $\p\Sigma$ and for $\xi_\pm$ along the orbits $\gamma_z$,
and define the \emph{total Maslov index}
$$
\mu^\Phi(u ; \mathbf{c}) = \mu^\Phi(u^*TW,\Lambda) + 
\sum_{z \in \Gamma^+} \muCZ^\Phi(\gamma_z + \mathbf{c}_z) -
\sum_{z \in \Gamma^-} \muCZ^\Phi(\gamma_z - \mathbf{c}_z).
$$
The trivializations of $\xi_\pm$ extend naturally to trivializations of
$TW = T(\RR\times M_\pm)$ along the orbits via the splitting
\begin{equation}
\label{eqn:splitting}
T(\RR\times M_\pm) = (\RR \oplus \RR X_\pm) \oplus \xi_\pm,
\end{equation}
so that one can also define the relative Chern number $c_1^\Phi(u^*TW)$.

\begin{prop}
\label{prop:indexD}
$$
\ind(\mathbf{D}_u) = n\chi(\dot{\Sigma}) + 2 c_1^\Phi(u^*TW) + 
\mu^\Phi(u ; \mathbf{c}) + \#\Gamma.
$$
\end{prop}
\begin{proof}
Denote by $\mathbf{D}_0$ the restriction of $\mathbf{D}_u$ to
$W^{1,p,\boldsymbol{\delta}}_\Lambda(u^*TW)$, so
\begin{equation*}
\begin{split}
\ind(\mathbf{D}_u) &= \ind(\mathbf{D}_0) + \dim V_\Gamma + \dim X_\Gamma \\
&= \ind(\mathbf{D}_0) + 2\#\Gamma + \sum_{z\in\Gamma_U} \dim\ker(\mathbf{A}_z).
\end{split}
\end{equation*}
Then $\mathbf{D}_0$ is a Cauchy-Riemann type operator asymptotic at 
$z \in \Gamma$ to the operators
$$
\mathbf{B}_z := \mathbf{C} \oplus \mathbf{A}_z,
$$
where we use the splitting \eqref{eqn:splitting} and define on the first
summand the ``trivial'' asymptotic operator $\mathbf{C} = -J_0 \frac{d}{dt}$.
The latter is degenerate, but we have
\begin{equation}
\label{eqn:CZofC}
\muCZ(\mathbf{C} \pm \boldsymbol{\delta}) = \mp 1
\end{equation}
if $\boldsymbol{\delta} > 0$ is sufficiently small.
By the discussion of exponential weights in \S\ref{subsec:general},
$\mathbf{D}_0$ is now conjugate to a Cauchy-Riemann operator
$W^{1,p}_\Lambda(u^*TW) \to L^p(\overline{\Hom}_\CC(T\dot{\Sigma},u^*TW))$
with nondegenerate asymptotics, so by \eqref{eqn:RiemannRoch}:
\begin{equation*}
\begin{split}
\ind(\mathbf{D}_0) &= n\chi(\dot{\Sigma}) + 2 c_1^\Phi(u^*TW) + 
\mu^\Phi(u^*TW,\Lambda) \\
&\qquad + \sum_{z \in \Gamma^+} \muCZ^\Phi(\mathbf{B}_z + \boldsymbol{\delta}) -
\sum_{z \in \Gamma^-} \muCZ^\Phi(\mathbf{B}_z - \boldsymbol{\delta}) \\
&= n\chi(\dot{\Sigma}) + 2 c_1^\Phi(u^*TW) + \mu^\Phi(u^*TW,\Lambda) \\
&\qquad + \sum_{z\in\Gamma^+} \muCZ^\Phi(\mathbf{A}_z + \boldsymbol{\delta})
- \sum_{z\in\Gamma^-} \muCZ^\Phi(\mathbf{A}_z - \boldsymbol{\delta}) - \#\Gamma,
\end{split}
\end{equation*}
where in the last line we've used the splitting $\mathbf{B}_z =
\mathbf{C} \oplus \mathbf{A}_z$ and \eqref{eqn:CZofC}.
Now using \eqref{eqn:spectralFlow}, we have
\begin{equation}
\label{eqn:aCalculation}
\begin{split}
\sum_{z\in\Gamma^+} &\muCZ^\Phi(\mathbf{A}_z + \boldsymbol{\delta})
- \sum_{z \in \Gamma^-} \muCZ^\Phi(\mathbf{A}_z - \boldsymbol{\delta}) - \#\Gamma \\
&\qquad + 2\#\Gamma + \sum_{z \in \Gamma_U} \dim\ker(\mathbf{A}_z) \\
&= \sum_{z\in\Gamma^+_U} \muCZ^\Phi(\mathbf{A}_z - \boldsymbol{\delta}) +
\sum_{z \in \Gamma^+_C} \muCZ^\Phi(\mathbf{A}_z + \boldsymbol{\delta}) \\
&\qquad - \sum_{z \in \Gamma^-_U} \muCZ^\Phi(\mathbf{A}_z + \boldsymbol{\delta})
- \sum_{z \in \Gamma^-_C} \muCZ^\Phi(\mathbf{A}_z - \boldsymbol{\delta}) + \#\Gamma \\
&= \sum_{z \in \Gamma^+} \muCZ^\Phi(\mathbf{A}_z + \mathbf{c}_z) 
- \sum_{z \in \Gamma^-} \muCZ^\Phi(\mathbf{A}_z - \mathbf{c}_z)
+ \#\Gamma,
\end{split}
\end{equation}
and the result follows.
\end{proof}

For the following lemma, recall that $\mathbf{D}^\Sigma : \Gamma(T\dot{\Sigma})
\to \Gamma(\overline{\End}_\CC(T\dot{\Sigma}))$ denotes the natural
linear Cauchy-Riemann operator on $\Gamma(T\dot{\Sigma})$ determined by
the holomorphic structure of $T\dot{\Sigma}$; it is the linearization at the
identity of the
operator $\dbar^\Sigma_j \varphi 
= T\varphi + j\circ T\varphi \circ j$ acting on
maps $\varphi : \dot{\Sigma} \to \dot{\Sigma}$.  We use the bundle map
$du : T\dot{\Sigma} \to u^*TW$ to define linear maps
\begin{equation}
\label{eqn:duCorrespondence}
\begin{split}
\Gamma(T\dot{\Sigma}) &\xrightarrow{du}
\Gamma(u^*TW), \\
\Gamma(\overline{\End}_\CC(T\dot{\Sigma})) &\xrightarrow{du}
\Gamma(\overline{\Hom}_\CC(T\dot{\Sigma},u^*TW).
\end{split}
\end{equation}

\begin{lemma}
\label{lemma:DuOfTangents}
For any smooth vector field $v \in \Gamma(T\dot{\Sigma})$,
$$
\mathbf{D}_u (du(v)) = du( \mathbf{D}^\Sigma v).
$$
\end{lemma}
\begin{proof}
Choose any open subset $\uU \subset \dot{\Sigma}$ with compact support.  
On this neighborhood, the flow 
$\varphi^{\tau}$ of $v$ is well defined for $\tau$ sufficiently close to~$0$,
and by definition, if $z \in \uU$ and $Y \in T_z\dot{\Sigma}$,
$$
(\mathbf{D}^\Sigma v) Y = \left. \nabla_{\tau} \left[ 
\dbar^\Sigma_j \varphi^\tau(Y) \right] \right|_{\tau=0},
$$
where $\nabla$ is any connection on $\dot{\Sigma}$.
Similarly, using the fact that $u : (\dot{\Sigma},j) \to (W,J)$ and
$\1 : (\dot{\Sigma},j) \to (\dot{\Sigma},j)$ are both holomorphic,
\begin{equation*}
\begin{split}
\mathbf{D}_u (du(v))(Y) &= \left. \nabla_\tau 
\left[ \dbar_J (u \circ \varphi^\tau)(Y)\right]\right|_{\tau=0} \\
&= \left. \nabla_\tau \left[ T(u \circ \varphi^\tau)(Y) + J \circ
T(u \circ \varphi^\tau) \circ j(Y) \right]\right|_{\tau=0} \\
&= \left. \nabla_\tau \left[ Tu \circ T\varphi^\tau(Y) + J \circ
Tu \circ T\varphi^\tau \circ j(Y) \right]\right|_{\tau=0} \\
&= \left. \nabla_\tau \left[ Tu \circ T\varphi^\tau(Y) + 
Tu \circ j \circ T\varphi^\tau \circ j(Y) \right]\right|_{\tau=0} \\
&= \left. \nabla_\tau \left[ du(\varphi^\tau(z)) \cdot
\dbar_j^\Sigma \varphi^\tau(Y) \right]\right|_{\tau=0} \\
&= \nabla_{v(z)}(du) \cdot \dbar_j^\Sigma(\1)(Y) + 
du \left( \left. \nabla_\tau \left[ 
\dbar_j^\Sigma \varphi^\tau(Y)\right] \right|_{\tau=0} \right) \\
&= du (\mathbf{D}^\Sigma v(Y)).
\end{split}
\end{equation*}
\end{proof}

Varying complex structures on the domain 
can be incorporated into the picture by fixing
$j_0 \in \jJ(\Sigma)$ and choosing
a Teichm\"uller slice $\tT$ through~$j_0$ 
(see Def.~\ref{defn:TeichmuellerSlice}).
We can now redefine the Banach space bundle $\eE$ over 
$\tT \times \bB$ so that
$$
\eE_{(j,u)} = L^{p,\boldsymbol{\delta}}\left(\overline{\Hom}_\CC((T\dot{\Sigma},j),
(u^*TW,J))\right),
$$
and extend the section $\dbar_J$ over this bundle by
$$
\dbar_J : \tT \times \bB \to \eE : (j,u) \mapsto Tu + J\circ Tu \circ j.
$$
The linearization at $(j,u) \in \dbar_J^{-1}(0)$ can now be expressed via
its ``partial derivatives,''
$$
D\dbar_J(j,u) : T_j \tT \oplus T_u \bB \to \eE_{(j,u)} :
(y,v) \mapsto \mathbf{G}_u y + \mathbf{D}_u v
$$
where 
\begin{equation*}
\begin{split}
\mathbf{G}_u : T_j \tT \subset \Gamma(\overline{\End}_\CC(T\dot{\Sigma})) &\to
\Gamma(\overline{\Hom}_\CC((T\dot{\Sigma},j),(u^*TW,J))) \\
y &\mapsto J \circ Tu \circ y.
\end{split}
\end{equation*}

We can now present the precise definition of regularity.

\begin{defn}
\label{defn:regular}
The curve $(\Sigma,j,\Gamma,u) \in \mM^{\mathbf{c}}$ 
is called \emph{regular} if there exists a Teichm\"uller slice $\tT$ 
through $j$ such that the operator
$D\dbar_J(j,u) : T_j \tT \oplus T_u\bB \to \eE_{(j,u)}$
is surjective.
\end{defn}

\begin{remark}
\label{remark:reparam}
This condition clearly doesn't depend on the choice of map 
$u : (\dot{\Sigma},j) \to (W,J)$ representing a given equivalence class in
$\mM^{\mathbf{c}}$; if $\varphi : (\Sigma',j') \to (\Sigma,j)$ is a 
biholomorphic map and $u' = u \circ \varphi$, one can use the pullback 
$\varphi^*$ to construct a Teichm\"uller slice $\tT'$ through $j'$ so that 
the operators $D\dbar_J(j',u')$ and $D\dbar_J(j,u)$ are conjugate.
The next lemma shows also that the surjectivity of
$D\dbar_J(j,u)$---and in fact the codimension of its image---does not depend
on the choice of Teichm\"uller slice.
\end{remark}

\begin{lemma}
\label{lemma:allSlices}
For $(\Sigma,j,\Gamma,u) \in \mM^{\mathbf{c}}$ and
any two Teichm\"uller slices $\tT$ and $\tT'$ through $j$, denote by
$\mathbf{L} : T_j\tT \oplus T_u\bB \to \eE_{(j,u)}$ and 
$\mathbf{L}' : T_j\tT' \oplus T_u\bB \to \eE_{(j,u)}$ the
corresponding linearizations of $\dbar_J$ at $(j,u)$.  Then
$\im(\mathbf{L}) = \im(\mathbf{L}')$.
\end{lemma}
\begin{proof}
Using the inclusion $T_j\tT \subset 
L^{p,\boldsymbol{\delta}}(\overline{\End}_\CC(T\dot{\Sigma}))$, extend
$\mathbf{L}$ to
\begin{equation*}
\begin{split}
\overline{\mathbf{L}} : L^{p,\boldsymbol{\delta}}
(\overline{\End}_\CC(T\dot{\Sigma})) \oplus T_u\bB &\to \eE_{(j,u)} \\
(y,v) &\mapsto J \circ Tu \circ y + \mathbf{D}_u v.
\end{split}
\end{equation*}
For $y = \mathbf{D}^\Sigma Y \in \im(\mathbf{D}^\Sigma) \subset 
L^{p,\boldsymbol{\delta}}(\overline{\End}_\CC(T\dot{\Sigma}))$,
we use the fact that $u$ is $J$--holomorphic and write
$$
\overline{\mathbf{L}}(y,0) = J \circ Tu \circ y = du(jy) =
du(\mathbf{D}^\Sigma (jY)).
$$
Then by Lemma~\ref{lemma:DuOfTangents}, this equals $\mathbf{D}_u(du(jY))$
if $Y$ is smooth, and the same holds for general
$Y$ in the domain of $\mathbf{D}^\Sigma$ by a density argument, thus
the restriction of $\overline{\mathbf{L}}$ to $\im(\mathbf{D}^\Sigma)$
has image contained in $\im(\mathbf{D}_u)$.  Since
$L^{p,\boldsymbol{\delta}}(\overline{\End}_\CC(T\dot{\Sigma})) =
\im(\mathbf{D}^\Sigma) \oplus T_j\tT$ by Cor.~\ref{cor:directSum},
this implies $\im(\overline{\mathbf{L}}) = \im(\mathbf{L})$.
Now using the same argument for $\tT'$, we have
$\im(\mathbf{L}) = \im(\overline{\mathbf{L}}) = \im(\mathbf{L}')$.
\end{proof}

Since $\mathbf{D}_u$ is Fredholm and $\tT$ is finite dimensional,
$D\dbar_J(j,u)$ is also Fredholm.
Recalling \eqref{eqn:TeichmuellerFredholm} and the definition of
$\ind(u ; \mathbf{c}) = \virdim \mM_u^{\mathbf{c}}$ in \eqref{eqn:index},
we have
$$
\ind D\dbar_J(j,u) = \dim \tT + \ind(\mathbf{D}_u) =
\ind(u ; \mathbf{c}) + \dim \Aut(\dot{\Sigma},j).
$$
For completeness, we now prove the fact that regularity gives
$\mM^{\mathbf{c}}$ the structure of a smooth orbifold of dimension 
$\ind(u ; \mathbf{c})$ near~$u$.

\begin{proof}[Proof of Theorem~\ref{thm:orbifoldFolk}]
Assume $(j_0,u_0) \in \dbar_J^{-1}(0)$ is regular and let
$G = \Aut(\dot{\Sigma},j_0)$.  By Lemma~\ref{lemma:allSlices}, the regularity
condition is independent of the choice of
Teichm\"uller slice, so if $\dot{\Sigma}$ is stable, then 
using Lemma~\ref{lemma:Ginvariant} we can
pick a slice $\tT$ through $j_0$ that is invariant under the natural
$G$--action.  Similarly if $\dot{\Sigma}$ is
$\AA$ or $\TT^2$, then without loss of generality we can compose with
a diffeomorphism such that $j_0$ belongs to one of the special Teichm\"uller
slices $\tT_{\AA}$ or $\tT_{\TT}$ constructed in \S\ref{subsec:Teichmueller} 
(which also admit a natural $G$--action), and choose this for~$\tT$.
There is now a $G$--action on $\tT \times \bB$ defined by
$$
\varphi \cdot (j,u) = (\varphi^*j , u \circ \varphi).
$$
This clearly preserves $\dbar_J^{-1}(0)$, and the stabilizer of any
$(j,u) \in \dbar_J^{-1}(0)$ is the finite subgroup
$\{ \varphi \in \Aut(\dot{\Sigma},j_0) |\ \text{$\varphi^*j = j$,
$u \circ \varphi = u$} \} \subset \Aut(u)$.
Since $D\dbar_J(j_0,u_0)$ is surjective, Remark~\ref{remark:reparam} implies
that the same is true for all $(j,u)$ in the $G$--orbit of $(j_0,u_0)$,
thus by the implicit function theorem, a neighborhood $\uU \subset
\dbar_J^{-1}(0)$ of this orbit admits a natural smooth manifold structure, with
dimension $\ind(u_0 ; \mathbf{c}) + \dim\Aut(\dot{\Sigma},j_0)$.  Starting
from a small neighborhood of $(j_0,u_0)$ in $\dbar_J^{-1}(0)$ and extending
this under the $G$--action, we may assume $\uU$ to be $G$--invariant.
The quotient $\uU / G$ then inherits the structure of a smooth
orbifold of dimension $\ind(u_0 ; \mathbf{c})$, with isotropy group
$\Aut(u_0)$ at $(j_0,u_0)$ and a natural isormorphism
$$
T_{(j_0,u_0)} \left( \uU / G \right) = \ker D\dbar_J(j_0,u_0) /
\aut(\dot{\Sigma},j_0).
$$
One can adapt the argument in
\cite{Dragnev} to show that charts constructed in this way are always
smoothly compatible.

To complete the proof, we show that $\uU / G$ is homeomorphic to a neighborhood
of $(\Sigma,j_0,\Gamma,u_0)$ in $\mM$.  
Clearly $\uU$ contains a representative of
every $J$-holomorphic curve near $(j_0,u_0)$, so the point is to show that any
two such curves $(j,u)$ and $(j',u')$ that are equivalent in $\mM^{\mathbf{c}}$
are related by the $G$--action.  

Suppose first that $\dot{\Sigma}$ is non-stable and is not $\AA$ or $\TT^2$:
then $\tT$ contains only $j_0$, and $(j_0,u) \sim (j_0,u')$ if and only if
$u' = u \circ \varphi$ for some $\varphi \in \Aut(\dot{\Sigma},j_0) = G$, so
we are done.  The case $\AA$ is hardly more complicated: now $\tT$ is
$1$--dimensional and $M(\dot{\Sigma})$ is trivial, so $\mM(\dot{\Sigma}) 
= \tT(\dot{\Sigma})$ and $j, j' \in \tT$ are equivalent in 
$\mM(\dot{\Sigma})$ if and only if $j = j'$.  Thus $(j,u) \sim (j',u')$
implies $j = j'$ and $u' = u \circ \varphi$ for some $\varphi \in
\Aut(\dot{\Sigma},j)$.  But our construction of $\tT_{\AA}$ identifies
$\Aut(\dot{\Sigma},j)$ with $\Aut(\dot{\Sigma},j_0) = G$, so again we are
done.

Consider now the stable cases and $\TT^2$, for which $M(\dot{\Sigma})$ 
is nontrivial.  For these, the groups 
$\Aut_0(\dot{\Sigma},j) := \Aut(\dot{\Sigma},j) \cap \Diff_0(\Sigma,\Gamma)$
for every $j \in \tT$ are identified with $\Aut_0(\dot{\Sigma},j_0)$;
this is a nontrivial statement only for $\dot{\Sigma} = \TT^2$, where our
explicit construction of $\tT = \tT_{\TT^2}$ identifies every
$\Aut_0(\TT^2,j)$ with $\TT^2$, acting by translations.
Meanwhile, for each $j \in \tT$ there is a finite subgroup $G_j \subset
\Aut(\dot{\Sigma},j)$ ($G_j = \Aut(\dot{\Sigma},j)$ in the stable 
cases) naturally isomorphic to the stabilizer of $[j]
\in \tT(\dot{\Sigma})$ under the $M(\dot{\Sigma})$--action, such that
$\Aut(\dot{\Sigma},j)$ is the semidirect product of $G_j$ with 
$\Aut_0(\dot{\Sigma},j_0)$.
Now if $(j,u)$ and $(j',u')$ are two elements of $\uU$ that
represent equivalent curves, so $j' = \psi^*j$ and $u' = u\circ \psi$
for some $\psi \in \Diff_+(\Sigma,\Gamma)$, we need to show that
$\psi \in G$.  In terms of the $M(\dot{\Sigma})$--action on
$\tT(\dot{\Sigma})$, $[\psi] \cdot [j] = [j']$ implies that if $j$ and
$j'$ are both sufficiently close to $j_0$, then $[\psi]$ belongs to the
stabilizer of $[j_0]$, i.e.~$[\psi] \cdot [j_0] = [j_0]$.  Thus there
is a unique $\varphi \in G_{j_0}$ such that $[\varphi] = [\psi] \in
M(\dot{\Sigma})$, and by construction, $\varphi^*j = j'$.  It follows
that $(\psi \circ \varphi^{-1})^*j = j$, so $\psi \circ \varphi^{-1} \in
\Aut_0(\dot{\Sigma},j) = \Aut_0(\dot{\Sigma},j_0)$, and
$\psi$ is thus a product of two maps in $\Aut(\dot{\Sigma},j_0)$.
\end{proof}

\subsection{The generalized normal bundle}
\label{subsec:normal}

For the remainder of this section, we shall consider a fixed
non-constant holomorphic curve $(j,u) \in \dbar_J^{-1}(0) 
\subset \tT \times \bB$ and examine the operator 
$D\dbar_J(j,u) = \mathbf{G}_u + \mathbf{D}_u$ more closely.  When we refer to
$\dot{\Sigma}$ or $\Sigma$ as a Riemann surface, we will always mean with
complex structure~$j$.  

As was observed in \cite{IvashkovichShevchishin},
the operator $\mathbf{D}_u$ defines a natural holomorphic structure on the
bundle $u^*TW \to \dot{\Sigma}$: indeed, the complex linear part of
$\mathbf{D}_u$ is also a Cauchy-Riemann type operator, so there is a
unique holomorphic structure whose local holomorphic sections vanish
under this operator.  This induces a holomorphic structure on
$\Hom_\CC(T\dot{\Sigma},u^*TW)$, and one can then show
(cf.~\cite{IvashkovichShevchishin}) that
$$
du \in \Gamma(\Hom_\CC(T\dot{\Sigma},u^*TW))
$$
is a holomorphic section.  Thus if $z_0$ is an interior critical point 
of~$u$, we can choose a holomorphic trivialization of 
$\Hom_\CC(T\dot{\Sigma},u^*TW)$ 
near $z_0$ and express $du$ as a $\CC^n$--valued function of the form
\begin{equation}
\label{eqn:criticalOrder}
(z - z_0)^k F(z)
\end{equation}
for some $k \in \NN$ and $\CC^n$--valued holomorphic function~$F$ with
$F(z_0) \ne 0$.
In this case we define the \emph{order} of the critical point by
$$
\ord(du ; z_0) = k.
$$
A similar definition is possible for
$z_0 \in \Crit(u) \cap \p\Sigma$ since $du$ satisfies the totally real
boundary condition $du(\p\Sigma) \subset \lL$, where
$$
\lL = \{ A \in \Hom_\CC(T\dot{\Sigma},u^*TW)|_{\p\Sigma} \ |\ 
A(T(\p\Sigma)) \subset \Lambda \}.
$$
Indeed, one can then choose a trivialization near $z_0$ such that
$du$ satisfies the Schwartz reflection principal, and define
$\ord(du ; z_0)$ again via \eqref{eqn:criticalOrder} after reflection.
Define the $\frac{1}{2}\ZZ$--valued algebraic count of critical points by
\begin{equation}
\label{eqn:totalCrit}
Z(du) = \sum_{z \in du^{-1}(0) \cap \interior{\dot{\Sigma}}} \ord(du ; z)
+ \frac{1}{2} \sum_{z \in du^{-1}(0) \cap \p\Sigma} \ord(du ; z).
\end{equation}

The expression \eqref{eqn:criticalOrder} has a second important purpose:
the complex subspace of $T_{u(z)}W$ spanned in the trivialization by
$F(z) \in \CC^n \setminus \{0\}$ allows us to define a smooth 
rank~$1$ subbundle
$$
T_u \subset u^*TW
$$
such that for any $z \in \dot{\Sigma}\setminus\Crit(u)$,
$(T_u)_z = \im du(z)$.  We will call this the \emph{generalized tangent
bundle} to $u$.

\begin{lemma}
The intersection $(T_u)_z \cap \Lambda_z$ is $1$--dimensional for all
$z \in\p\Sigma$.
\end{lemma}
\begin{proof}
It can never be $2$--dimensional since $(T_u)_z$ is a complex subspace and
$\Lambda_z$ is totally real.  Moreover it is clearly at least $1$--dimensional 
whenever $du(z) \ne 0$, as then $Tu(Y) \in T_{u(z)}L = \Lambda_z$ for any
$Y \in T_z\p\Sigma$.  Since critical points are isolated and the condition
$\dim (T_u)_z \cap \Lambda_z = 0$ is open, the result follows.
\end{proof}

By the lemma, we can define a totally real subbundle
$$
\ell^T = \Lambda \cap T_u|_{\p\Sigma} \subset T_u|_{\p\Sigma},
$$
and by construction $du$ now defines a section of the complex line bundle
$$
\Hom_\CC(T\dot{\Sigma},T_u) \to \dot{\Sigma}
$$
with totally real boundary condition $du(\p\Sigma) \subset \lL^T$, where
$$
\lL^T = \{ A \in \Hom_\CC(T\dot{\Sigma},T_u)|_{\p\Sigma}\ |\ 
A(T(\p\Sigma)) \subset \ell^T \}.
$$
As defined in the appendix, the algebraic count of zeroes for this section
is precisely $Z(du)$.

Observe that both $T\dot{\Sigma}$ and $T_u$ admit natural extensions over
the compactified surface $\overline{\Sigma}$; we define this extension
for $T_u$ via its natural identification with $T\dot{\Sigma}$ under $du$
since $u$ is immersed near infinity.  There is also a natural trivialization
$\tau$ of $T\dot{\Sigma}$ at infinity defined by the cylindrical coordinates
$(s,t) \in Z_\pm$, and we can define $\tau$ also over $\p\Sigma$ such that
the Maslov index $\mu^\tau(T\dot{\Sigma},T(\p\Sigma))$ vanishes.  Then
$$
c_1^\tau(T\dot{\Sigma}) = \chi(\dot{\Sigma}).
$$
Now choose any trivialization $\Phi$ of $T_u$ over $\p\Sigma$ and 
define it at infinity to be the same as~$\tau$.
The combination of $\tau$ and $\Phi$
induces a trivialization of $\Hom_\CC(T\dot{\Sigma},T_u)$ over $\p\Sigma$
and at infinity, which we will also denote by~$\Phi$.  Then we can apply
Prop.~\ref{prop:halfc1}, noting that the winding terms are zero by
construction, and obtain
\begin{equation}
\label{eqn:Zdu2}
Z(du) = c_1^\Phi(\Hom_\CC(T\dot{\Sigma},T_u)) +
\frac{1}{2}\mu^\Phi(\Hom_\CC(T\dot{\Sigma},T_u),\lL^T).
\end{equation}
To break this down further, note that the natural bundle isomorphism
$T\dot{\Sigma} \otimes \Hom_\CC(T\dot{\Sigma},T_u) \to T_u :
v \otimes A \mapsto Av$ sends $T(\p\Sigma) \otimes \lL^T$ to $\ell^T$, thus
$$
c_1^\Phi(T_u) = c_1^\tau(T\dot{\Sigma}) + c_1^\Phi(\Hom_\CC(T\dot{\Sigma},
T_u)) = \chi(\dot{\Sigma}) + c_1^\Phi(\Hom_\CC(T\dot{\Sigma},T_u)),
$$
and
\begin{equation*}
\begin{split}
\mu^\Phi(T_u,\ell^T) &= \mu^\tau(T\dot{\Sigma},T(\p\Sigma)) +
\mu^\Phi(\Hom_\CC(T\dot{\Sigma},T_u),\lL^T) \\
&= \mu^\Phi(\Hom_\CC(T\dot{\Sigma},T_u),\lL^T),
\end{split}
\end{equation*}
so \eqref{eqn:Zdu2} implies
\begin{equation}
\label{eqn:c1Tu}
c_1^\Phi(T_u) = \chi(\dot{\Sigma}) - \frac{1}{2}\mu^\Phi(T_u,\ell^T) +
Z(du).
\end{equation}

We next choose a \emph{generalized normal bundle} $N_u \to \dot{\Sigma}$,
which we define to be any rank~$n-1$ subbundle of $u^*TW$ such that
$$
u^*TW = T_u \oplus N_u,
$$
and the following conditions are satisfied:
\begin{enumerate}
\item
On the cylindrical neighborhoods $\uU_z$ for $z \in \Gamma^\pm$, 
$N_u$ matches the hyperplane distributions $\xi_\pm$, and thus extends to
infinity as $N_u|_{\delta_z} = \xi_\pm|_{\gamma_z}$.
\item
For $z \in \p\Sigma$, there is always a real $(n-1)$--dimensional intersection
$\ell_z^N := (N_u)_z \cap \Lambda_z$, thus defining a totally real
subbundle
$$
\ell^N \subset N_u|_{\p\Sigma}
$$
such that $\ell^T \oplus \ell^N = \Lambda$.
\end{enumerate}

\subsection{Splitting the linearization}
\label{subsec:splitting}

The splitting $u^*TW = T_u \oplus N_u$ defines projection maps
$\pi_T \in \Gamma(\Hom_\CC(u^*TW,T_u))$ and $\pi_N \in 
\Gamma(\Hom_\CC(u^*TW,N_u))$, both of which are smooth and satisfy
exponential decay conditions due to the asymptotic behavior of~$u$.
It follows that these define bounded linear projection operators
\begin{equation*}
\begin{split}
W^{1,p,\boldsymbol{\delta}}_\Lambda(u^*TW) \oplus
V_\Gamma \oplus X_\Gamma &\xrightarrow{\pi_T}
W^{1,p,\boldsymbol{\delta}}_{\ell^T}(T_u) \oplus V^T_\Gamma, \\
W^{1,p,\boldsymbol{\delta}}_\Lambda(u^*TW) \oplus
V_\Gamma \oplus X_\Gamma &\xrightarrow{\pi_N}
W^{1,p,\boldsymbol{\delta}}_{\ell^N}(N_u) \oplus X_\Gamma, \\
\end{split}
\end{equation*}
where $V_\Gamma^T \subset \Gamma(T_u)$ is the isomorphic image of
$V_\Gamma^\Sigma \subset \Gamma(T\dot{\Sigma})$ under the map 
$du : \Gamma(T\dot{\Sigma}) \to \Gamma(T_u) : v \mapsto Tu \circ v$,
and without loss of generality $X_\Gamma \in \Gamma(N_u)$.
There is thus a Banach space splitting
$$
W^{1,p,\boldsymbol{\delta}}_\Lambda(u^*TW) \oplus V_\Gamma \oplus X_\Gamma =
\left( W^{1,p,\boldsymbol{\delta}}_{\ell^T}(T_u) \oplus V_\Gamma^T \right) 
\oplus\left( W^{1,p,\boldsymbol{\delta}}_{\ell^N}(N_u) \oplus X_\Gamma \right),
$$
and a similar splitting
$$
L^{p,\boldsymbol{\delta}}(\overline{\Hom}_\CC(T\dot{\Sigma},u^*TW)) =
L^{p,\boldsymbol{\delta}}(\overline{\Hom}_\CC(T\dot{\Sigma},T_u)) \oplus
L^{p,\boldsymbol{\delta}}(\overline{\Hom}_\CC(T\dot{\Sigma},N_u)),
$$
so that with respect to these splittings, the operator
$$
\mathbf{D}_u : 
W^{1,p,\boldsymbol{\delta}}_\Lambda(u^*TW) \oplus V_\Gamma \oplus X_\Gamma
\to L^{p,\boldsymbol{\delta}}(\overline{\Hom}_\CC(T\dot{\Sigma},u^*TW))
$$ 
can be written in matrix form as
$$
\mathbf{D}_u = \begin{pmatrix}
\mathbf{D}_u^T & \mathbf{D}_u^{NT} \\
\mathbf{D}_u^{TN} & \mathbf{D}_u^N
\end{pmatrix}.
$$
It is trivial to show that 
$$\mathbf{D}_u^T : W^{1,p,\boldsymbol{\delta}}_{\ell^T}(T_u) \oplus V^T_\Gamma
\to L^{p,\boldsymbol{\delta}}(\overline{\Hom}_\CC(T\dot{\Sigma},T_u))
$$ 
and
$$\mathbf{D}_u^N : W^{1,p,\boldsymbol{\delta}}_{\ell^N}(N_u) \oplus X_\Gamma
\to L^{p,\boldsymbol{\delta}}(\overline{\Hom}_\CC(T\dot{\Sigma},N_u))
$$
each satisfy the appropriate Leibnitz rule for a Cauchy-Riemann type 
operator.  The former is asymptotic at each puncture $z \in \Gamma^\pm$
to the degenerate
asymptotic operator $-J_0 \frac{d}{dt}$ on a trivial complex line bundle; 
removing the exponential weight as in \S\ref{subsec:general},
this operator becomes $-J_0 \frac{d}{dt} \pm \boldsymbol{\delta}$, giving Conley-Zehnder
index $\mp 1$ with respect to the natural trivialization~$\tau$.
Thus the restriction of $\mathbf{D}_u^T$ to $W^{1,p,\boldsymbol{\delta}}_{\ell^T}(T_u)$
has index
$$
\chi(\dot{\Sigma}) + 2 c_1^\Phi(T_u) + \mu^\Phi(T_u,\ell^T) - \#\Gamma
$$
and adding the dimension of $V^T_\Gamma$ we find
\begin{equation}
\begin{split}
\label{eqn:indexTangent}
\ind(\mathbf{D}_u^T) &= \chi(\dot{\Sigma}) + 
2 c_1^\Phi(T_u) + \mu^\Phi(T_u,\ell^T) + \#\Gamma \\
&= 3\chi(\dot{\Sigma}) + \#\Gamma + 2 Z(du) \\
&= \dim\Aut(\dot{\Sigma},j) - \dim \tT(\dot{\Sigma}) + 2 Z(du),
\end{split}
\end{equation}
where the second line follows from \eqref{eqn:c1Tu}.

We call $\mathbf{D}_u^N$ the \emph{normal Cauchy-Riemann operator} at $u$.
It is also Fredholm; from the asymptotic identification of $N_u$ 
with $\xi_\pm$ along
orbits, we see that $\mathbf{D}_u^N$ is asymptotic to $\mathbf{A}_z$ at each
puncture $z \in \Gamma$.  We can use $\eqref{eqn:c1Tu}$ to relate its
index to $\ind(u ; \mathbf{c})$.  Abbreviate 
$\muCZ^\Phi(\gamma_\Gamma \pm \mathbf{c}_\Gamma) =
\sum_{z \in \Gamma^+} \muCZ^\Phi(\gamma_z + \mathbf{c}_z) 
- \sum_{z \in \Gamma^-} \muCZ^\Phi(\gamma_z - \mathbf{c}_z)$.  
Then removing the exponential weights as in \S\ref{subsec:general},
we apply the Riemann-Roch formula \eqref{eqn:RiemannRoch} and
repeat the calculation in \eqref{eqn:aCalculation} to find
\begin{equation}
\label{eqn:indexNormalCR}
\begin{split}
\ind(\mathbf{D}_u^N) &= (n-1)\chi(\dot{\Sigma}) + 2 c_1^\Phi(N_u) +
\mu^\Phi(N_u,\ell^N) \\
& \qquad + \sum_{z\in\Gamma^+} \muCZ^\Phi(\mathbf{A}_z + \boldsymbol{\delta}) -
\sum_{z\in\Gamma^-} \muCZ^\Phi(\mathbf{A}_z - \boldsymbol{\delta}) + \dim X_\Gamma \\
&= (n-1)\chi(\dot{\Sigma}) + 2 c_1^\Phi(N_u) +
\mu^\Phi(N_u,\ell^N) + \muCZ^\Phi(\gamma_\Gamma \pm \mathbf{c}_z) \\
&= (n - 1)\chi(\dot{\Sigma}) + 2\left[c_1^\Phi(u^*TW) - c_1^\Phi(T_u)\right] \\
&\qquad +
\left[ \mu^\Phi(u^*TW,\Lambda) - \mu^\Phi(T_u,\ell^T) \right] 
+ \muCZ^\Phi(\gamma_\Gamma \pm \mathbf{c}_\Gamma) \\
&= (n - 1)\chi(\dot{\Sigma}) + 2 c_1^\Phi(u^*TW) - 2\chi(\dot{\Sigma})
 - 2Z(du) + \mu^\Phi(u ; \mathbf{c}) \\
&= \ind(u ; \mathbf{c}) - 2Z(du).
\end{split}
\end{equation}

The main goal for this section is the following:

\begin{thm}
\label{thm:normalCR}
Assume $(\Sigma,j,\Gamma,u) \in \mM^\mathbf{c}$ is a non-constant curve
with Morse-Bott asymptotic orbits
and $\tT$ is any Teichm\"uller slice through $j$.
Then $\ker D\dbar_J(j,u)$ contains a subspace
$\ker(\mathbf{G}_u + \mathbf{D}_u^T) \subset T_j\tT \oplus
W^{1,p,\boldsymbol{\delta}}_{\ell^T}(T_u) \oplus V^T_\Gamma$
of dimension $2Z(du) + \dim\Aut(\dot{\Sigma},j)$ 
such that the normal projection induces a natural isomorphism
$$
\ker D\dbar_J(j,u) / \ker(\mathbf{G}_u + \mathbf{D}_u^T)
= \ker\mathbf{D}_u^N,
$$
and
$$
\im D\dbar_J(j,u) = L^{p,\boldsymbol{\delta}}(\overline{\Hom}_\CC(T\dot{\Sigma},
T_u)) \oplus \im\mathbf{D}_u^N .
$$
In particular, we have
\begin{equation*}
\begin{split}
\dim\ker D\dbar_J(j,u) &= 2Z(du) + \dim\Aut(\dot{\Sigma},j) + 
\dim\ker\mathbf{D}_u^N,\\
\dim\coker D\dbar_J(j,u) &= \dim\coker\mathbf{D}_u^N.
\end{split}
\end{equation*}
\end{thm}

\begin{cor}
\label{cor:regularity}
$(\Sigma,j,\Gamma,u) \in \mM^\mathbf{c}$ is regular if and only if 
$\mathbf{D}_u^N$ is surjective.
\end{cor}

The reason for this result is essentially that the analysis of the map
$(y,v) \mapsto \mathbf{G}_u y + \mathbf{D}_u v$ when $v$ is a section
of $T_u$ can be reduced to Lemma \ref{lemma:Teichmueller}, which one can
regard as an analytical statement about the smoothness of Teichm\"uller space.  
To achieve this reduction, we
introduce certain special Banach spaces of sections: 
for each $z_0 \in \Crit(u)$, choose holomorphic coordinates and
corresponding trivializations of $T\dot{\Sigma}$ and $T_u$ near $z_0$ so that
the bundle map $du : T\dot{\Sigma} \to T_u$ locally takes the form
$z \mapsto z^k$, where $k = \ord(du ; z_0)$.  Now for any function 
$d : \Crit(u) \to \ZZ$, define the Banach space
$$
W^{k,p,\boldsymbol{\delta},d}(T\dot{\Sigma})
$$
to consist of sections $v$ that are of class $W^{k,p}_{\text{loc}}$ on
$\dot{\Sigma}\setminus\Crit(u)$, class $W^{k,p,\boldsymbol{\delta}}$ near infinity,
and such that near each $z_0 \in \Crit(u)$, using the coordinates and
trivialization chosen above, the map
$$
z^{d(z_0)} v(z)
$$
is of class $W^{k,p}$.  Note that $v(z_0)$ may or may not be well defined:
if $d(z_0) > 0$ then $v$ is allowed to blow up at $z_0$, e.g.~it could
be meromorphic with a pole of order $\le d(z_0)$.
A suitable Banach space norm can be defined
using weighting functions supported near~$\Crit(u)$, and the subspace
$W^{k,p,\boldsymbol{\delta},d}_{T(\p\Sigma)}(T\dot{\Sigma})$ is defined by adding
the usual boundary condition; similarly we can define such spaces on
the bundles $T_u$, $\overline{\End}_\CC(T\dot{\Sigma})$
and $\overline{\Hom}_\CC(T\dot{\Sigma},T_u)$.
These are naturally isomorphic to our original Banach spaces
if $d(z) = 0$ for all $z \in \Crit(u)$.

The usefulness of this notion lies in the fact that if we choose
$d(z) := \ord(du ; z)$, then the correspondence $v \mapsto Tu \circ v$
defines Banach space isomorphisms
\begin{equation*}
\begin{split}
W^{1,p,\boldsymbol{\delta},d}_{T(\p\Sigma)}(T\dot{\Sigma}) &\xrightarrow{du}
W^{1,p,\boldsymbol{\delta}}_{\ell^T}(T_u), \\
V^\Sigma_\Gamma &\xrightarrow{du} V^T_\Gamma, \\
L^{p,\boldsymbol{\delta},d}(\overline{\End}_\CC(T\dot{\Sigma})) 
&\xrightarrow{du} L^{p,\boldsymbol{\delta}}(\overline{\Hom}_\CC(T\dot{\Sigma},T_u)).
\end{split}
\end{equation*}
We will stick with this choice of $d$ henceforward.

Using the fact that $z^k$ is holomorphic on the punctured disk for any
$k \in\ZZ$, it's easy to show that the natural linear Cauchy-Riemann operator
on $\Gamma(T\dot{\Sigma})$ defines a bounded linear map
$$
\mathbf{D}^\Sigma_d : 
W^{1,p,\boldsymbol{\delta},d}_{T(\p\Sigma)}(T\dot{\Sigma}) \oplus
V^\Sigma_\Gamma \to
L^{p,\boldsymbol{\delta},d}(\overline{\End}_\CC(T\dot{\Sigma})).
$$
The next result then follows from Lemma~\ref{lemma:DuOfTangents} by a
density argument.

\begin{lemma}
\label{lemma:conjugation}
For any $v \in W^{1,p,\boldsymbol{\delta},d}_{T(\p\Sigma)}(T\dot{\Sigma}) 
\oplus V^\Sigma_\Gamma$,
$\mathbf{D}_u (du(v)) = du( \mathbf{D}^\Sigma_d v)$.
\end{lemma}

\begin{lemma}
\label{lemma:surjective}
The operator 
\begin{equation*}
\begin{split}
\mathbf{L}_d : T_j \tT \oplus 
\left( W^{1,p,\boldsymbol{\delta},d}_{T(\p\Sigma)}(T\dot{\Sigma}) \oplus
V^\Sigma_\Gamma \right) &\to L^{p,\boldsymbol{\delta},d}
(\overline{\End}_\CC(T\dot{\Sigma})) \\
(y,v) &\mapsto j y + \mathbf{D}^\Sigma_d v
\end{split}
\end{equation*}
is surjective and has kernel of dimension 
$$
\dim\ker(\mathbf{L}_d) = 2Z(du) + \dim\Aut(\dot{\Sigma},j).
$$
\end{lemma}
\begin{proof}
We claim first that the result doesn't depend on the choice of Teichm\"uller
slice~$\tT$.  Indeed, in light of the splitting
$L^{p,\boldsymbol{\delta}}(\overline{\End}_\CC(T\dot{\Sigma})) =
\im(\mathbf{D}^\Sigma) \oplus T_j \tT$ and the natural
inclusion of this space in
$L^{p,\boldsymbol{\delta},d}(\overline{\End}_\CC(T\dot{\Sigma}))$,
an argument analogous to that in the proof of Lemma~\ref{lemma:allSlices}
shows that $\mathbf{L}_d$ has the same image as its natural extension
to
$L^{p,\boldsymbol{\delta}}(\overline{\End}_\CC(T\dot{\Sigma})) \oplus
\left( 
W^{1,p,\boldsymbol{\delta},d}_{T(\p\Sigma)}(T\dot{\Sigma}) \oplus
V^\Sigma_\Gamma \right)$.  We are thus free to change $\tT$: in
particular, we shall use Lemma~\ref{lemma:vanishAtCritical} to assume
in the following that all $y \in T_j\tT$ vanish on some fixed
neighborhood of $\Crit(u) \cup \Gamma$.

The subscript $d$ is meant to distinguish $\mathbf{L}_d$ and
$\mathbf{D}^\Sigma_d$ from the operators that appeared in
Lemma~\ref{lemma:Teichmueller}; we'll continue to denote the latter simply by
$$
\mathbf{L} : T_j \tT \oplus
\left( W^{1,p,\boldsymbol{\delta}}_{T(\p\Sigma)}(T\dot{\Sigma}) \oplus
V^\Sigma_\Gamma \right) \to L^{p,\boldsymbol{\delta}}
(\overline{\End}_\CC(T\dot{\Sigma})),
$$
with $\mathbf{D}^\Sigma$ denoting the restriction to
$W^{1,p,\boldsymbol{\delta}}_{T(\p\Sigma)}(T\dot{\Sigma}) \oplus V^\Sigma_\Gamma$.
The latter has index
$\ind(\mathbf{D}^\Sigma) = \dim\Aut(\dot{\Sigma},j) -\dim \tT$, 
whereas Lemma~\ref{lemma:conjugation}
implies that $\mathbf{D}^\Sigma_d$ is conjugate to $\mathbf{D}_u^T$, hence
\begin{equation*}
\begin{split}
\ind(\mathbf{D}^\Sigma_d) &= \ind(\mathbf{D}_u^T) =
\dim\Aut(\dot{\Sigma},j) - \dim \tT + 2Z(du) \\
&= \ind(\mathbf{D}^\Sigma) + 2Z(du)
\end{split}
\end{equation*}
and $\ind(\mathbf{L}_d) = \ind(\mathbf{L}) + 2Z(du) = 2Z(du) +
\dim\Aut(\dot{\Sigma},j)$.
The result will follow if we can show that $\dim\ker(\mathbf{L}_d) \le
\ind(\mathbf{L}) + 2Z(du)$.

To this end, define a $2Z(du)$--dimensional subspace $P \subset 
W^{1,p,\boldsymbol{\delta},d}_{T(\p\Sigma)}(T\dot{\Sigma})$ as follows: 
$P$ shall consist
of smooth sections on $\dot{\Sigma} \setminus \Crit(u)$, supported near
$\Crit(u)$, which in our chosen holomorphic trivializations near any
given $z_0 \in \Crit(u)$ take the form
\begin{equation}
\label{eqn:principal}
\frac{c_1}{z} + \frac{c_2}{z^2} + \ldots + \frac{c_{d(z_0)}}{z^{d(z_0)}}
\end{equation}
for $c_i \in \CC$ if $z_0 \in \interior{\Sigma}$, or $c_i \in \RR$
if $z_0 \in\p\Sigma$.  Since every section in $P$ is holomorphic near 
$\Crit(u)$, there is an obvious extension of $\mathbf{L}$,
$$
\mathbf{L}' : T_j \tT \oplus \left( W^{1,p,\boldsymbol{\delta}}_{T(\p\Sigma)}(T\dot{\Sigma})
\oplus V^\Sigma_\Gamma \oplus P \right) \to 
L^{p,\boldsymbol{\delta}}(\overline{\End}_\CC(T\dot{\Sigma})),
$$
which has $\ind(\mathbf{L}') = \ind(\mathbf{L}) + \dim P = 
\ind(\mathbf{L}) + 2Z(du) = \ind(\mathbf{L}_d)$.  Now since $\mathbf{L}$
is surjective by Lemma~\ref{lemma:Teichmueller}, so is $\mathbf{L}'$,
and thus $\dim\ker(\mathbf{L}') = \ind(\mathbf{L}) + 2Z(du)$.

To finish, we claim that $\ker(\mathbf{L}_d) \subset \ker(\mathbf{L}')$.
Indeed, suppose $y \in T_j \tT$ and $v \in W^{1,p,\boldsymbol{\delta},d}(T\dot{\Sigma})
\oplus V^\Sigma_\Gamma$ such that
$$
\mathbf{D}^\Sigma v = - j \circ y.
$$
Then by our assumption on $\tT$, $y$ vanishes near $\Crit(u)$ and
$v$ is therefore a holomorphic section in this neighborhood, except 
possibly at points of $\Crit(u)$.  Near $z_0 \in \Crit(u)$, the principal
part of $v$ in our holomorphic trivialization must have the form of 
\eqref{eqn:principal}: in particular
there cannot be an essential singularity or pole of order higher than
$d(z_0)$ since $z^{d(z_0)} v(z)$ is of class $W^{1,p}$.  There is thus
a unique section in $P$ that equals the principal part near
$\Crit(u)$, and subtracting this off we obtain a section in
$W^{1,p,\boldsymbol{\delta}}_{T(\p\Sigma)}(T\dot{\Sigma}) \oplus V^\Sigma_\Gamma$,
showing that $v$ belongs to the domain of $\mathbf{L}'$.
\end{proof}

\begin{proof}[Proof of Theorem~\ref{thm:normalCR}]
Any $v \in T_u\bB$ can be decomposed uniquely as $v = du(v_\Sigma) + v_N$ where
$v_\Sigma \in W^{1,p,\boldsymbol{\delta},d}_{T(\p\Sigma)}(T\dot{\Sigma}) \oplus
V^\Sigma_\Gamma$ and $v_N \in W^{1,p,\boldsymbol{\delta}}_{\ell^N}(N_u) \oplus X_\Gamma$.
Then for $y \in T_j\tT$, Lemma~\ref{lemma:conjugation} implies
\begin{equation*}
\begin{split}
D\dbar_J(j,u)(y,v) &= J\circ Tu \circ y + 
du(\mathbf{D}^\Sigma_d v_\Sigma) + \mathbf{D}_u^{NT} v_N +
\mathbf{D}_u^N v_N \\
&= du\left( j\circ y + \mathbf{D}^\Sigma_d v_\Sigma \right)
+ \mathbf{D}_u^{NT} v_N + \mathbf{D}_u^N v_N.
\end{split}
\end{equation*}
The first term is $\mathbf{G}_u y + \mathbf{D}_u^T du(v_\Sigma)$, and
by Lemma~\ref{lemma:surjective} this maps onto 
$L^{p,\boldsymbol{\delta}}(\overline{\Hom}_\CC(T\dot{\Sigma},T_u))$ with
$\dim\ker(\mathbf{G}_u + \mathbf{D}_u^T) = 2Z(du) + \dim\Aut(\dot{\Sigma},j)$.
The desired description of $\ker D\dbar_J(j,u)$ and $\im D\dbar_J(j,u)$ now 
follows easily from this expression since
$\mathbf{D}_u^{NT} v_N \in L^{p,\boldsymbol{\delta}}
(\overline{\Hom}_\CC(T\dot{\Sigma},T_u))$.
\end{proof}

\begin{example}
\label{ex:branchedCovers}
Though we've generally assumed $n \ge 2$, Theorem~\ref{thm:normalCR} also
applies to the case $n=1$: then the normal bundle has rank zero and
$D\dbar_J(j,u) = \mathbf{G}_u + \mathbf{D}_u^T$, so the theorem says
that $D\dbar_J(j,u)$ is a surjective operator of index
$2Z(du) + \dim\Aut(\dot{\Sigma},j)$.  One can apply this to understand
the moduli space $\mM(\dot{\Sigma},\dot{\Sigma}')$ of asymptotically 
cylindrical holomorphic maps
$$
\varphi : (\dot{\Sigma},j) \to (\dot{\Sigma}',j')
$$
between two punctured Riemann surfaces $\Sigma\setminus\Gamma$
and $\Sigma'\setminus\Gamma'$, up to equivalence by automorphisms on the
domain.  Such maps are equivalent to
holomorphic maps $(\Sigma,j) \to (\Sigma',j')$ that send $\Gamma$ to
$\Gamma'$.  Combining Theorem~\ref{thm:normalCR} and 
Theorem~\ref{thm:orbifoldFolk}, we see that for any
$\varphi \in \mM(\dot{\Sigma},\dot{\Sigma}')$, the connected component
$\mM_\varphi(\dot{\Sigma},\dot{\Sigma}')$ containing $\varphi$
is a smooth orbifold with
$$
\dim \mM_\varphi(\dot{\Sigma},\dot{\Sigma}') = 2Z(d\varphi),
$$
where of course the right hand side can be computed from the
Riemann-Hurwitz formula.
This fact is classical, but it will be useful in the proof of
Theorem~\ref{thm:orbifold} to view it
in our particular analytical setup.
\end{example}

Before restricting to the four-dimensional case, we mention one more simple 
application of Theorem~\ref{thm:normalCR}.  It gives namely an
upper bound on the algebraic number of critical points in terms of the
dimension of the kernel.  For somewhere injective curves in the generic
case this is simply the index, and we obtain:
\begin{cor}
\label{cor:criticalBound}
For generic $J$, all somewhere injective curves $u \in \mM$ satisfy
$$
2Z(du) \le \ind(u ; \mathbf{c}).
$$
\end{cor}
So for instance, if $\p\Sigma = \emptyset$ then 
somewhere injective curves of index~$0$ or~$1$ are necessarily
immersed for generic~$J$.  This is a simple version of the folk theorem that
generically, spaces of curves with at least a certain number of critical
points have positive codimension.

\subsection{The transversality criterion in dimension four}
\label{subsec:four}

We will now show that Theorem~\ref{thm:normalCR} implies
Theorem~\ref{thm:criterion} in the case $\dim W = 4$.  The key is
the fact that $N_u \to \dot{\Sigma}$ is now a complex line bundle,
so $\mathbf{D}_u^N$ will be subject to the constraints of 
Prop.~\ref{prop:transversality}.

Recall from
\eqref{eqn:normalChernIndex} the definition of the normal first Chern
number $c_N(u ; \mathbf{c})$.  An easy exercise combining the index formula 
with the
relations \eqref{eqn:CZwinding} between winding numbers and Conley-Zehnder 
indices yields the following alternative definition, reminiscent of
\eqref{eqn:c1adjusted}:

\begin{multline}
\label{eqn:cN}
c_N(u ; \mathbf{c}) = c_1^\Phi(u^*TW) - \chi(\dot{\Sigma}) + 
\frac{1}{2}\mu^\Phi(u^*TW,\Lambda) \\
+ \sum_{z \in \Gamma^+} \alpha_-^\Phi(\gamma_z + \mathbf{c}_z) -
\sum_{z \in \Gamma^-} \alpha_+^\Phi(\gamma_z - \mathbf{c}_z).
\end{multline}

\begin{prop}
\label{prop:c1cN}
If $u \in \mM^{\mathbf{c}}$ is not constant, then
$$
c_1(N_u,\ell^N,\mathbf{A}_\Gamma \pm \mathbf{c}_\Gamma) 
= c_N(u ; \mathbf{c}) - Z(du).
$$
\end{prop}
\begin{proof}
Choosing appropriate trivializations $\Phi$,
the relation follows by a simple calculation using the definitions 
\eqref{eqn:cN} and
\begin{multline*}
c_1(N_u,\ell^N,\mathbf{A}_\Gamma \pm \mathbf{c}_\Gamma) = c_1^\Phi(N_u) + 
\frac{1}{2} \mu^\Phi(N_u,\ell^N) \\
+ \sum_{z \in \Gamma^+} \alpha_-^\Phi(\mathbf{A}_z + \mathbf{c}_z) -
\sum_{z \in \Gamma^-} \alpha_+^\Phi(\mathbf{A}_z - \mathbf{c}_z)
\end{multline*}
and plugging in \eqref{eqn:c1Tu}.
\end{proof}

To finish the proof of Theorem~\ref{thm:criterion}, 
we relate $\mathbf{D}_u^N$ to a similar operator
on a larger weighted domain: for $z \in \Gamma$, regard the numbers
$\mathbf{c}_z = \pm\boldsymbol{\delta}$ now as exponential weights and,
recalling the notation for weighted Sobolev spaces from
\S\ref{sec:linear}, extend $\mathbf{D}_u^N$ to a new operator
$$
\widetilde{\mathbf{D}}_u^N : W^{1,p,\mathbf{c}_\Gamma}(N_u) \to
L^{p,\mathbf{c}_\Gamma}(\overline{\Hom}_\CC(T\dot{\Sigma},N_u)).
$$
The extended operator is conjugate to an operator on non-weighted spaces 
asymptotic to $\mathbf{A}_\Gamma \pm \mathbf{c}_\Gamma$, so
\eqref{eqn:RiemannRoch} gives $\ind(\widetilde{\mathbf{D}}_u^N) =
\ind(\mathbf{D}_u^N)$.  Moreover Prop.~\ref{prop:transversality} together
with Prop.~\ref{prop:c1cN} above implies
for $\widetilde{\mathbf{D}}_u^N$ precisely the transversality
criterion and kernel bound that we would desire for $\mathbf{D}_u^N$.
The result then follows because the domain of $\widetilde{\mathbf{D}}_u^N$
contains that of $\mathbf{D}_u^N$, hence
$\ker\mathbf{D}_u^N \subset \ker\widetilde{\mathbf{D}}_u^N$.

\section{Application to spaces of embedded curves}
\label{sec:application}

As an application of the transversality theory, we shall in this section
state and prove a stronger 
version of Theorem~\ref{thm:orbifold}.  For preparation, we review in
\S\ref{subsec:intersection} some
basic facts from the intersection theory of asymptotically cylindrical
holomorphic curves in four dimensions.  This theory has been developed
by R.~Siefring \cite{Siefring:intersection} for curves with fixed asymptotic
orbits, and is generalized to the Morse-Bott case in \cite{SiefringWendl}.
We expand on this in \S\ref{subsec:covering} by proving 
some useful formulas involving the intersection theory for multiple 
covers of orbits and holomorphic curves.
The proof of Theorem~\ref{thm:orbifold} then
appears in \S\ref{subsec:immersed}.  

For the remainder of this paper
we consider only pseudoholomorphic curves \emph{without boundary}.

\begin{notation}
In the following, we will often abbreviate the notation by printing
summations with $\pm$--signs in their index sets, e.g.
$$
\sum_{(z_1,z_2) \in \Gamma_1^\pm \times \Gamma_2^\pm} F_\pm(z_1,z_2).
$$
The intended meaning is then literally,
$$
\sum_{(z_1,z_2) \in \Gamma_1^+ \times \Gamma_2^+} F_+(z_1,z_2) +
\sum_{(z_1,z_2) \in \Gamma_1^- \times \Gamma_2^-} F_-(z_1,z_2).
$$
Several variations on this scheme will appear.
\end{notation}

\subsection{Intersection theory for punctured holomorphic curves}
\label{subsec:intersection}

This section will consist only of definitions and statements; we refer to
\cites{Siefring:intersection,SiefringWendl} for all proofs.

Throughout the following, $(W,J)$ is a $4$--dimensional almost complex 
manifold with cylindrical ends $(M_\pm,\hH_\pm)$, whose vector fields 
$X_\pm$ are Morse-Bott.
Suppose $u : \dot{\Sigma} \to W$ and $u' : \dot{\Sigma}' \to W$ are 
asymptotically cylindrical holomorphic curves belonging to
moduli spaces $\mM^{\mathbf{c}}$ and $\mM^{\mathbf{c}'}$ respectively
for some choices $\mathbf{c}$, $\mathbf{c}'$ 
of asymptotic constraints.  One of the goals of the intersection theory is to 
define an integer $i(u ; \mathbf{c} \ |\ u' ; \mathbf{c}')$ that is invariant 
as $u$ and $u'$ move continuously
through $\mM^{\mathbf{c}}$ and $\mM^{\mathbf{c}'}$
respectively, and can be interpreted as an algebraic intersection count
for the two curves.  One can show (see \cite{Siefring:asymptotics})
that if $u$ and $u'$ are geometrically
distinct, meaning they do not both cover the same somewhere injective curve,
then their intersections occur only within some compact subset, so the 
algebraic intersection count
$u \bullet u'$ is indeed finite and nonnegative.  It is \emph{not} however
homotopy invariant in general, as intersections can run out to infinity 
under homotopies.  There is nonetheless a well defined notion of an 
\emph{asymptotic intersection number}
$$
i_\infty(u ; \mathbf{c} \ |\ u' ; \mathbf{c}') \in \ZZ
$$
which is also nonnegative, such that the sum
\begin{equation}
\label{eqn:iAsSum}
i(u ; \mathbf{c} \ |\ u' ; \mathbf{c}') := u \bullet u' + 
i_\infty(u ; \mathbf{c} \ |\ u' ; \mathbf{c}')
\end{equation}
depends only on the respective connected components $\mM_{u}^{\mathbf{c}}$ 
and $\mM_{u'}^{\mathbf{c}'}$.  With some additional effort and 
(as yet unpublished)
analysis, one can show that $i_\infty(u ; \mathbf{c}\ | \ u' ; \mathbf{c}') = 0$ 
for generic somewhere 
injective curves and generic $J$: more precisely, the spaces of curves for 
which $i_\infty(u ; \mathbf{c}\ |\ u' ; \mathbf{c}') > 0$ have positive 
codimension, and so $u \bullet u'$ attains the maximal possible value 
$i(u ; \mathbf{c}\ |\ u' ; \mathbf{c}')$ generically.

It is useful to phrase the definition of $i(u ; \mathbf{c}\ |\ u' ;\mathbf{c}')$ 
in terms of the
\emph{relative} intersection number $u \bullet_\Phi u'$, where $\Phi$
is an arbitrary choice of trivialization for $\xi_\pm$ along the asymptotic
orbits of $u$ and~$u'$.  One computes $u \bullet_\Phi u'$ by counting the
intersections of $u'$ with a small perturbation of $u$ that is offset
in the $\Phi$--direction at infinity: the resulting integer is homotopy
invariant and depends on $\Phi$ up to homotopy.  Then as shown in
\cites{Siefring:intersection,SiefringWendl}, for each pair of orbits 
$\gamma, \gamma'$ and numbers $\boldsymbol{\epsilon},\boldsymbol{\epsilon}'
\in \RR$, there are integers $\Omega^\Phi_\pm(\gamma + \boldsymbol{\epsilon},
\gamma' + \boldsymbol{\epsilon}')$ such that
\begin{equation}
\label{eqn:i}
i(u ; \mathbf{c}\ |\ u' ; \mathbf{c}') = u \bullet_\Phi u' - 
\sum_{(z,z') \in \Gamma^\pm \times (\Gamma')^\pm} 
\Omega_\pm^\Phi(\gamma_z \pm \mathbf{c}_z,\gamma_{z'} \pm \mathbf{c}'_{z'}),
\end{equation}
with the dependence on $\Phi$ canceling out on the right hand side.
The actual definitions of 
$\Omega^\Phi_\pm(\gamma + \boldsymbol{\epsilon},\gamma' + 
\boldsymbol{\epsilon}')$ are as follows.  
We set $\Omega^\Phi_\pm(\gamma + \boldsymbol{\epsilon},\gamma' + 
\boldsymbol{\epsilon}') = 0$ if $\gamma$ and $\gamma'$ 
are geometrically distinct orbits, and for any simply covered orbit
$\gamma$ and $m , n \in \NN$, if $\gamma^m$ and $\gamma^n$ denote the
corresponding covers of $\gamma$, let
\begin{equation}
\label{eqn:Omega}
\Omega^\Phi_\pm(\gamma^m + \boldsymbol{\epsilon},\gamma^n
+ \boldsymbol{\epsilon}') = m n \cdot
\min\left\{ \frac{\mp \alpha_{\mp}^\Phi(\gamma^m + \boldsymbol{\epsilon})}{m}, 
\frac{\mp \alpha_{\mp}^\Phi(\gamma^n + \boldsymbol{\epsilon}')}{n} \right\}.
\end{equation}
We'll use the abbreviated notation $\Omega^\Phi_\pm(\gamma,\gamma')$ when
$\boldsymbol{\epsilon} = \boldsymbol{\epsilon}' = 0$.  Observe that
the right hand side of \eqref{eqn:i} makes sense even when $u$ and $u'$
are not geometrically distinct; in particular, we can use it to define
$i(u ; \mathbf{c}\ |\ u ; \mathbf{c})$, 
which is the appropriate generalization of a ``self-intersection number'' for 
punctured holomorphic curves.

If $u$ and $u'$ are geometrically distinct, then the
asymptotic contribution $i_\infty(u ; \mathbf{c}\ |\ u' ; \mathbf{c}')$
can be defined directly, thus giving a more conceptually
revealing definition of $i(u ; \mathbf{c}\ |\ u' ; \mathbf{c}')$ via
\eqref{eqn:iAsSum}.  Indeed, any pair of punctures for $u$ and $u'$ that
have the same sign and indistinct orbits offers a potential for intersections
to be ``hidden at infinity''.  For two such punctures $z \in \Gamma^\pm$
and $z' \in (\Gamma')^\pm$, denote by
$$
i_\infty^\Phi(u_z , u'_{z'})
$$
the \emph{relative asymptotic intersection}: this is computed by restricting
both curves to suitably small cylindrical neighborhoods of the
respective punctures and counting any
intersections that appear near infinity after perturbing $u$ in the 
$\Phi$--direction.  It turns out that whenever both curves are 
$J$--holomorphic, $i_\infty^\Phi(u_z, u'_{z'})$ is \textit{a priori}
bounded from below by $\Omega_\pm^\Phi(\gamma_z,\gamma_{z'})$: thus
the integer
$$
i_\infty(u_z , u'_{z'}) := i_\infty^\Phi(u_z , u'_{z'}) -
\Omega^\Phi_\pm(\gamma_z , \gamma_{z'})
$$
is nonnegative and independent of~$\Phi$.  Intuitively, it counts the
potential intersections of these two ends that can ``emerge from infinity''
under homotopies of $u$ and $u'$ that fix their asymptotic orbits.
Additional intersections may appear if either orbit is unconstrained and
allowed to move in a Morse-Bott family: the number of these is also nonnegative
and turns out to depend only on the orbits and their constraints.
Thus for any orbits $\gamma,\gamma'$ and numbers 
$\boldsymbol{\epsilon}, \boldsymbol{\epsilon}' \in \RR$, define
$$
i_{\MB}^\pm(\gamma + \boldsymbol{\epsilon} , \gamma' + \boldsymbol{\epsilon}') =
\Omega^\Phi_\pm(\gamma,\gamma') - 
\Omega^\Phi_\pm(\gamma + \boldsymbol{\epsilon}, 
\gamma' + \boldsymbol{\epsilon}').
$$
We will interpret $i_{\MB}^\pm(\gamma_z \pm \mathbf{c}_z ,
\gamma_{z'} \pm \mathbf{c}'_{z'})$ as the number of ``extra'' hidden
intersections not counted by $i_\infty(u_z , u'_{z'})$ that can emerge
as these two ends move generically according to their respective
constraints, potentially shifting the asymptotic orbits.  The total
asymptotic intersection number is then
\begin{equation}
\label{eqn:iInfty}
i_\infty(u ; \mathbf{c}\ |\ u' ; \mathbf{c}') := 
\sum_{(z,z') \in \Gamma^\pm \times (\Gamma')^\pm} \left[ 
i_\infty(u_z , u'_{z'}) + i_{\MB}^\pm(\gamma_z \pm \mathbf{c}_z ,
\gamma_{z'} \pm \mathbf{c}'_{z'}) \right].
\end{equation}
Each individual term in this sum is nonnegative, and can be expected to
vanish under generic perturbations of $u$ and $u'$ as ``potential''
intersections become real.

If $u : \dot{\Sigma} \to W$ is somewhere injective, we recall from 
\cite{MicallefWhite} that $u$ admits a local description near any critical 
point allowing one to 
define a nonnegative \emph{singularity index} $\delta(u)$: it gives an 
algebraic count of self-intersections $u(z) = u(z')$ for $z \ne z'$ after 
making local perturbations so that $u$ becomes immersed.  
As shown in \cite{Siefring:asymptotics}, this still makes
sense in the punctured case because $u$ is necessarily embedded outside of
some compact subset: then $\delta(u) \ge 0$, with equality if and only 
if $u$ is embedded.  
It is however possible for self-intersections to escape to infinity under 
homotopies, thus $\delta(u)$ is not homotopy invariant, but as with the
intersection number, one can add a nonnegative \emph{asymptotic singularity
index} $\delta_\infty(u ; \mathbf{c})$ so that the sum
\begin{equation}
\label{eqn:singu}
\sing(u ; \mathbf{c}) := \delta(u) + \delta_\infty(u ; \mathbf{c})
\end{equation}
depends only on the connected component $\mM_u^{\mathbf{c}}$, and equals 
$\delta(u)$ generically but not always.  The condition
$\sing(u ; \mathbf{c}) = 0$ is then necessary and sufficient so that all 
somewhere injective curves in $\mM_u^{\mathbf{c}}$ should be embedded for
generic~$J$; note 
that one still may have $u$ embedded if $\sing(u ; \mathbf{c}) > 0$, but 
then generic curves close to $u$ will not be.  The asymptotic contribution
is a sum of the form
\begin{multline}
\label{eqn:singMB}
2\delta_\infty(u ; \mathbf{c}) = \sum_{z \ne z' \in \Gamma^\pm}
\left[ i_\infty(u_z,u_{z'}) + i_{\MB}^\pm(\gamma_z \pm \mathbf{c}_z,
\gamma_{z'} \pm \mathbf{c}_{z'}) \right] \\+
\sum_{z \in \Gamma^\pm} \left[ 2\delta_\infty(u_z) + 
2\delta_{\MB}^\pm(\gamma_z \pm \mathbf{c}_z) \right],
\end{multline}
in which every term is nonnegative if $u$ is $J$--holomorphic.
We interpret $\delta_\infty(u_z)$ as the number of self-intersections near
the puncture $z$ that may emerge from infinity under generic homotopies
fixing the orbit $\gamma_z$; this can happen if $\gamma_z$ is
multiply covered, as distinct branches of the cylinder approaching
$\gamma_z$ run into each other under perturbation.  Define
$$
2\delta_\infty(u_z) = i_\infty^\Phi(u_z,u_z) - \Omega^\Phi_\pm(\gamma_z),
$$
where $\Omega^\Phi_\pm(\gamma) \in \ZZ$ is the ``self-intersection analogue'' 
of $\Omega^\Phi_\pm(\gamma,\gamma')$, giving a different theoretical minimum for
$i_\infty^\Phi(u_z,u_z)$ since the two ends are identical.
To write it down explicitly for the $k$--fold cover of a simple orbit $\gamma$,
choose any nontrivial eigenfunction $e_\mp$ of $\mathbf{A}_{\gamma^k}$ whose
winding about $\gamma^k$ equals $\alpha_\mp^\Phi(\gamma^k)$, and note that its
covering number $\cov(e_\mp) \in \NN$ depends on
$k$ and $\alpha_\mp^\Phi(\gamma^k)$ but not on the choice $e_\mp$
(cf.~Lemma~\ref{lemma:eigenCover}).  Thus we denote
$$
\cov_\mp(\gamma^k) := \cov(e_\mp),
$$
and then define
\begin{equation}
\label{eqn:Omega0}
\Omega^\Phi_\pm(\gamma^k) = \mp (k-1) \alpha_\mp^\Phi(\gamma^k) + 
\left[ \cov_\mp(\gamma^k) - 1 \right].
\end{equation}

Similarly, $\delta_{\MB}^\pm(\gamma_z \pm \mathbf{c}_z)$ counts further
self-intersections that may emerge if $\gamma_z$ is allowed to move in a
Morse-Bott family.  This doesn't happen if every orbit in the family has
the same minimal period, but if $\gamma_z$ converges to an orbit with smaller
minimal period (and thus higher covering number), the existence of
additional branches can hide extra intersections at infinity.  The following
characterization of Morse-Bott manifolds will be useful.

\begin{prop}
\label{prop:MBcharacterization}
If $M$ is a $3$--manifold with a Morse-Bott vector field $X$,
then every Morse-Bott submanifold  $P \subset M$ can be
described as follows.  There exists a number $\tau > 0$ such that all but a
discrete set of orbits in $P$ have minimal period~$\tau$; we shall call these
\emph{generic} orbits.  The other orbits will be called \emph{exceptional}:
any such orbit with period~$\tau$ is an $m$--fold cover of a simply
covered orbit $\gamma$ for some $m \ge 2$ (called the
\emph{isotropy}), and $\gamma^k$ is nondegenerate for all
$k \in \{1,\ldots,m-1\}$.  The isotropy of an exceptional orbit 
is always~$2$ if $\dim P = 2$.
\end{prop}

Now, define $\delta_{\MB}^\pm(\gamma \pm \boldsymbol{\delta}) = 0$
if $\boldsymbol{\delta} > 0$; recall this case is associated with
a constraint that fixes $\gamma_z$, thus there can be no ``extra''
self-intersections appearing due to Morse-Bott considerations.  The
definition is as follows if $\boldsymbol{\delta} < 0$: given an orbit
$\gamma$, set $\gamma_\epsilon = \gamma$ if it's nondegenerate, otherwise
let $\gamma_\epsilon$ denote any nearby \emph{generic} orbit in the
same Morse-Bott family as~$\gamma$.  If $\gamma_\epsilon$ is
simply covered, $k \in \NN$, and $\gamma$ has isotropy $m \in \NN$, then set
\begin{equation}
\label{eqn:deltaMB}
2\delta_{\MB}^\pm(\gamma^k \pm \boldsymbol{\delta}) = k(m-1) \nu_\mp(\gamma^k) 
+ \cov_\mp(\gamma^k) - \cov_\mp(\gamma_\epsilon^k),
\end{equation}
where $\nu_\mp(\gamma) \in \{0,1\}$ is defined in \eqref{eqn:nu}. 
Observe that the two covering terms refer to homotopic
eigenfunctions $e$ of
$\mathbf{A}_{\gamma^k}$ and $e_\epsilon$ of $\mathbf{A}_{\gamma_\epsilon^k}$,
so if $e_\epsilon$ is an $n$--fold cover then
$e$ is as well, hence $\cov_\mp(\gamma^k) \ge \cov_\mp(\gamma_\epsilon^k)$.  
The inequality may
sometimes be strict, because $e$ is attached to a $km$--covered orbit, while
the orbit of $e_\epsilon$ is only $k$--covered.  In any case, clearly
$\delta_{\MB}^\pm(\gamma^k \pm \boldsymbol{\delta}) \ge 0$
in general, and it vanishes whenever $\gamma$ is a generic orbit.

For the curve $u \in \mM^{\mathbf{c}}$, choose for each $z \in \Gamma$ a
generic perturbation $\gamma_z^\epsilon$ of the orbit $\gamma_z$,
setting $\gamma_z^\epsilon = \gamma_z$ if either $\gamma_z$ is nondegenerate
or $z \in \Gamma_C$.  Then let
$$
\cov_\infty(\gamma_\Gamma ; \mathbf{c}) = 
\sum_{z \in \Gamma^\pm} \left[\cov_\mp(\gamma^\epsilon_z) - 1 \right],
$$
and
$$
\cov_{\MB}(\gamma_\Gamma ; \mathbf{c}) =
\sum_{z \in \Gamma^\pm_U} \left[\cov(\gamma_z^\epsilon) - 1\right] \cdot
\nu_\mp(\gamma_z),
$$
with $\cov(\gamma_z^\epsilon)$ denoting the covering number of
$\gamma_z^\epsilon$.

\begin{thmu}[Adjunction formula \cites{Siefring:intersection,
SiefringWendl}]
For any somewhere injective curve $u \in \mM^{\mathbf{c}}$,
$$
i(u ; \mathbf{c}\ |\ u ; \mathbf{c}) = 2\sing(u ; \mathbf{c}) + 
c_N(u ; \mathbf{c}) +
\cov_\infty(\gamma_\Gamma ; \mathbf{c}) + 
\cov_{\MB}(\gamma_\Gamma ; \mathbf{c}).
$$
\end{thmu}

\subsection{Some covering relations}
\label{subsec:covering}

It will be useful to have formulas relating the intersection invariants
of holomorphic curves and their multiple covers.  A prerequisite for this is
to have corresponding covering formulas for periodic orbits, so to start with,
assume $M$ is a $3$--manifold with stable Hamiltonian structure
$\hH = (X,\xi,\omega,J)$.
Given an orbit $\gamma$, we shall denote the corresponding asymptotic
operator by $\mathbf{A}_\gamma$ and the $k$--fold cover of $\gamma$ by
$\gamma^k$.  Then if $\mathbf{A}_\gamma e = \lambda e$, the eigenfunction 
has a $k$--fold cover $e^k$ such that $\mathbf{A}_{\gamma^k} e^k = k\lambda
e^k$.  In general, we say that an eigenfunction $f$ of 
$\mathbf{A}_{\gamma^k}$ \emph{is a $k$--fold cover} if there exists an 
eigenfunction $e$ of $\mathbf{A}_\gamma$ such that $f = e^k$.

In the following, whenever a trivialization $\Phi$ along an orbit $\gamma$
appears, we will use the same notation $\Phi$ to denote the resulting
induced trivializations along all covers of $\gamma$.

\begin{lemma}[\cite{Wendl:compactnessRinvt}*{Lemma~3.5}]
\label{lemma:eigenCover}
Suppose $\Phi$ is a trivialization along $\gamma$. Then
a nontrivial eigenfunction $e$ of $\mathbf{A}_{\gamma^k}$ is a $k$--fold
cover if and only if $\wind^\Phi(e) \in k \ZZ$.
\end{lemma}

\begin{lemma}
\label{lemma:coversAreEven}
Suppose $\gamma$ is a periodic orbit of $X$ and $\boldsymbol{\epsilon} \in \RR$.
If $\mathbf{A}_\gamma + \boldsymbol{\epsilon}$ is nondegenerate and 
$p(\gamma + \boldsymbol{\epsilon}) = 0$, then $\mathbf{A}_{\gamma^k} + k
\boldsymbol{\epsilon}$ is nondegenerate and
$p(\gamma^k + k\boldsymbol{\epsilon}) = 0$ for all $k \in \NN$.
\end{lemma}
\begin{proof}
If $p(\gamma + \boldsymbol{\epsilon}) = 0$, then $\sigma(\mathbf{A}_\gamma +
\boldsymbol{\epsilon})$ contains a pair of neighboring eigenvalues with
opposite signs and eigenfunctions of the same winding number.  The $k$--fold
covers of these are eigenfunctions of $\mathbf{A}_{\gamma^k} 
+ k\boldsymbol{\epsilon}$ with the same properties, thus
$\mathbf{A}_{\gamma^k} + k\boldsymbol{\epsilon}$ is nondegenerate and has
even parity.
\end{proof}
\begin{remark}
\label{remark:dropTheK}
If $\boldsymbol{\epsilon} \in \RR$ is sufficiently close but not equal to 
zero, then we may
always assume that for all $k \in \NN$ up to some arbitrarily large (but finite)
bound, $\mathbf{A}_{\gamma^k} + \boldsymbol{\epsilon}$ is nondegenerate and
$\alpha^\Phi_\pm(\gamma^k + k\boldsymbol{\epsilon}) = \alpha^\Phi_\pm(\gamma^k
+ \boldsymbol{\epsilon})$.  We can thus replace $k\boldsymbol{\epsilon}$ with
$\boldsymbol{\epsilon}$ in the statement above whenever $\boldsymbol{\epsilon}$
is assumed close to zero, and the same applies to several statements below.
\end{remark}

\begin{cor}
\label{cor:exceptionalsAreOdd}
For any exceptional orbit in a Morse-Bott family, the underlying
simple orbit and all of its nondegenerate covers are odd.
\end{cor}

\begin{prop}
For any periodic orbit $\gamma$ of $X$, $k \in \NN$ and $\boldsymbol{\epsilon} \in \RR$, 
there exist
integers $q_\pm(\gamma + \boldsymbol{\epsilon} ; k) \in \{0,\ldots,k-1\}$ such that
\begin{equation}
\label{eqn:q}
\alpha_\pm^\Phi(\gamma^k + k\boldsymbol{\epsilon}) = 
k \alpha_\pm^\Phi(\gamma + \boldsymbol{\epsilon}) \mp
q_\pm(\gamma + \boldsymbol{\epsilon} ; k).
\end{equation}
\end{prop}
\begin{proof}
The integer $q_\pm(\gamma + \boldsymbol{\epsilon} ; k) :=
\mp \left[ \alpha^\Phi_\pm(\gamma^k + \boldsymbol{\epsilon}) - k 
\alpha^\Phi_\pm(\gamma + \boldsymbol{\epsilon}) \right]$
is well defined
after observing that all dependence on $\Phi$ in the right hand side
cancels, so it remains only to show that this number is between $0$ and~$k-1$.
Consider first the case $\boldsymbol{\epsilon}=0$, and choose a trivialization $\Phi_0$ 
along $\gamma$ such that $\alpha^{\Phi_0}_-(\gamma) = 0$.  Then there exists an
eigenfunction $e_-$ of $\mathbf{A}_\gamma$ with negative eigenvalue and
$\wind^{\Phi_0}(e_-) = 0$, and another eigenfunction $e_+$ with nonnegative 
eigenvalue and $\wind^{\Phi_0}(e_+) = 1$; moreover there are no eigenfunctions with
eigenvalue strictly between that of $e_-$ and~$0$.  Moving to the 
$k$--fold cover,
we obtain eigenfunctions $e_-^k$ and
$e_+^k$ of $\mathbf{A}_{\gamma^k}$ with
\begin{equation*}
\begin{split}
\wind^{\Phi_0}(e_-^k) = 0 & \qquad \text{ eigenvalue $< 0$,} \\
\wind^{\Phi_0}(e_+^k) = k & \qquad \text{ eigenvalue $\ge 0$,}
\end{split}
\end{equation*}
and there is no $k$--fold covered 
eigenfunction with eigenvalue strictly between that of $e_-^k$ and~$0$.
Then by Lemma~\ref{lemma:eigenCover}, this range of the spectrum of
$\mathbf{A}_{\gamma^k}$ contains no eigenfunctions with winding~$k$.
Since the winding depends monotonically on the eigenvalue, this implies
$\alpha_-^{\Phi_0}(\gamma^k) \in \{0,\ldots,k-1\}$.
An analogous argument gives the corresponding result for $\alpha_+$.
Finally if $\boldsymbol{\epsilon} \ne 0$, 
the arguments above give the same relation between 
the eigenfunctions of $\mathbf{A}_\gamma + \boldsymbol{\epsilon}$ and 
$\mathbf{A}_{\gamma^k} + k\boldsymbol{\epsilon}$.
\end{proof}

In preparation for the next lemma, for any orbit $\gamma$, numbers
$m,n,k \in \NN$ and $\boldsymbol{\delta},\boldsymbol{\epsilon} \in \RR$,
define the nonnegative integers
\begin{multline}
\label{eqn:qtilde}
\tilde{q}_\pm(\gamma^m + \boldsymbol{\delta},\gamma^n + \boldsymbol{\epsilon}; k) = k m n \cdot
\min\left\{ \frac{\mp \alpha_\mp^\Phi(\gamma^m + \boldsymbol{\delta})}{m}, 
\frac{\mp\alpha_\mp(\gamma^n + \boldsymbol{\epsilon})}{n} \right\} \\
- k m n \cdot \min\left\{ \frac{\mp \alpha_\mp^\Phi(\gamma^m + \boldsymbol{\delta})}{m} - 
\frac{q_\mp(\gamma^m + \boldsymbol{\delta} ; k)}{km}, 
\frac{\mp\alpha_\mp(\gamma^n + \boldsymbol{\epsilon})}{n} \right\}.
\end{multline}
Then a simple computation using the definitions of
$\Omega^\Phi_\pm(\gamma + \boldsymbol{\epsilon},\gamma' + \boldsymbol{\epsilon}')$
and $q_\pm(\gamma + \boldsymbol{\epsilon} ; k)$ implies:
\begin{lemma}
\label{lemma:OmegaCover}
For any simply covered orbit $\gamma$, $m,n,k \in \NN$ and $\boldsymbol{\delta},
\boldsymbol{\epsilon} \in \RR$,
$$
\Omega^\Phi_\pm(\gamma^{km} + k\boldsymbol{\delta},\gamma^n + \boldsymbol{\epsilon}) =
k \cdot \Omega^\Phi_\pm(\gamma^m + \boldsymbol{\delta},\gamma^n + \boldsymbol{\epsilon}) 
- \tilde{q}_\pm(\gamma^m + \boldsymbol{\delta},\gamma^n + \boldsymbol{\epsilon} ; k).
$$
\end{lemma}

Returning now to the context of a $4$--manifold $W$ with Morse-Bott 
cylindrical ends $(M_\pm,\hH_\pm)$, let us fix the following notation: 
$u \in \mM^{\mathbf{c}}$ is a holomorphic curve with domain 
$(\Sigma\setminus\Gamma,j)$, $(\Sigma,j)$ and 
$(\widetilde{\Sigma},\tilde{\jmath})$
are closed Riemann surfaces, and $\varphi : (\widetilde{\Sigma},\tilde{\jmath})
\to (\Sigma,j)$ is a 
holomorphic branched cover of degree~$\deg(\varphi) \in \NN$.  This restricts 
to a branched cover of punctured surfaces
$$
\dot{\varphi} : \widetilde{\Sigma} \setminus \widetilde{\Gamma}
\to \Sigma \setminus\Gamma,
$$
where $\widetilde{\Gamma} := \varphi^{-1}(\Gamma)$, and there is a resulting 
holomorphic curve $u \circ \varphi : \widetilde{\Sigma} \setminus
\widetilde{\Gamma} \to W$.  Its asymptotic orbits
are related to those of $u$ by
$$
\gamma_z = \gamma_{\varphi(z)}^{k_z}
$$
at each $z \in \widetilde{\Gamma}$, where $k_z := \ord(d\varphi ; z) + 1$, so that
$\varphi$ is $k_z$--to--$1$ near~$z$.  The constraints $\mathbf{c}$ on
$\Gamma$ can then be pulled back to constraints $\varphi^*\mathbf{c}$ on
$\widetilde{\Gamma}$ like so:
for any $\zeta \in \Gamma$ constrained to the orbit $\gamma_{\zeta}$,
define $\varphi^*\mathbf{c}$ by fixing the orbit $\gamma_{\zeta}^{k_z}$ at each
$z \in \varphi^{-1}(\zeta)$.  Then $u\circ\varphi 
\in \mM^{\varphi^*\mathbf{c}}$.

\begin{prop}
\label{prop:cNcover}
For $u \in \mM^\mathbf{c}$ and the cover
$u \circ \varphi \in \mM^{\varphi^*\mathbf{c}}$ defined above,
$$
c_N(u \circ \varphi ; \varphi^*\mathbf{c}) = \deg(\varphi) \cdot 
c_N(u ; \mathbf{c}) + Z(d\dot{\varphi}) + Q
$$
where $Q$ is a nonnegative integer.  Specifically,
$$
Q = \sum_{z \in \widetilde{\Gamma}^\pm} q_\mp(\gamma_{\varphi(z)} 
\pm \mathbf{c}_{\varphi(z)}; k_z).
$$
\end{prop}
\begin{proof}
Denote $\tilde{u} := u\circ \varphi$, $\tilde{\mathbf{c}} := 
\varphi^*\mathbf{c}$ and observe that for each $\zeta \in \Gamma$,
$\sum_{z \in \varphi^{-1}(\zeta)} k_z = k := \deg(\varphi)$.  Note also that
by extending $\dot{\varphi}$ to the circle compactifications of
$\widetilde{\Sigma}\setminus\widetilde{\Gamma}$ and 
$\Sigma\setminus\Gamma$, one can apply the Riemann-Hurwitz formula and obtain
$$
Z(d\dot{\varphi}) = -\chi(\widetilde{\Sigma}\setminus\widetilde{\Gamma}) +
k \chi(\Sigma\setminus\Gamma).
$$
Then using \eqref{eqn:q} and Remark~\ref{remark:dropTheK},
\begin{equation*}
\begin{split}
c_N(\tilde{u} ; \tilde{\mathbf{c}}) &= c_1^\Phi(\varphi^*u^*TW) - 
\chi(\widetilde{\Sigma} \setminus\widetilde{\Gamma}) 
+ \sum_{\zeta \in \Gamma^\pm} \sum_{z \in \varphi^{-1}(\zeta)}
\pm \alpha^\Phi_\mp(\gamma_{\zeta}^{k_z} \pm \tilde{\mathbf{c}}_z) \\
&= k c_1^\Phi(u^*TW) - k\chi(\Sigma \setminus\Gamma) + Z(d\dot{\varphi})
 + k \sum_{\zeta \in \Gamma^\pm} 
\pm \alpha^\Phi_\mp(\gamma_{\zeta} \pm \mathbf{c}_{\zeta}) \\
&\qquad + \sum_{\zeta \in \Gamma^\pm}
 \sum_{z \in \varphi^{-1}(\zeta)} \pm \left( 
  \alpha^\Phi_\mp(\gamma_{\zeta}^{k_z} \pm \tilde{\mathbf{c}}_z) 
- k_z \alpha^\Phi_\mp(\gamma_{\zeta} \pm \mathbf{c}_\zeta)\right) \\
&= k \cdot c_N(u ; \mathbf{c}) + Z(d\dot{\varphi}) 
+ \sum_{z \in \widetilde{\Gamma}^\pm}
q_\mp( \gamma_{\varphi(z)} \pm \mathbf{c}_{\varphi(z)} ; k_z).
\end{split}
\end{equation*}
\end{proof}

\begin{prop}
\label{prop:iCovering}
For the cover $u \circ \varphi \in \mM^{\varphi^*\mathbf{c}}$ as in 
Prop.~\ref{prop:cNcover} and any other curve $v \in \mM^{\mathbf{c}'}$,
$$
i(u \circ \varphi ; \varphi^*\mathbf{c}\ |\  v ; \mathbf{c}')
\ge 
\deg(\varphi) \cdot i(u ; \mathbf{c}\ |\ v ; \mathbf{c}').
$$
\end{prop}
\begin{proof}
Again denote $k := \deg(\varphi)$, $k_z := \ord(d\varphi ; z) + 1 \in \NN$ 
for each $z \in \widetilde{\Gamma}$, $\tilde{u} := u \circ \varphi$ and
$\tilde{\mathbf{c}} := \varphi^*\mathbf{c}$.
The relative intersection number satisfies
$\tilde{u} \bullet_\Phi v = k (u \bullet_\Phi v)$.  
Writing the puncture set of $v$ as $\Gamma'$, we apply
Lemma~\ref{lemma:OmegaCover} with Remark~\ref{remark:dropTheK} in mind and find
\begin{equation*}
\begin{split}
i(\tilde{u} &; \tilde{\mathbf{c}} \ |\ v ; \mathbf{c}') = 
\tilde{u} \bullet_\Phi v - \sum_{(z,z') \in \widetilde{\Gamma}^\pm \times 
(\Gamma')^\pm}
\Omega^\Phi_\pm(\gamma_z \pm \tilde{\mathbf{c}}_z,\gamma_{z'} \pm \mathbf{c}'_{z'}) \\
&= k \cdot ( u \bullet_\Phi v) - \sum_{(\zeta,z') \in \Gamma^\pm \times 
(\Gamma')^\pm}
\left( \sum_{z \in \varphi^{-1}(\zeta)} \Omega^\Phi_\pm(\gamma_\zeta^{k_z}
 \pm \tilde{\mathbf{c}}_z, \gamma_{z'} \pm \mathbf{c}'_{z'}) \right) \\
&= k \cdot ( u \bullet_\Phi v) - \sum_{(\zeta,z') \in \Gamma^\pm \times 
(\Gamma')^\pm}
\Bigg( \sum_{z \in \varphi^{-1}(\zeta)} 
\Big[ k_z \Omega^\Phi_\pm(\gamma_\zeta \pm \mathbf{c}_\zeta,\gamma_{z'} \pm 
\mathbf{c}'_{z'}) \\
 &\qquad - \tilde{q}_\pm(\gamma_\zeta \pm \mathbf{c}_{\zeta},\gamma_{z'} \pm 
\mathbf{c}'_{z'} ; k_z) 
\Big] \Bigg) \\
&= k \cdot ( u \bullet_\Phi v) - k \sum_{(\zeta,z') \in \Gamma^\pm\times 
(\Gamma')^\pm}
\Omega^\Phi_\pm(\gamma_\zeta \pm \mathbf{c}_\zeta\ |\ \gamma_{z'} \pm 
\mathbf{c}'_{z'}) \\
&\qquad + \sum_{(z,z') \in \widetilde{\Gamma}^\pm \times (\Gamma')^\pm} 
\tilde{q}_\pm(\gamma_{\varphi(z)} \pm \mathbf{c}_{\varphi(z)} \ |\ 
\gamma_{z'} \pm \mathbf{c}'_{z'} \ |\  k_z) \\
&= k \cdot i(u ; \mathbf{c} \ |\ v ; \mathbf{c}') +
\sum_{(z,z') \in \widetilde{\Gamma}^\pm\times (\Gamma')^\pm} 
\tilde{q}_\pm(\gamma_{\varphi(z)} \pm\tilde{\mathbf{c}}_{\varphi(z)} 
\ |\ \gamma_{z'} \pm\mathbf{c}'_{z'} \ |\ k_z).
\end{split}
\end{equation*}
The last term is nonnegative.
\end{proof}

\begin{defn}
\label{defn:partialOrder}
For a given set of punctures $\Gamma$, the set of all choices of asymptotic
constraints on $\Gamma$ admits a partial order defined as follows.  We say
$\mathbf{c}^- \le \mathbf{c}^+$ if for every $z \in \Gamma$ at which the
asymptotic orbit $\gamma_z$ is constrained by $\mathbf{c}^-$, it is also
constrained by $\mathbf{c}^+$ to the same orbit.
\end{defn}

Observe that if $\mathbf{c}^- \le \mathbf{c}+$, then $\mM^{\mathbf{c}+} \subset
\mM^{\mathbf{c}^-}$ and $\mathbf{c}^-_z \le \mathbf{c}^+_z$ for each 
$z \in \Gamma$.  One expects in general that weaker constraints should lead 
to larger intersection numbers, as intersections can more easily emerge from 
infinity under more general homotopies.  Indeed, using $\mathbf{c}_z^- \le
\mathbf{c}_z^+$ together with the fact that 
$\alpha_\mp^\Phi(\gamma + \boldsymbol{\epsilon})$ always has monotone
decreasing dependence on $\boldsymbol{\epsilon}$, we easily derive 
the following:

\begin{prop}
\label{prop:partialOrder}
If $\mathbf{c}^- \le \mathbf{c}^+$ and $u \in \mM^{\mathbf{c}^+}$, then
$$
c_N(u ; \mathbf{c}^-) \ge c_N(u ; \mathbf{c}^+).
$$
Moreover for any other curve $v \in \mM^{\mathbf{c}}$,
$$
i(u ; \mathbf{c}^- \ |\ v ; \mathbf{c}) \ge 
i(u ; \mathbf{c}^+ \ |\ v ; \mathbf{c}).
$$
\end{prop}

\subsection{Multiply covered limits are immersed}
\label{subsec:immersed}

We shall now state and prove a parametrized version of
Theorem~\ref{thm:orbifold}.

\begin{defn}
\label{defn:nice}
We will say that $u \in \mM^\mathbf{c}$ is a \emph{stable, nicely embedded} 
curve (with respect to the constraints $\mathbf{c}$) 
if it is somewhere injective and satisfies the following relations:
\begin{enumerate}
\item $i(u ; \mathbf{c} \ |\ u ; \mathbf{c}) \le 0$,
\item $\ind(u ; \mathbf{c}) \ge 0$,
\item $\ind(u ; \mathbf{c}) > c_N(u ; \mathbf{c})$.
\end{enumerate}
\end{defn}
Before going further, let us consider the properties of such curves and the
motivation for the definition.  Observe first that the combination of
$\ind(u ; \mathbf{c}) \ge 0$ and the relation
\begin{equation}
\label{eqn:cNindexAgain}
2 c_N(u ; \mathbf{c}) = \ind(u ; \mathbf{c}) - 2 + 2g + \#\Gamma_0(\mathbf{c})
\end{equation}
gives the lower bound $c_N(u ; \mathbf{c}) \ge -1$.  
Then the adjunction formula 
together with $i(u ; \mathbf{c}\ |\ u ; \mathbf{c}) \le 0$ 
implies $\sing(u ; \mathbf{c}) = 0$, so every somewhere 
injective curve in $\mM_u^\mathbf{c}$ is embedded.  
We can also deduce from the adjunction
formula that $c_N(u ; \mathbf{c}) \le 0$, 
and then \eqref{eqn:cNindexAgain} implies $\ind(u ; \mathbf{c}) \le 2$.  
The index~$1$ and $2$ cases are of particular interest: 
since $\#\Gamma_0(\mathbf{c})$ and $\ind(u ; \mathbf{c})$ always have the 
same parity due to the index formula, it follows from 
\eqref{eqn:cNindexAgain} that curves of 
index~$1$ or~$2$ satisfying our conditions have $c_N(u ; \mathbf{c}) = 0$ 
and thus $i(u ; \mathbf{c}\ |\ u ; \mathbf{c}) = 0$.
The transversality criterion $\ind(u ; \mathbf{c}) > c_N(u ; \mathbf{c})$
is clearly satisfied, and thus
$u$ lives in a $1$ or $2$--dimensional family of embedded curves that
never intersect each other.  These are precisely the curves that appear in
$J$--holomorphic foliations of $W$, or in the case where $W$ is a 
symplectization $\RR\times M$, the \emph{finite energy foliations} of
Hofer, Wysocki and Zehnder~\cite{HWZ:foliations}.  Isolated curves with 
$\ind(u ; \mathbf{c}) = 0$ can also occur in such foliations 
(surrounded by families
of larger index): we'll show for instance that stable, nicely embedded
index~$0$ curves appear as the underlying somewhere injective curves
when families of larger index degenerate to multiple covers.

It will be useful to note that due to \eqref{eqn:cNindexAgain}, all stable
nicely embedded curves also have the following properties:
\begin{enumerate}
\item $g = 0$,
\item $\#\Gamma_0(\mathbf{c}) = 1$ if $\ind(u ; \mathbf{c}) = 1$, 
and otherwise $\#\Gamma_0(\mathbf{c}) = 0$.
\end{enumerate}
In the cases $\ind(u ; \mathbf{c}) = 1$ or~$2$, we've observed that
$c_N(u ; \mathbf{c}) = 0$ and thus the adjunction formula 
also implies $\cov_\infty(\gamma_\Gamma ; \mathbf{c}) = 
\cov_{\MB}(\gamma_\Gamma ; \mathbf{c}) = 0$.  We will use this shortly
to prove the following consequence for the 
unique even puncture $z \in \Gamma_0(\mathbf{c})$ in the index~$1$ case:
\begin{prop}
\label{prop:uniqueEven}
If $u \in \mM^{\mathbf{c}}$ is a stable, nicely embedded curve with
$\ind(u ; \mathbf{c}) = 1$, then the unique even puncture 
$z \in \Gamma_0^\pm(\mathbf{c})$ satisfies one of the following:
\begin{enumerate}
\item $\gamma_z$ is nondegenerate and even,
\item $\gamma_z$ belongs to a $2$--dimensional Morse-Bott manifold,
and $\nu_\mp(\gamma_z) = 0$ if and only if $z \in \Gamma_C$.
\end{enumerate}
Moreover, $\gamma_z$ is either simply covered, or is doubly covered such
that the underlying simple orbit is nondegenerate and odd.
\end{prop}

\begin{defn}
\label{defn:badPuncture}
Adapting some terminology from Symplectic Field Theory \cite{SFT},
we will call $z \in \Gamma^\pm$ a \emph{bad puncture} if
$z \in \Gamma_0(\mathbf{c})$ and $\gamma_z = \gamma^2$ for some
nondegenerate odd orbit~$\gamma$.
\end{defn}
\begin{remark}
In this terminology, Prop.~\ref{prop:uniqueEven} says that the unique even
puncture has an orbit of covering number $1$ or $2$, and is bad in the latter
case.  In SFT of course, ``bad'' also means ``to be ignored'': moduli spaces
of curves with such punctures cannot be oriented, but they also need not be
counted in constructing the algebra of the theory.
\end{remark}

This is enough preparation to state the strong version of 
Theorem~\ref{thm:orbifold}.
In the following, we use expressions such as ``for generic $J$\ldots'' or
``$J$ is generic'' to mean more precisely:
``there exists a Baire subset
$\jJ \subset \jJ_\omega(W,\hH_+,\hH_-)$ such that the following is true
if $J \in \jJ$.''  Similarly, ``for generic homotopies\ldots'' means that
there exists a Baire subset in the space of smooth homotopies in
$\jJ_\omega(W,\hH_+,\hH_-)$ for which the statement is true.

\begin{thm}
\label{thm:orbifoldStrong}
Assume $\{ J_\tau \}_{\tau \in [0,1]}$ is a smooth $1$--parameter family
of almost complex structures in $\jJ_\omega(W,\hH_+,\hH_-)$
such that either
\begin{enumerate}
\item the homotopy $\tau\mapsto J_\tau$ is generic, or
\item $J_\tau = J$ is independent of $\tau$ and is generic.
\end{enumerate}
Suppose $\tau_n \to \tau_\infty \in [0,1]$
and $u_n : \dot{\Sigma} \to W$ is a sequence of asymptotically cylindrical
$J_{\tau_n}$--holomorphic curves, which are stable and nicely embedded with
respect to some fixed asymptotic constraints $\mathbf{c}$ and converge to
a smooth $J_{\tau_\infty}$--holomorphic curve $u : \dot{\Sigma} \to W$.  Then:
\begin{itemize}
\item
If $\ind(u ; \mathbf{c}) = 0$ or $\ind(u ; \mathbf{c}) = 1$ 
with $\gamma_z$ simply covered for the
unique even puncture $z \in \Gamma_0(\mathbf{c})$,
then $u$ is a stable, nicely embedded curve.
\item
If $\ind(u ; \mathbf{c}) = 1$ and the 
unique even puncture $z \in \Gamma_0(\mathbf{c})$ 
is bad (with $\gamma_z$ doubly covered), or $\ind(u ; \mathbf{c}) = 2$,
then $u$ is either a stable, nicely embedded curve or an unbranched
multiple cover of a stable, nicely embedded index~$0$ curve.
\end{itemize}
In all cases, $u$ is regular.
\end{thm}

Note that since $\sing(u ; \mathbf{c}) = 0$ for all stable, nicely embedded 
curves, the Morse-Bott contribution $\delta_{\MB}^\pm(\gamma_z)$ also 
vanishes at each puncture.  Plugging in \eqref{eqn:deltaMB} leads 
immediately to the following consequence:
\begin{lemma}
\label{lemma:deltaMBis0}
For any stable, nicely embedded curve $u$, if $z \in \Gamma_U^\pm$ is an
unconstrained puncture with a degenerate orbit $\gamma_z$ which is
exceptional in the sense of Prop.~\ref{prop:MBcharacterization}, then
$\nu_\mp(\gamma_z) = 0$ and $\cov_\mp(\gamma_z) = \cov_\mp(\gamma_z^\epsilon)$.
\end{lemma}

\begin{proof}[Proof of Prop.~\ref{prop:uniqueEven}]
The first alternative follows easily from
the formula $p(\gamma_z \pm \mathbf{c}_z) = \alpha^\Phi_+(\gamma_z \pm
\mathbf{c}_z) - \alpha^\Phi_-(\gamma_z \pm \mathbf{c}_z)$ and the definition
of~$\nu_\mp(\gamma_z)$.  We have also
$\cov_\mp(\gamma_z^\epsilon) = 1$ and 
$\left[\cov(\gamma_z^\epsilon) - 1 \right] \cdot \nu_\mp(\gamma_z) = 0$,
implying the same statements for $\gamma_z$ due to 
Lemma~\ref{lemma:deltaMBis0}.

If $\gamma_z$ is nondegenerate, we claim now that it cannot be a multiple
cover of any \emph{even} orbit $\gamma'$.  Otherwise there are eigenfunctions
$e_\pm$ of $\mathbf{A}_{\gamma'}$ with identical winding numbers and 
eigenvalues of opposite sign, so the corresponding covers give a pair of
neighboring eigenfunctions in the spectrum of $\mathbf{A}_{\gamma_z}$;
their eigenvalues are therefore the largest negative and smallest positive
elements of $\sigma(\mathbf{A}_{\gamma_z})$, implying $\cov_\mp(\gamma_z) > 1$,
a contradiction.  This leaves two possibilities 
for the simply covered orbit underlying $\gamma_z$: it is
either even (and thus is $\gamma_z$ itself) or is odd but hyperbolic, in
which case $\gamma_z$ can only be its double cover.

In the Morse-Bott case, suppose first that $\nu_\mp(\gamma_z) = 0$.
Then any section whose covering number is
counted by $\cov_\mp(\gamma_z^\epsilon)$ has the same winding and hence the
same covering number as a section in $\ker\mathbf{A}_{\gamma_z^\epsilon}$,
thus also the same covering number as $\gamma_z^\epsilon$ itself.  We conclude
that $\gamma_z^\epsilon$ is simply covered, so either $\gamma_z$ is as well or
it is an exceptional orbit with isotropy~$2$ as described in
Prop.~\ref{prop:MBcharacterization}.  The same result follows if
$\nu_\mp(\gamma_z) = 1$ because $\left[ \cov(\gamma_z) - 1 \right] \cdot
\nu_\mp(\gamma_z)$.
\end{proof}

In the proof of Theorem~\ref{thm:orbifoldStrong}, 
we'll need the following small variation on the
usual implicit function theorem:
\begin{lemma}
\label{lemma:IFT}
Suppose $f : X \to Y$ is a smooth Fredholm map between Banach spaces 
with $f(0) = 0$, and
$Q \subset X$ is a smooth finite dimensional submanifold of $X$ that contains
$0$, is contained in $f^{-1}(0)$ and satisfies
$$
\dim\ker df(0) = \dim Q.
$$
Then $Q$ also contains a neighborhood of $0$ in $f^{-1}(0)$; in particular
this neighborhood is a smooth manifold of dimension $\dim\ker df(0)$.
\end{lemma}
\begin{proof}
Let $V = \im df(0) \subset Y$ and choose a linear projection map
$\pi_V : Y \to V$ along some closed complement.  Then $\pi_V \circ f :
X \to V$ is also Fredholm and is regular at~$0$, so the implicit function
theorem gives $(\pi_V \circ f)^{-1}(0)$ near~$0$ the structure of a smooth
manifold of dimension~$\dim \ker df(0)$.  Now
$$
Q \subset f^{-1}(0) \subset (\pi_V \circ f)^{-1}(0),
$$
where the spaces on the left and right are manifolds of the same
dimension containing~$0$; the result follows.
\end{proof}

To every connected component of the moduli space $\mM^\mathbf{c}$, one can 
associate the data
$(\Sigma,\Gamma,P_\Gamma)$, where $\Sigma\setminus\Gamma$ is the domain of
any curve in the component (well defined up to diffeomorphism) 
and $P_\Gamma$ is the collection of orbits and/or Morse-Bott submanifolds
$\{ P_z \}_{z\in\Gamma}$ that determine the asymptotic behavior of such
a curve.  Let us introduce the notation
$$
\mM(\Sigma,\Gamma,P_\Gamma) \subset \mM^\mathbf{c}
$$ 
to indicate the union 
of all connected components of $\mM^\mathbf{c}$ that have this 
particular domain and asymptotic behavior.

\begin{lemma}
\label{lemma:coveringApriori}
For any component $\mM(\Sigma,\Gamma,P_\Gamma) \subset \mM^\mathbf{c}$,
there exists a finite set $\cC$ containing tuples
$(\Sigma',\Gamma',P_{\Gamma'},\mathbf{c}')$ such that the following is true: if
$u = v \circ \varphi \in \mM(\Sigma,\Gamma,P_\Gamma)$ is a multiple cover and
$v$ is the underlying somewhere injective curve, then there exists
$(\Sigma',\Gamma',P_{\Gamma'},\mathbf{c}') \in \cC$ such that
$v \in \mM(\Sigma',\Gamma',P_{\Gamma'}) \subset \mM^{\mathbf{c}'}$ and 
$\mathbf{c} \le \varphi^*\mathbf{c}'$ in the
sense of Def.~\ref{defn:partialOrder}.
\end{lemma}
\begin{proof}
The Riemann-Hurwitz formula constrains the genus of $\Sigma'$ to be less than
or equal to that of $\Sigma$, allowing only finitely many different
closed surfaces.  Having chosen $\Sigma'$, the relation $\gamma_z =
\gamma_{\varphi(z)}^{k_z}$ for $z \in \Gamma$ and $k_z := \ord(d\varphi;z) + 1$
allows $k_z$ to vary between $1$ and $\cov(\gamma_z)$, thus giving a finite
range of choices for each puncture.  After making this choice, we can
also decide which punctures $z,z' \in \Gamma$ might have the same image
under $\varphi$: this is allowed only when $P_z$ and $P_{z'}$ belong to the
same Morse-Bott manifold, and again presents a finite range of choices.
The number of punctures $\Gamma'$ and their asymptotic limits $P_{\Gamma'}$
are uniquely determined by this choice.  Finally, the constraints
$\mathbf{c}'$ can be defined as follows: for any constrained $z \in \Gamma$,
define $\zeta := \varphi(z)$ to be a constrained puncture, fixed at the
unique orbit $\gamma_\zeta$ such that $\gamma_z = \gamma_\zeta^{k_z}$.
Any puncture $\zeta \in \Gamma'$ not touched by this algorithm will be
considered unconstrained.  By construction now, $\mathbf{c} \le
\varphi^*\mathbf{c}'$.
\end{proof}

\begin{proof}[Proof of Theorem~\ref{thm:orbifoldStrong}]
We will carry out the proof in several steps
assuming $\{J_\tau\}_{\tau\in[0,1]}$ is a
generic homotopy; the proof for a fixed generic $J$ is the same but
slightly simpler in a few details.

If $u$ is somewhere injective there's nothing to prove, so assume
$u = v \circ \varphi$ for a somewhere injective curve
$v : \dot{\Sigma}' \to W$ and a holomorphic branched cover
$\varphi : \Sigma \to \Sigma'$ of degree $k \ge 2$.
By Lemma~\ref{lemma:coveringApriori}, the domain of $v$ is one out of a
finite set of choices and satisfies constraints $\mathbf{c}'$ with
$\mathbf{c} \le \varphi^*\mathbf{c}'$.  For each such choice, there exists
a generic set of homotopies $\{J_\tau\}$ such that we can assume
$\ind(v ; \mathbf{c}') \ge -1$, and the intersection of all these generic
sets is also generic, hence the genericity assumption implies
$\ind(v ; \mathbf{c}') \ge -1$.  By
\eqref{eqn:cNindexAgain} then, $c_N(v ; \mathbf{c}') \ge -1$.

\textbf{Step~1}: We show that $c_N(v ; \mathbf{c}') = -1$.  Combining
Prop.~\ref{prop:cNcover} and Prop.~\ref{prop:partialOrder} yields
$$
0 \ge c_N(u ; \mathbf{c}) \ge c_N(u ; \varphi^*\mathbf{c}') =
k c_N(v ; \mathbf{c}') + Z(d\dot{\varphi}) + Q,
$$
so the only other alternative is $c_N(v ; \mathbf{c}') = 0$, in which case
$c_N(u ; \mathbf{c}) = Z(d\dot{\varphi}) = Q = 0$.  Then all critical points of
$\varphi : \Sigma \to \Sigma'$ are at punctures, and the Riemann-Hurwitz
formula gives $2k - 2$ of them (counting multiplicity)
since both $\Sigma$ and $\Sigma'$ necessarily have genus zero.
Denote $k_z = \ord(d\varphi ; z) + 1 \in \NN$ 
for each $z \in \Gamma$.  Now $Q = 0$ implies that for each $z \in \Gamma^\pm$,
if $\zeta = \varphi(z)$,
$$
\alpha_\mp^\Phi(\gamma_z \pm (\varphi^*\mathbf{c}')_z) = 
\alpha_\mp^\Phi(\gamma_{\zeta}^{k_z} \pm (\varphi^*\mathbf{c}')_z) =
k_z \alpha_\mp^\Phi(\gamma_\zeta \pm \mathbf{c}'_\zeta),
$$
so depending on the constraints, we have either
$\alpha_\mp^\Phi(\gamma_z) \in k_z \ZZ$ or 
$\alpha_\mp^\Phi(\gamma_z \mp \boldsymbol{\delta}) \in k_z\ZZ$, the
latter only if $\zeta \in \Gamma'_U$, which implies $z \in \Gamma_U$.  
In either case, Lemma~\ref{lemma:eigenCover} then implies that a certain 
eigenfunction $e$ of $\mathbf{A}_{\gamma_z}$ is a $k_z$--fold cover.
If it's the first case, then $\cov_\mp(\gamma_z) \ge k_z$, and this
equals $\cov_\mp(\gamma_z^\epsilon)$ by Lemma~\ref{lemma:deltaMBis0}.
In the second case, we have $\cov(\gamma_z) \ge k_z$ and 
$e \in \ker\mathbf{A}_{\gamma_z}$.  If $\nu_\mp(\gamma_z) = 0$, then
$\wind^\Phi(e) = \alpha^\Phi_\mp(\gamma_z) \in k_z \ZZ$, so
Lemmas~\ref{lemma:eigenCover} and~\ref{lemma:deltaMBis0} again imply
$\cov_\mp(\gamma_z^\epsilon) = \cov_\mp(\gamma_z) \ge k_z$.
Otherwise $\nu_\mp(\gamma_z) = 1$, so Lemma~\ref{lemma:deltaMBis0} implies
that $\gamma_z$ is generic and thus 
$\nu_\mp(\gamma_z) \cdot \left[ \cov(\gamma_z^\epsilon) - 1 \right] =
\nu_\mp(\gamma_z) \cdot \left[ \cov(\gamma_z) - 1 \right] \ge k_z - 1$.
Putting all of these cases together and summing over $z \in \Gamma$, we find
$$
\cov_\infty(\gamma_\Gamma ; \mathbf{c}) + 
\cov_{\MB}(\gamma_\Gamma ; \mathbf{c})
\ge \sum_{z \in \Gamma} (k_z - 1) = 2 k - 2 \ge 2.
$$
Thus for large $n$, 
$i(u_n ; \mathbf{c}\ |\ u_n ; \mathbf{c}) = 2\sing(u_n ; \mathbf{c}) + 
c_N(u_n ; \mathbf{c}) + \cov_\infty(\gamma_\Gamma ; \mathbf{c}) +
\cov_{\MB}(\gamma_\Gamma ; \mathbf{c}) \ge 2$, a contradiction.

In light of this result and \eqref{eqn:cNindexAgain}, we have either
$\ind(v ; \mathbf{c}') = 0$ with all punctures odd or 
$\ind(v ; \mathbf{c}') = -1$ with exactly one even puncture.

\textbf{Step~2}: Claim $\ind(v ; \mathbf{c}') = 0$.  If not, 
then $\ind(v ; \mathbf{c}') = -1$ and
$\#\Gamma_0'(\mathbf{c}') = 1$, and since covers of even orbits are always even
(Lemma~\ref{lemma:coversAreEven}),
$\#\Gamma_0(\mathbf{c}) \ge 1$, implying $\ind(u ; \mathbf{c}) = 1$.  
In this case $u$ also
has exactly one even puncture $z \in \Gamma_0(\mathbf{c})$, so
$\gamma_z = \gamma_{\zeta}^k$ with $\zeta := \varphi(z)
\in \Gamma_0'(\mathbf{c}')$.
There are now three cases to consider:
\begin{enumerate}
\item
If $\gamma_{\zeta}$ is nondegenerate, then so is $\gamma_z$ and its
extremal eigenfunctions are the $k$--fold covers of those of
$\gamma_{\zeta}$, giving $\cov_\mp(\gamma_z^\epsilon) \ge k$.
\item
If $\gamma_{\zeta}$ is Morse-Bott with $\nu_\mp(\gamma_{\zeta}) = 0$,
then the extremal eigenfunction of a generic perturbation
$\gamma_{\zeta}^\epsilon$ has the same winding and thus same covering
number as a section in $\ker\mathbf{A}_{\gamma_{\zeta}^\epsilon}$,
and the same is true for the $k$--fold cover.  Moreover the Morse-Bott family
containing $\gamma_z^\epsilon$ is at least $k$--fold covered, which implies 
the same for sections in $\ker\mathbf{A}_{\gamma_z^\epsilon}$.  So again,
$\cov_\mp(\gamma_z^\epsilon) = \cov(\gamma_z^\epsilon) \ge k$.
\item
If $\nu_\mp(\gamma_{\zeta}) = 1$, then Lemma~\ref{lemma:deltaMBis0}
implies $\gamma_z$ is generic, so we can take $\gamma_z^\epsilon =
\gamma_z$ without loss of generality and conclude
$\left[ \cov(\gamma_z^\epsilon) - 1 \right] \cdot \nu_\mp(\gamma_z) \ge k$.
\end{enumerate}
The conclusion from all of these cases is that
$\cov_\infty(\gamma_\Gamma ; \mathbf{c}) + \cov_{\MB}(\gamma_\Gamma ; \mathbf{c})
\ge k - 1 \ge 1$, and since $c_N(u ; \mathbf{c}) = 0$
for the index~$1$ case, a contradiction arises again from the
adjunction formula: $i(u_n ; \mathbf{c}\ |\ u_n ; \mathbf{c}) = 
2\sing(u_n ; \mathbf{c}) + c_N(u_n ; \mathbf{c}) + 
\cov_\infty(\gamma_\Gamma ; \mathbf{c}) + \cov_{\MB}(\gamma_\Gamma ; \mathbf{c}) 
\ge 1$.

\textbf{Step~3}: Since $0 = \ind(v ; \mathbf{c}') > c_N(v ; \mathbf{c}') = -1$, 
it now follows immediately from Prop.~\ref{prop:iCovering} and
Prop.~\ref{prop:partialOrder} that $v$ is a stable, nicely embedded curve, as
$$
0 \ge i(u ; \mathbf{c} \ |\ u ; \mathbf{c}) \ge
i(u ; \varphi^*\mathbf{c}' \ |\ u ; \varphi^*\mathbf{c}')
\ge k^2 \cdot i(v ; \mathbf{c}'\ |\ v ; \mathbf{c}').
$$

\textbf{Step~4}: If $\ind(u ; \mathbf{c}) = 1$, then its unique even orbit 
cannot be simply covered since $v$ has only odd orbits.  Thus the even orbit 
must be a doubly covered orbit at a bad puncture.

\textbf{Step~5}: We claim $u$ is immersed and has $\ind(u ; \mathbf{c}) > 0$.  
Suppose not, i.e.~that either $Z(du) > 0$ or $\ind(u ; \mathbf{c}) = 0$.  
Then $\ind(u ; \mathbf{c}) \le
2Z(du)$, so Theorem~\ref{thm:criterion} gives
\begin{equation}
\label{eqn:crucialKernel}
\begin{split}
2Z(du) &\le \dim\ker\left( D\dbar_{J_{\tau_\infty}}(j,u) / \aut(\dot{\Sigma},j) \right)\\
&\le 2Z(du) + K(c_N(u ; \mathbf{c}) - Z(du),
\#\Gamma_0(\mathbf{c})) = 2Z(du)
\end{split}
\end{equation}
since $c_N(u ; \mathbf{c}) - Z(du) < 0$.  Extending the usual bundle on
$\tT \times \bB$ to allow parametrized $J$, we can now consider a nonlinear
operator
$$
\dbar : [0,1] \times \tT \times \bB \to \eE : (\tau,j',u') \mapsto
\dbar_{J_\tau}(j',u').
$$
Since $v$ is embedded, every critical point of $u$ arises
as a branch point of $\dot{\varphi} : \dot{\Sigma} \to \dot{\Sigma}'$, thus
$Z(du) = Z(d\dot{\varphi})$, and \eqref{eqn:crucialKernel} now implies
\begin{equation}
\label{eqn:dimPlus1}
\dim \ker D\dbar(\tau_\infty,j,u) \le 2Z(d\dot{\varphi}) 
+ \dim\Aut(\dot{\Sigma},j) + 1.
\end{equation}
To apply Lemma~\ref{lemma:IFT}, we shall now find a smooth manifold of
precisely this dimension that is contained in $\dbar^{-1}(0)$.  The key
is to look at the space of holomorphic branched covers 
$\dot{\Sigma} \to \dot{\Sigma}'$ close to $\dot{\varphi}$.
Observe that since $v$ is embedded and satisfies the transversality criterion
$0 = \ind(v ; \mathbf{c}') > c_N(v ; \mathbf{c}') = -1$, for $\tau$ close
to $\tau_\infty$ we obtain from the implicit function theorem a
smooth $1$--parameter family of asymptotically cylindrical 
pseudoholomorphic maps
$$
v_\tau : (\dot{\Sigma}',j_\tau) \to (W,J_\tau)
$$
satisfying the constraints $\mathbf{c}'$, with $v_{\tau_\infty} = v$.  
The holomorphic maps from $\dot{\Sigma}$ to $(\dot{\Sigma}',j_\tau)$ can
then be identified with the zero-set of a section
$$
\dbar_{\Sigma'} : [0,1] \times \tT \times \bB^{\Sigma'} \to \eE^{\Sigma'} :
(\tau,j',\varphi') \mapsto T\varphi' + j_\tau \circ T\varphi' \circ j',
$$
and by the remarks in Example~\ref{ex:branchedCovers}, a neighborhood of
$\dot{\varphi}$ in $\dbar_{\Sigma'}^{-1}(0)$ is a smooth manifold
of dimension $2Z(\dot{\varphi}) + \Aut(\dot{\Sigma},j) + 1$.
Now for $(\tau,j',\varphi') \in \dbar_{\Sigma'}^{-1}(0)$, 
we have $(\tau,j',v_\tau \circ \varphi') \in \dbar^{-1}(0)$,
thus embedding $\dbar_{\Sigma'}^{-1}(0)$ as a smooth submanifold
of dimension $2Z(d\dot{\varphi}) +
\dim\Aut(\dot{\Sigma},j) + 1$ in $\dbar^{-1}(0)$.  It follows that
\eqref{eqn:dimPlus1} is an equality, and
Lemma~\ref{lemma:IFT} now implies that every element of
$\dbar^{-1}(0)$ near $(\tau_\infty,j,u)$ belongs to this submanifold;
this is a contradiction, as it implies that for large $n$,
$u_n$ must also be a multiple cover.

\textbf{Step~6}: Since we now know that $u$ is immersed, it is
immediate from Theorem~\ref{thm:criterion} that $u$ is regular.
\end{proof}

\begin{cor}
\label{cor:orbifold}
For generic $J$, if $u \in \mM^\mathbf{c}$ is a stable, nicely embedded 
$J$--holomorphic curve then $\mM_u^\mathbf{c}$ is a smooth orbifold of 
dimension $\ind(u ; \mathbf{c})$ with only isolated singularities.
Moreover, the images of any two curves in $\mM_u^\mathbf{c}$ are either
identical or disjoint, and they are all stable and
nicely embedded except for a discrete set of unbranched multiple covers
of embedded curves.
\end{cor}
\begin{proof}
The statement about the images follows from positivity of intersections and
the condition $i(u ; \mathbf{c}\ |\ u ; \mathbf{c}) \le 0$.
The only remaining part not immediate from 
Theorem~\ref{thm:orbifoldStrong} is that
the multiple covers are isolated; this is related to the fact that the
orbifold singularities must be isolated for orientation reasons
(see Remark~\ref{remark:orientations} below), but doesn't follow from it.
So, we claim that for any
multiple cover $u = v \circ \varphi$ arising as a limit in the theorem,
every other curve close to $u$ is somewhere injective.  Here we can
assume $\ind(u ; \mathbf{c})$ is $1$ or $2$, so $c_N(u ; \mathbf{c}) = 
\cov_\infty(\gamma_\Gamma ; \mathbf{c}) = 
\cov_{\MB}(\gamma_\Gamma ; \mathbf{c}) = 0$.  Now note that 
Theorem~\ref{thm:normalCR} and the implicit function 
theorem give a natural isomorphism $T_u \mM^\mathbf{c} = 
\ker \mathbf{D}_u^N$ since
$u$ is immersed, and it will thus suffice to show that nontrivial
sections $\eta \in \ker \mathbf{D}_u^N \subset \Gamma(N_u)$ are not
multiply covered.  Otherwise, using the natural identification
$N_u = \varphi^*N_v$, there exists a nonzero $\eta \in \ker \mathbf{D}_u^N$
and a section $\eta' \in \Gamma(N_v)$
such that $\eta = \eta' \circ \varphi$.  We know that $\eta$ is zero
free (also at infinity), since using \eqref{eqn:c1counts},
Prop.~\ref{prop:c1cN} and the usual identification of $\mathbf{D}_u^N$
with an operator on non-weighted spaces,
$$
Z(\eta) + Z_\infty(\eta) = c_1(N_u,\ell_N,\mathbf{A}_\Gamma \pm 
\mathbf{c}_\Gamma) = c_N(u ; \mathbf{c}) = 0.
$$
This implies that the winding of $\eta$ near each puncture $z \in \Gamma^\pm$
attains the extremal value $\alpha_\mp^\Phi(\gamma_z \pm \mathbf{c}_z)$.  
But this is impossible if $\eta = \eta' \circ \varphi$: indeed, the fact that
$\dot{\varphi} : \dot{\Sigma} \to \dot{\Sigma}'$ is immersed but both
surfaces have genus zero implies that there exists a puncture $z \in \Gamma$
at which $\varphi$ has nontrivial branching order 
$k_z := \ord(d\varphi ; z) > 1$,
so the asymptotic winding $\wind_z^\Phi(\eta)$ of $\eta$ near $z$ satisfies
\begin{equation*}
\begin{split}
\pm \wind^\Phi_z(\eta) &\le
\pm k_z \alpha^\Phi_\mp(\gamma_{\varphi(z)} \pm \mathbf{c}'_{\varphi(z)}) \\
&= \pm \alpha_\mp^\Phi(\gamma_z \pm (\varphi^*\mathbf{c}')_{\varphi(z)}) -
q_\mp(\gamma_{\varphi(z)} \pm \mathbf{c}'_{\varphi(z)} ; k_z) \\
&\le \pm\alpha^\Phi_\mp(\gamma_z \pm \mathbf{c}_z) - 
q_\mp(\gamma_{\varphi(z)} \pm \mathbf{c}'_{\varphi(z)} ; k_z),
\end{split}
\end{equation*}
implying $q_\mp(\gamma_{\varphi(z)} \pm \mathbf{c}'_{\varphi(z)} ; k_z) = 0$ 
and thus $\alpha^\Phi_\mp(\gamma_z \pm \mathbf{c}_z) \in k_z\ZZ$.
Then repeating an argument that is by now familiar from the proof of
Theorem~\ref{thm:orbifoldStrong}, we find a contradiction in the form
$\cov_\infty(\gamma_\Gamma ; \mathbf{c}) + 
\cov_{\MB}(\gamma_\Gamma ; \mathbf{c}) > 0$.
\end{proof}

\begin{remark}
\label{remark:foliation}
It is shown in \cite{Wendl:BP1} that if $u : \dot{\Sigma} \to W$ is a 
stable nicely embedded index~$2$ curve, then the nearby curves in
$\mM^\mathbf{c}_u$ foliate a neighborhood of $u(\dot{\Sigma})$ in~$W$.
Now suppose $u$ is a multiply covered index~$2$ curve that is a limit of
stable nicely embedded curves.  Then since $u$ is immersed and 
$T_u\mM^\mathbf{c}$ consists of zero-free sections of its normal bundle, the
same argument shows that the
nearby curves in $\mM_u^\mathbf{c}$ again foliate a neighborhood of
$u(\dot{\Sigma})$.  In this foliation, $u(\dot{\Sigma})$ is
an exceptional leaf, being the embedded image of an isolated index~$0$
curve.  An explicit example is constructed below.
\end{remark}

\begin{remark}
\label{remark:orientations}
The fact that singularities in a $1$--dimensional orbifold are isolated is
obvious, and in two dimensions it's true if the orbifold is oriented, as
an oriented orbifold can only have singularities of codimension at least two.
By results in \cite{BourgeoisMohnke}, $\mM_u$ does admit an orientation if 
$u$ is a stable, nicely embedded curve of index~$2$,
and the same is true for index~$1$ if and only if the
unique even puncture is not a ``bad'' puncture.  This excludes
singularities in the index~$1$ case entirely unless the even puncture is bad,
and indeed, we've shown that multiple covers don't appear
in this case.  These remarks are not quite enough to prove
Cor.~\ref{cor:orbifold} however, as in general there can
be multiple covers with trivial automorphism groups, which therefore do not
cause singularities.
\end{remark}

\begin{example}
\label{ex:foliation}
We now consider a concrete situation in which nicely embedded curves of
index~$2$ are seen to converge to an isolated, unbranched multiple cover.

Identify $S^2$ with the extended complex plane and
let $W = (S^2 \times S^2) \setminus \{ (0,0),(\infty,\infty),(1,1) \}$,
choosing the standard complex structure $J = i \oplus i$.  This can be
regarded as a manifold with three negative cylindrical ends asymptotic
to the standard contact $3$--sphere, whose Reeb orbits are the fibers of
the Hopf fibration.  The asymptotics are therefore Morse-Bott: there is
a $2$--dimensional family of closed orbits at each end.  Now for
$\zeta \in \CC \setminus \{ 0,-1,1 \}$, consider the $2$--dimensional family
of $J$--holomorphic four-punctured spheres
$$
u_\zeta : S^2 \setminus \{0,1,-1,\infty\} \to W :
z \mapsto \left( z^3 \frac{z + \zeta}{\zeta z + 1}, z^2 \right).
$$
These are all proper and embedded, with asymptotic behavior as follows:
\begin{itemize}
\item
At $0$, $u_\zeta$ is asymptotic to a fixed doubly covered orbit in the end
$(0,0)$.
\item
At $\infty$, $u_\zeta$ is asymptotic to a fixed doubly covered orbit in the end
$(\infty,\infty)$.
\item
At $1$ and $-1$, $u_\zeta$ is asymptotic to an arbitrary (not fixed) simply
covered orbit in the end $(1,1)$.
\end{itemize}
One can use the setup we've described to show that the
moduli space of embedded holomorphic curves satisfying precisely these
asymptotic constraints and representing the same relative homology class
is indeed a smooth $2$--dimensional manifold:
indeed, $\ind(u ; \mathbf{c}) = 2$ and
$i(u ; \mathbf{c}\ |\ u ; \mathbf{c}) = c_N(u ; \mathbf{c}) = 
\cov_\infty(\gamma_\Gamma ; \mathbf{c}) = 
\cov_{\MB}(\gamma_\Gamma ; \mathbf{c}) = 0$.
Now as $\zeta \to 0$, the family converges
to the curve $u_0(z) = (z^4,z^2) =: v(z^2)$,
an unbranched double cover of the embedded $3$--punctured sphere
$$
v : S^2 \setminus \{0,1,\infty\} \to W : z \mapsto (z^2,z),
$$
and with the appropriate asymptotic constraints $\mathbf{c}'$ 
one can indeed show that $\ind(v ; \mathbf{c}') = 0$ and 
$i(v ; \mathbf{c}' \ |\ v ; \mathbf{c}') = c_N(v ; \mathbf{c}') = -1$.
Observe that the images of $u_\zeta$ for $\zeta$ 
near zero together with the image of $v$ foliate a neighborhood of $v$ in~$W$.
Due to the ordering of the punctures, $\Aut(u_0)$ is the trivial group,
so the moduli space remains a smooth manifold even with $u_0$ included.
If we take the quotient of this space by forgetting the order of the punctures,
it becomes a smooth orbifold in which $u_0$ has isotropy group~$\ZZ_2$.
\end{example}

\begin{remark}
\label{remark:bigCompactness}
It's also interesting to see what happens to the family $u_\zeta$ as
$\zeta \to \pm 1$ or $\zeta \to \infty$: here it turns out that $u_\zeta$ 
breaks into a $J$--holomorphic building (in the sense of the
SFT compactness theorem \cite{SFTcompactness}).  At $\zeta \to \infty$ in
particular, the building includes a component that is an unbranched
multiple cover of index~$-2$ over a nicely embedded index~$0$ curve.
It is work in progress by the author to generalize 
Theorem~\ref{thm:orbifoldStrong} in light of SFT compactness and
show that such behavior is
quite general: indeed, that only \emph{unbranched} multiple covers can arise in
such limits, and that there exists a well behaved gluing theory for
buildings of this type.
\end{remark}

\appendix

\section{Counting boundary zeros}
\label{app:zeros}

In this appendix we define a $\frac{1}{2}\ZZ$--valued count of zeroes for 
sections of a complex line bundle with totally real boundary condition.
Let $(E,J)\to S$ be a topological complex line bundle over a 
compact, connected and oriented surface with boundary.  
Partition the boundary into disjoint subsets
$\p S = \p_0 S \sqcup \p_1 S$, either of which may be empty.  
Now choose a totally real subbundle
$\ell \subset E|_{\p_0 S} \to \p_0 S$, and consider the space of
all continuous sections $\sigma : S \to E$ such that
$\sigma(\p_0 S) \subset \ell$ and $\sigma \ne 0$ on $\p_1 S$.  
We will call such sections
\emph{admissible}.  Suppose $\sigma$ is an admissible section with
a discrete zero set $\sigma^{-1}(0) \subset S$.  If 
$z_0\in \sigma^{-1}(0) \cap \interior{S}$, then it
is standard to define the \emph{order} of the zero $\ord(\sigma;z_0)$ as the
winding number of $\sigma$ over a small loop around $z_0$, computed
in any local trivialization.  The boundary condition makes it possible
to extend this definition to isolated zeros on $\p_0 S$ as well:
for $z_0\in Z(\sigma) \cap \p_0 S$, 
choose coordinates identifying a neighborhood
$\uU$ of $z_0$ with $\DD^+ = 
\{ z\in \CC \ |\ |z| \le 1 \text{ and } \Im{z} \ge 0\}$,
such that $z_0 = 0$ and $\uU \cap \p S = \DD^+ \cap \RR$.  Choose also
a local trivialization over $\uU$
that identifies $\ell$ 
with $(\DD^+ \cap \RR) \times \RR \subset \DD^+ \times \CC$.
Then $\sigma$ is represented on this neighborhood by a continuous function
$f : \DD^+ \to \CC$, satisfying the boundary condition
$f(\DD^+ \cap \RR) \subset \RR$.  We can therefore extend $f$ to a continuous
function $f^D : \DD\to \CC$ on the full disk, satisfying
$f^D(\bar{z}) = \overline{f^D(z)}$.  The order
$\ord(\sigma;z_0)$ 
is then the order of the isolated zero of $f^D$ at $0$, i.e.~the
winding number of $f^D$ for a small circle about~$0$.
This definition doesn't depend on the choices.

For an admissible section $\sigma$ with discrete zero set $\sigma^{-1}(0)$, 
we now define the algebraic count of zeros by
$$
Z(\sigma) = \sum_{z\in \sigma^{-1}(0) \cap \interior{S}} \ord(\sigma ; z)
 + \frac{1}{2} \sum_{z\in \sigma^{-1}(0) \cap \p_0 S} \ord(\sigma ; z).
$$

\begin{prop}
\label{prop:zeroCounting}
Suppose $\sigma_0$ and $\sigma_1$ are admissible sections with isolated
zeros, and are homotopic through a family of admissible sections.
Then $Z(\sigma_0) = Z(\sigma_1)$.
\end{prop}
\begin{proof}
This is clear if $\p_0 S = \emptyset$: then $Z(\sigma)$ is the Euler number
of $E$ if
$\p_1 S = \emptyset$, or more generally the homotopy invariant winding
number about $\p_1 S$ with respect to any global trivialization.

We reduce the general case to this by a doubling argument: define the
conjugate surface $S^C := S$ with the opposite orientation, and the
conjugate bundle $(E^C,J^C) := (E,-J) \to S^C$.  Then we can glue
$S$ to $S^C$ along $\p_0 S$ to define the doubled surface $S^D$, and
similarly form a bundle $(E^D,J^D) \to S^D$ by gluing $(E,J)$ to
$(E^C,J^C)$ via the unique complex bundle isomorphism
$E|_{\p_0 S} \to E^C|_{\p_0 S}$ that restricts to the identity on~$\ell$.
Now $\p S^D = \p_1 S^D = \p_1 S \cup \p_1 S^C$.  Any admissible
section $\sigma$ of $E$ defines an admissible section $\sigma^D$
of $E^D$, and the
same statement applies to homotopies, thus it suffices to prove the
following formula relating $Z(\sigma)$ to $Z(\sigma^D)$:
\begin{equation*}
\begin{split}
Z(\sigma^D) &= \sum_{z\in (\sigma^D)^{-1}(0)} \ord(\sigma^D ; z) \\
 &= 2\sum_{z\in \sigma^{-1}(0) \cap \interior{S}} \ord(\sigma; z) + 
 \sum_{z\in \sigma^{-1}(0) \cap \p_0 S} \ord(\sigma ; z) \\
&= 2 Z(\sigma).
\end{split}
\end{equation*}
This follows from two important facts which are easy to check: first,
if $z$ is a zero of $\sigma$ in $\interior{S}$, its order is the same as that
of the corresponding zero in $S^C$.
Secondly, if $z$ is a boundary zero of $\sigma$,
then its order equals its order as an interior zero of $\sigma^D$.
\end{proof}

The doubling formula $Z(\sigma^D) = 2 Z(\sigma)$ which emerged from this
proof is a useful fact in itself; we shall apply it now to express
$Z(\sigma)$ in terms of the relative first Chern number of $E$ and the
boundary Maslov index of the pair $(E|_{\p_0 S},\ell)$.

\begin{prop}
\label{prop:halfc1}
For any choice of trivialization $\Phi$ along $\p S$,
$$
Z(\sigma) = c_1^\Phi(E) + \frac{1}{2}\mu^\Phi(E,\ell) +
\wind_{\p_1 S}^\Phi(\sigma).
$$
\end{prop}
\begin{proof}
Label the right hand side $\widehat{Z}(\sigma)$ and observe that it does
not depend on $\Phi$ and depends on $\sigma$ only up to homotopy through
admissible sections.  Moreover it is clear that 
$Z(\sigma) = \widehat{Z}(\sigma)$ if $\p_0 S = \emptyset$, so it will
suffice to prove the doubling formula $\widehat{Z}(\sigma^D) =
2\widehat{Z}(\sigma)$.  Since the orientations of both $E^C$ and
$S^C$ are reversed, we have $\wind_{\p_1 S^D}^{\Phi^D}(\sigma^D) =
2\wind_{\p_1 S}^\Phi(\sigma)$ for the natural trivialization $\Phi^D$ induced
by~$\Phi$.  We claim also that
$$
c_1^{\Phi^D}(E^D) = 2c_1^\Phi(E) + \mu^\Phi(E,\ell),
$$
which will prove the result.  This can be reduced to the standard additivity
of the Maslov index under gluing.  Construct a new surface $\bar{S}
\supset S$ by gluing a disk to each component of $\p_1 S$, and glue in 
trivial bundles along $\Phi$ over these disks to produce a new bundle
$(\bar{E},\bar{J}) \to \bar{S}$, such that
$\bar{E}|_{S} = E$ and $c_1^\Phi(\bar{E}) = c_1^\Phi(E)$.
Now $\p\bar{S} = \p_0 S$, and $2 c_1^\Phi(E) + \mu^\Phi(E,\ell)$ is by
definition the \emph{absolute} Maslov index $\mu(\bar{E},\ell)$.
(Alternatively, one can define the latter as $\mu^\Phi(\bar{E},\ell)$
where $\Phi$ is any trivialization along $\p\bar{S}$ that extends
\emph{globally} over~$\bar{S}$.)  Now the gluing property for
$\mu(\bar{E},\ell)$ gives
$$
2\mu(\bar{E},\ell) = \mu\left(\bar{E}^D\right) 
= 2 c_1\left(\bar{E}^D\right)
$$
since $\bar{S}^D$ is closed.  But the latter is also equal to
$2 c_1^{\Phi^D}\left(E^D\right)$, proving the claim.
\end{proof}

\end{document}